\documentclass[11pt]{article}
\usepackage[letterpaper,hmargin=26mm,vmargin=1in]{geometry}
\usepackage{graphicx,amssymb,amsmath,amsthm}

\usepackage{latexsym}
\usepackage{epic}
\usepackage{eepic}

\newtheorem{theorem}{Theorem}[section]
\newtheorem{lemma}[theorem]{Lemma}

\newtheorem{corollary}[theorem]{Corollary}
\newtheorem{conjecture}[theorem]{Conjecture}
\theoremstyle{definition}
\newtheorem{definition}[theorem]{Definition}
\newtheorem{definitions}[theorem]{Definitions}

\newcommand{\ignore}[1]{}

\newcommand{\D}{D}

\begin{document}

\title{Overhang \thanks{A preliminary version of this paper \cite{PZ06} appeared in the
Proceedings of the 17th Annual ACM-SIAM Symposium on Discrete
Algorithms (SODA'06), pages 231--240. This full version is to
appear in the American Mathematical Monthly.}}

\author{
\em Mike Paterson
\thanks{DIMAP and Department of Computer Science, University of Warwick,
Coventry CV4 7AL, UK. E-mail: {\tt msp@dcs.warwick.ac.uk}}
\and \em Uri Zwick
\thanks{School of Computer Science, Tel Aviv
University, Tel Aviv 69978, Israel. E-mail: {\tt zwick@cs.tau.ac.il}}
}

\date{} 

\maketitle

\begin{figure}[bh]
\begin{center}
\centerline{
\includegraphics[height=70mm]{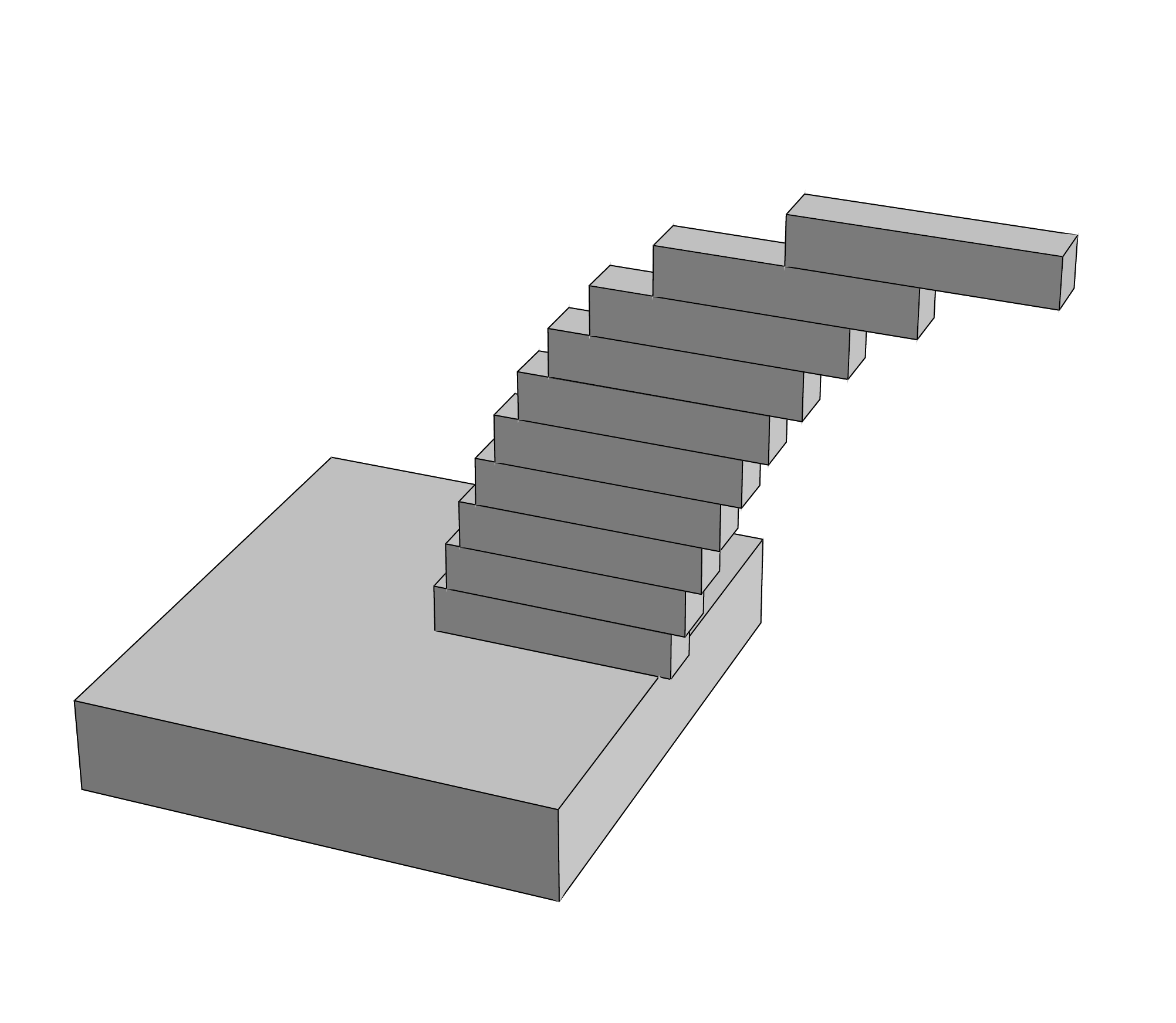}\hspace*{0.4cm}
\raisebox{0.7cm}{\includegraphics[height=50mm]{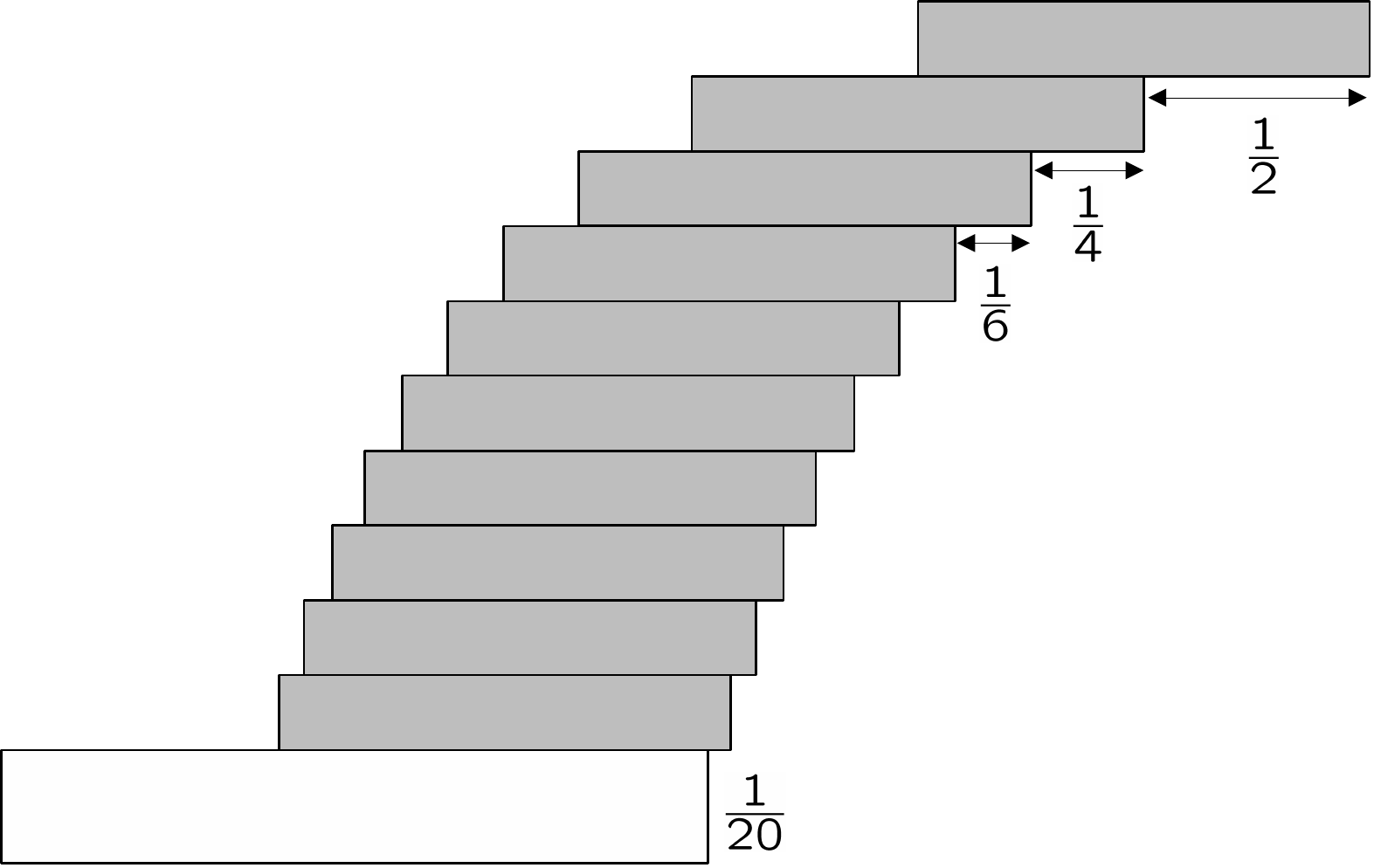}}}
\caption{A harmonic stack with 10 blocks.} \label{fig:h10}
\end{center}\vspace{-5mm}
\end{figure}

\section{Introduction} \label{sec:intro}

How far off the edge of the table can we reach by stacking $n$
identical, homogeneous, frictionless blocks of length~1? A classical
solution achieves an overhang asymptotic to $\frac{1}{2} \ln n$.
This solution is widely believed to be optimal. We show, however,
that it is exponentially far from optimality by constructing simple
$n$-block stacks that achieve an overhang of~$cn^{1/3}$, for some
constant $c>0$.

The problem of stacking a set of objects, such as bricks,
books, or cards, on a tabletop to maximize the
overhang is an attractive problem with
a long history. J. G. Coffin~\cite{C23} posed the problem in the
``Problems and Solutions" section of this
{\sc Monthly}, but no solution was given there. The problem recurred from
time to time over subsequent
years, e.g., \cite{S53,S54}, \cite{J55}, \cite{E59}. Either
deliberately or inadvertently, these authors all seem to have
introduced the further restriction that there can be at most one
object resting on top of another. Under this restriction, the
\emph{harmonic stacks}, described below, are easily seen to be
optimal.

The classical harmonic stack of size~$n$ is composed of $n$ blocks
stacked one on top of the other, with the $i$th block from the top
extending by $\frac{1}{2i}$ beyond the block below it. (We assume
that the length of each block is~$1$.) The overhang achieved by the
construction is clearly $\frac{1}{2} H_n$, where
$H_n=\sum_{i=1}^n\frac{1}{i}\sim \ln n$ is the $n$th harmonic
number. Both a 3D and a 2D view of the harmonic stack of size~10 are
given in Figure~\ref{fig:h10}. The harmonic stack of size~$n$ is
balanced since, for every $i< n$, the center of mass of the topmost
$i$ blocks lies exactly above the right-hand edge of the $(i+1)$st
block, as can be easily verified by induction. Similarly, the center
of mass of all the $n$ blocks lies exactly above the right edge of
the table. A formal definition of ``balanced'' is given in
Definition~\ref{def:balance}. A perhaps surprising and
counterintuitive consequence of the harmonic stacks construction is
that, given sufficiently many blocks, it is possible to obtain an
arbitrarily large overhang!

Harmonic stacks became widely known in the recreational math
community as a result of their appearance in the \emph{Puzzle-Math}
book of Gamow and Stern \cite{GS58} (Building-Blocks, pp.~90--93)
and in Martin Gardner's ``Mathematical Games'' section of the
November 1964 issue of Scientific American~\cite{G64} (see
also~\cite{G71}, Chapter~17: Limits of Infinite Series, p.~167).
Gardner refers to the fact that an arbitrarily large overhang can be
achieved, using sufficiently many blocks, as the
\emph{infinite-offset paradox}.
Harmonic stacks were subsequently used by countless authors as an
introduction to recurrence relations, the harmonic series, and
simple optimization problems; see, e.g.,~\cite{GKP88} pp.~258--260.
Hall~\cite{H05} notes that harmonic stacks started to appear in
textbooks on physics and engineering mechanics as early as the
mid-19th century (see, e.g.,~\cite{M07} p.~341, \cite{P50}
pp.~140--141, \cite{W55} p.~183).

It is perhaps surprising that none of the sources cited above
realizes how limiting is the one-on-one restriction under which the
harmonic stacks are optimal. Without this restriction, blocks can
be used as counterweights, to balance other blocks. The problem then
becomes vastly more interesting, and an exponentially larger
overhang can be obtained.

\begin{figure}[t]
\begin{center}
\includegraphics[height=50mm]{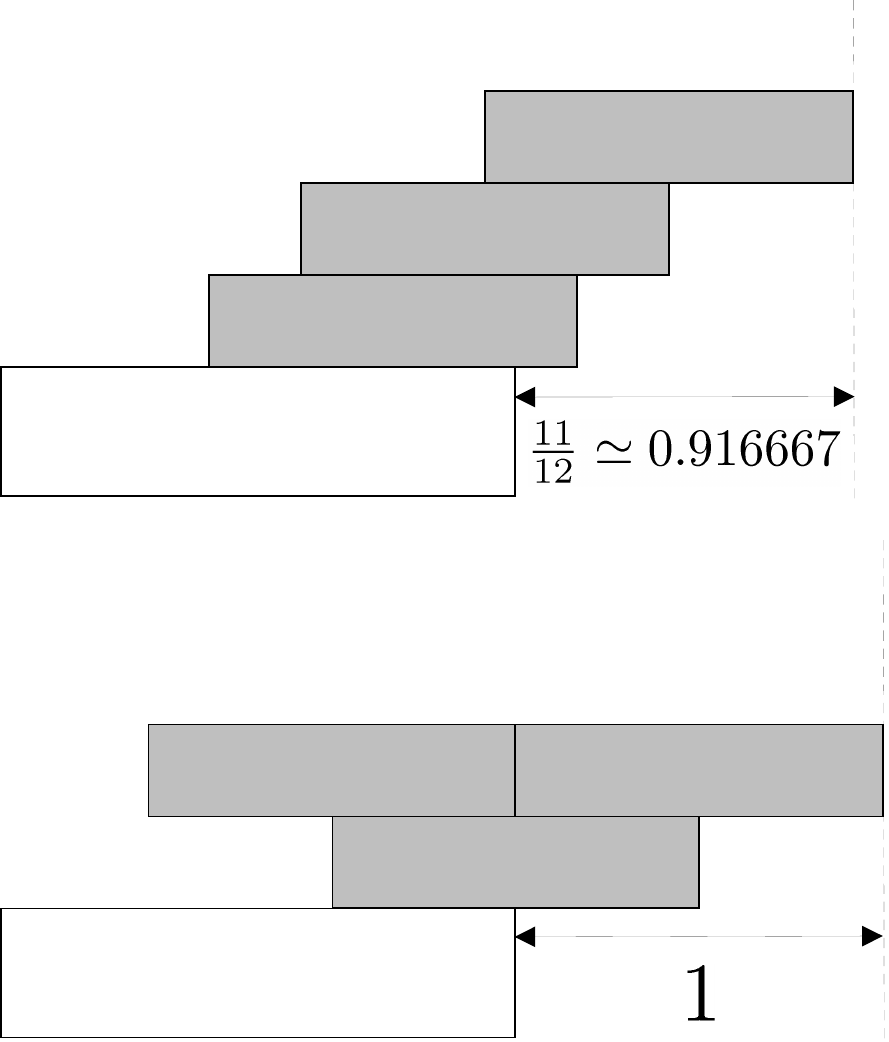}\hspace*{1cm}
\includegraphics[height=50mm]{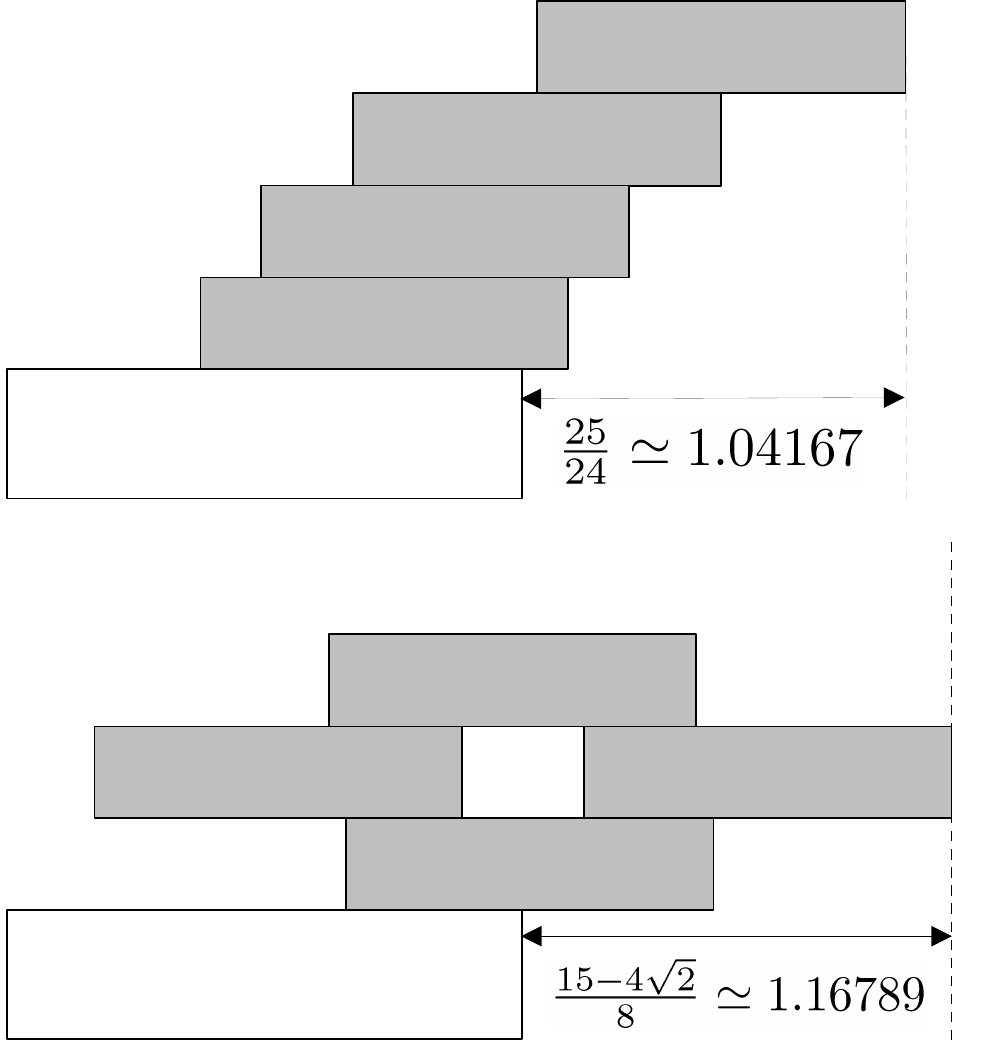}
\caption{Optimal stacks with 3 and 4 blocks compared with the
corresponding harmonic stacks.} \label{fig:opt34}
\end{center}\vspace{-5mm}
\end{figure}

Stacks with a specific small number of blocks that do not satisfy
the one-on-one restriction were considered before by several other
authors. Sutton \cite{S55}, for example, considered the case of
three blocks.
One of us
set a stacking problem with three uniform thin planks of lengths~2,~3,
and~4 for the Archimedeans Problems Drive in 1964 \cite{HP64}.
Ainley~\cite{A79} found the maximum overhang achievable with four
blocks to be $\frac{15-4\sqrt 2}{8}\sim 1.16789$.
The optimal stacks with~3 and~4 blocks are shown, together with the
corresponding harmonic stacks, in Figure~\ref{fig:opt34}.

Very recently, and independently of our work, Hall \cite{H05}
explicitly raises the problem of finding stacks of blocks that
maximize the overhang without the one-on-one restriction. (Hall
calls such stacks \emph{multiwide stacks}.) Hall gives a sequence of
stacks which he claims, without proof, to be optimal. We show,
however, that the stacks suggested by him are optimal only for $n\leqslant
19$. The stacks claimed by Hall to be optimal fall into a natural
class that we call \emph{spinal stacks}. We show in
Section~\ref{sec:spinal} that the maximum overhang achievable using
such stacks is only $\ln n+O(1)$. Thus, although spinal stacks
achieve, asymptotically, an overhang which is roughly twice the
overhang achieved by harmonic stacks, they are still exponentially
far from being optimal.

\begin{figure*}[t]

\centerline{
\includegraphics[height=3cm]{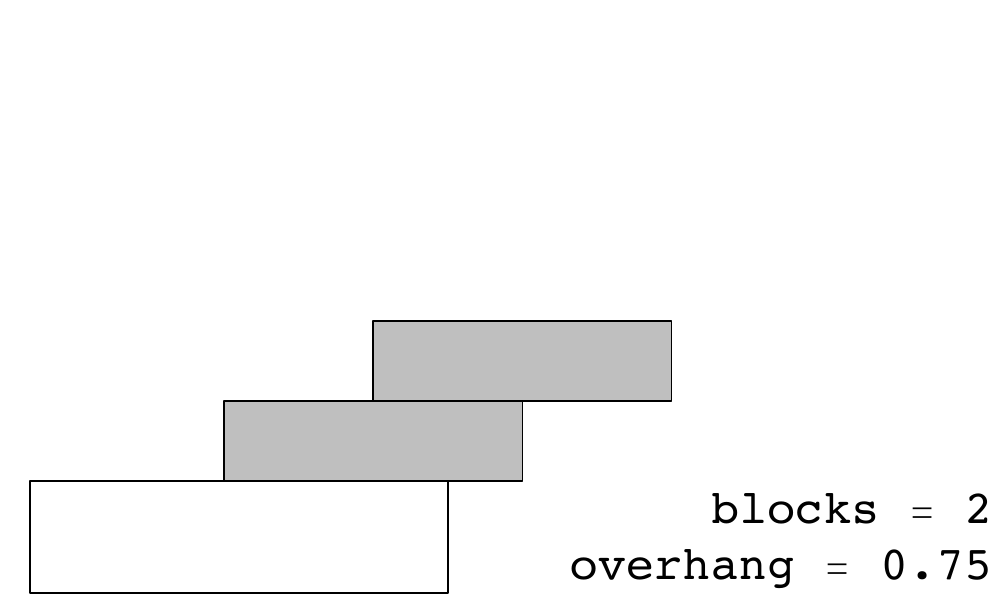}\hspace*{5mm}
\includegraphics[height=3cm]{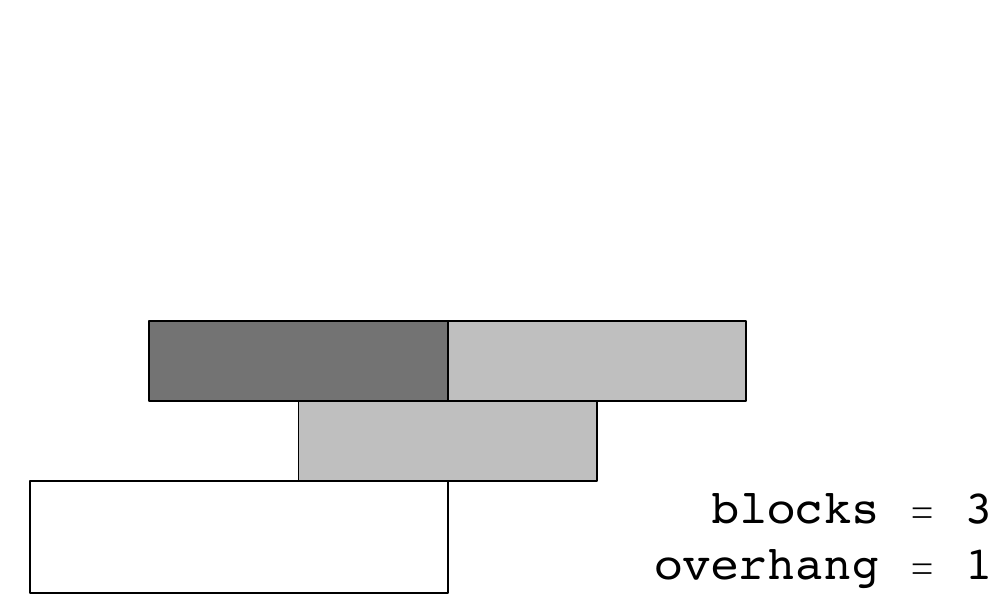}\hspace*{5mm}
\includegraphics[height=3cm]{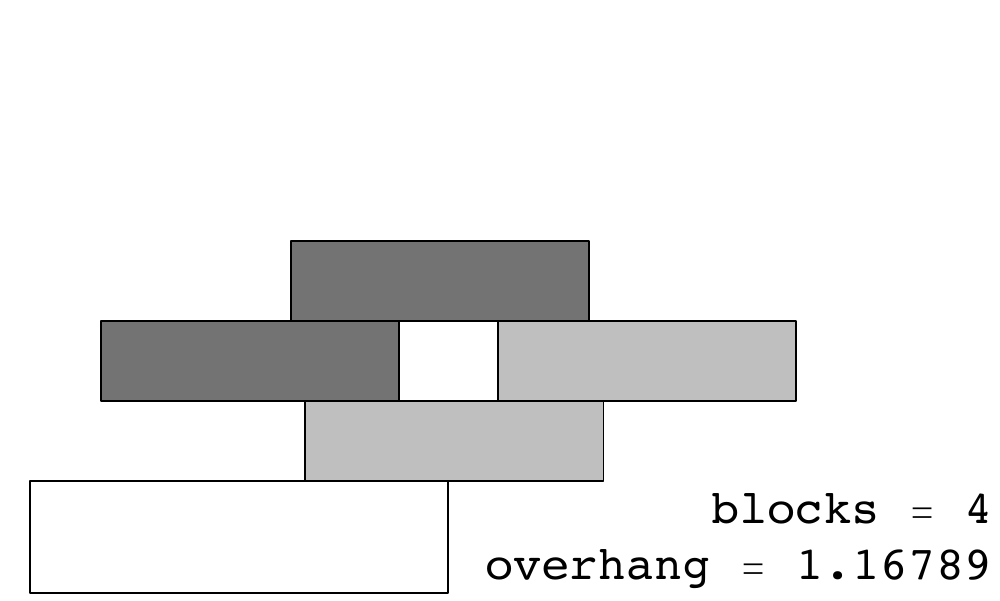}\hspace*{5mm}
}

\centerline{
\includegraphics[height=3cm]{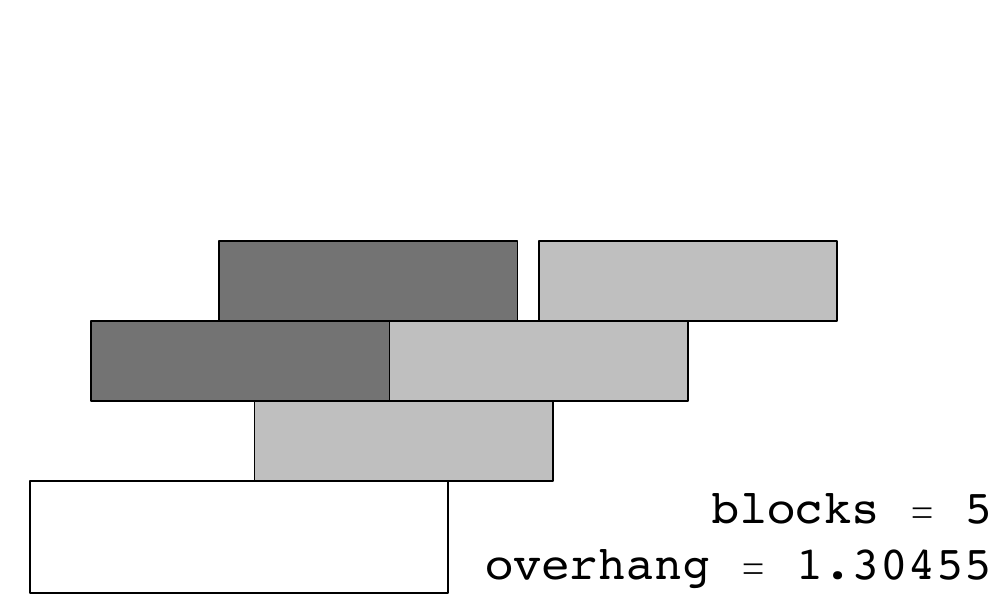}\hspace*{5mm}
\includegraphics[height=3cm]{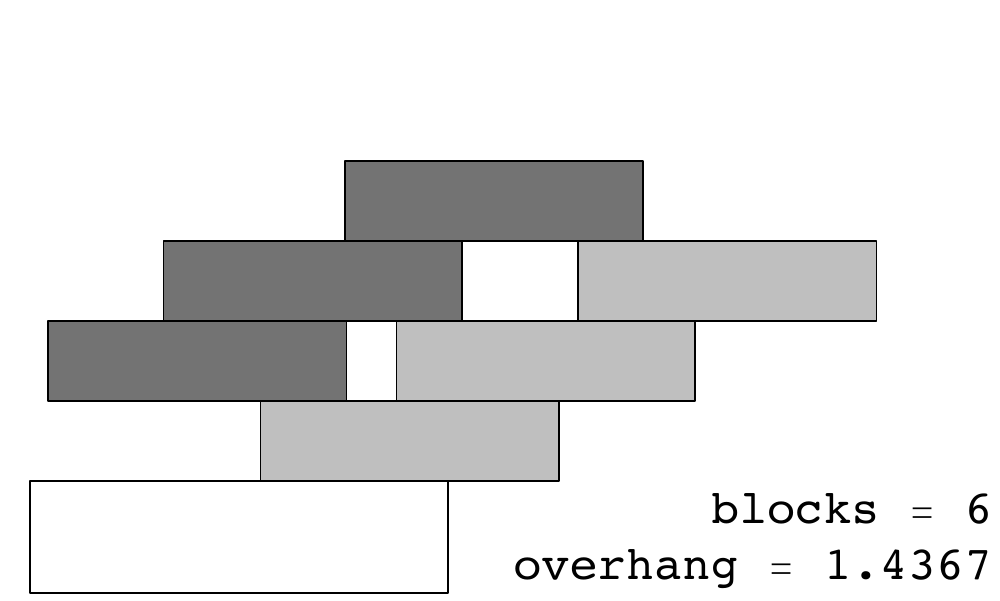}\hspace*{5mm}
\includegraphics[height=3cm]{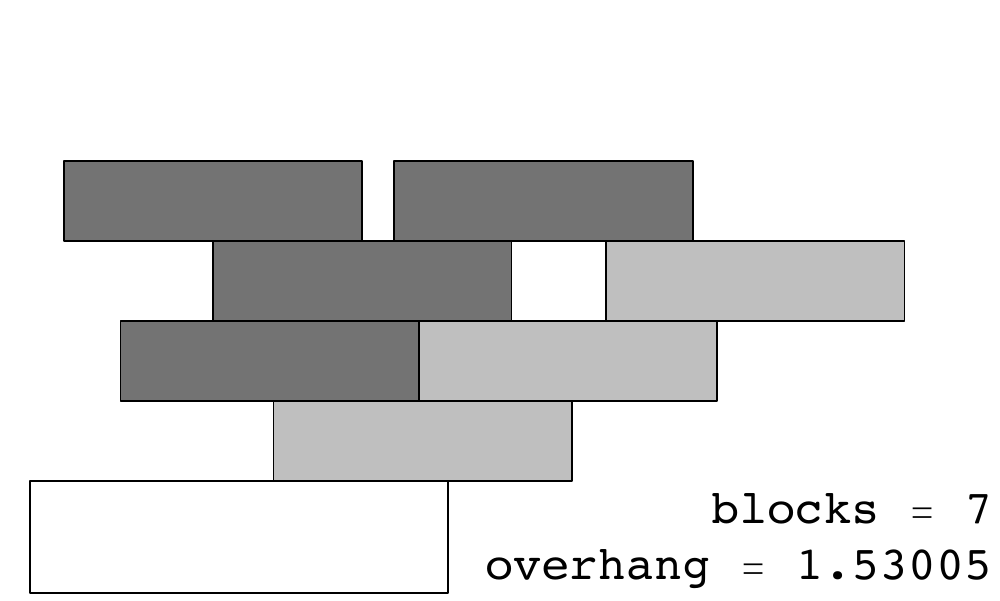}\hspace*{5mm}
}

\vspace*{0.5cm}

 \centerline{
\includegraphics[height=3cm]{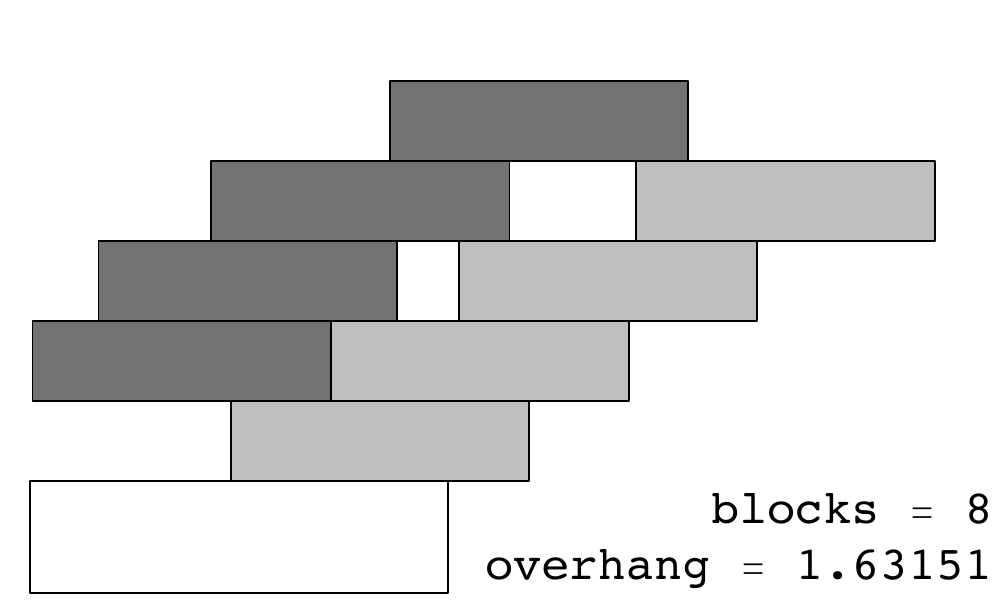}  \hspace*{3mm}
\includegraphics[height=3cm]{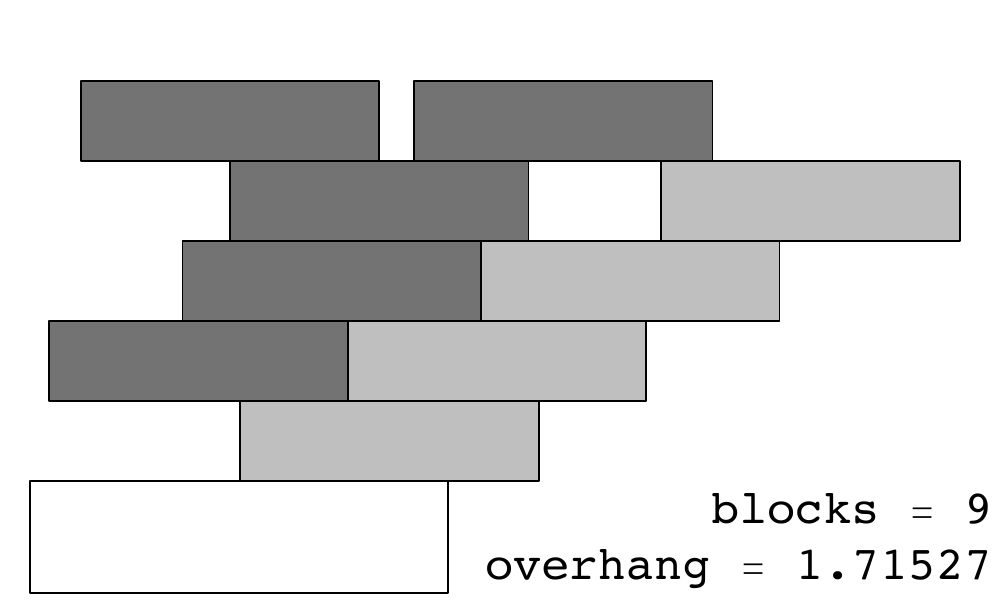}  \hspace*{3mm}
\includegraphics[height=3cm]{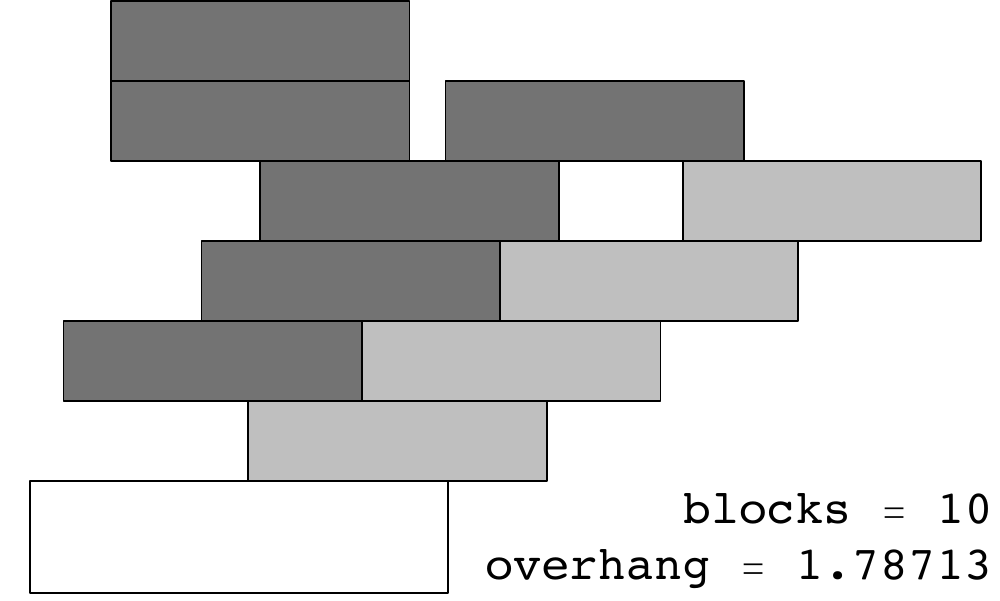} \hspace*{3mm}
}

\vspace*{0.5cm}

\caption{Optimal stacks with 2 up to 10 blocks.} \label{fig:2-10}
\end{figure*}

\begin{figure*}[t]

\centerline{
\includegraphics[height=4cm]{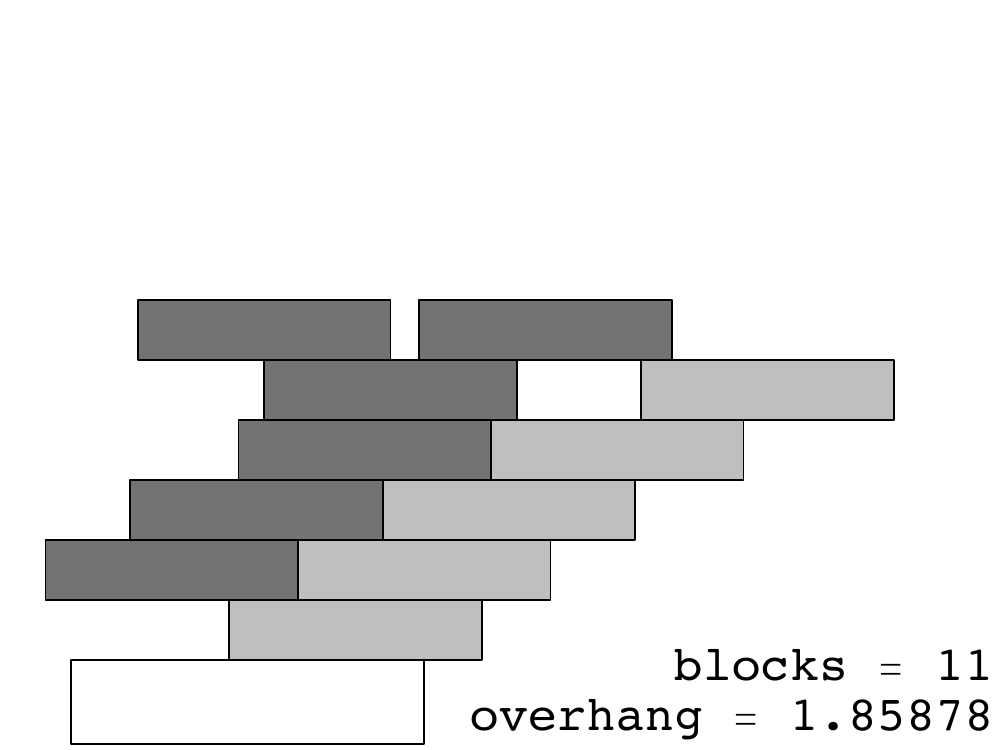}\hspace*{5mm}
\includegraphics[height=4cm]{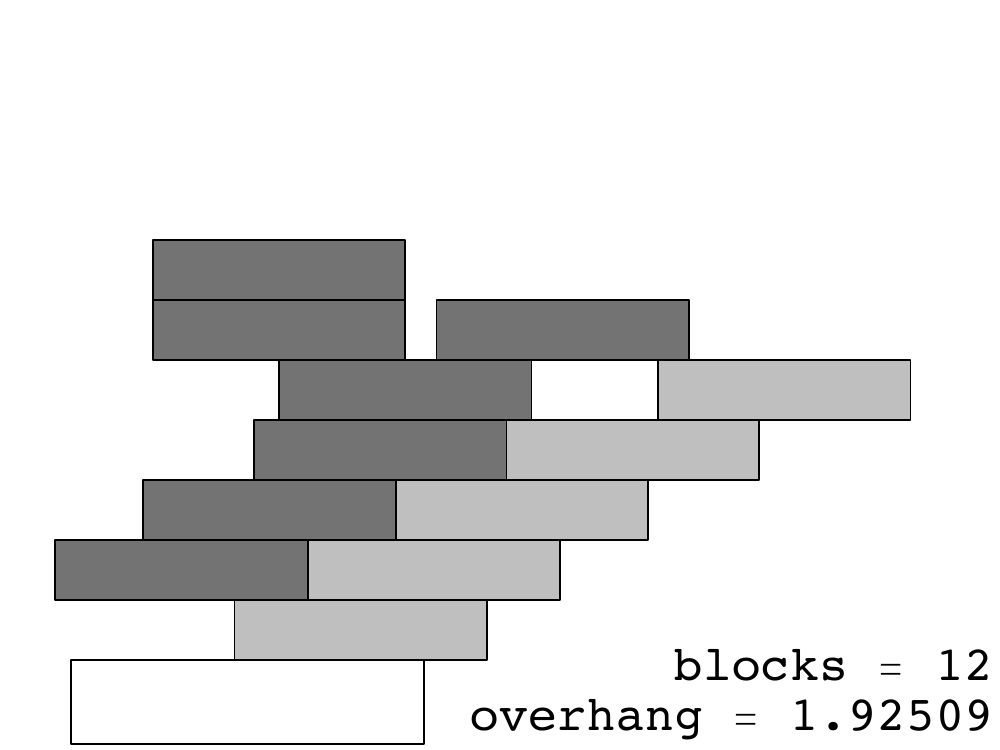}\hspace*{5mm}
\includegraphics[height=4cm]{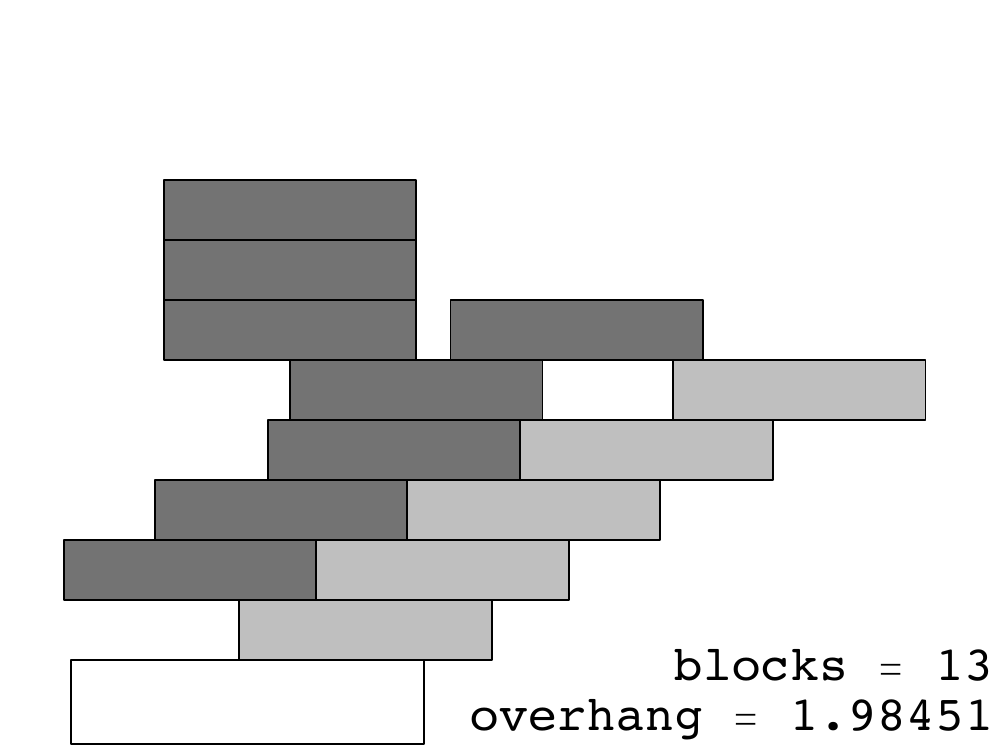}\hspace*{5mm}
}

\centerline{
\includegraphics[height=4cm]{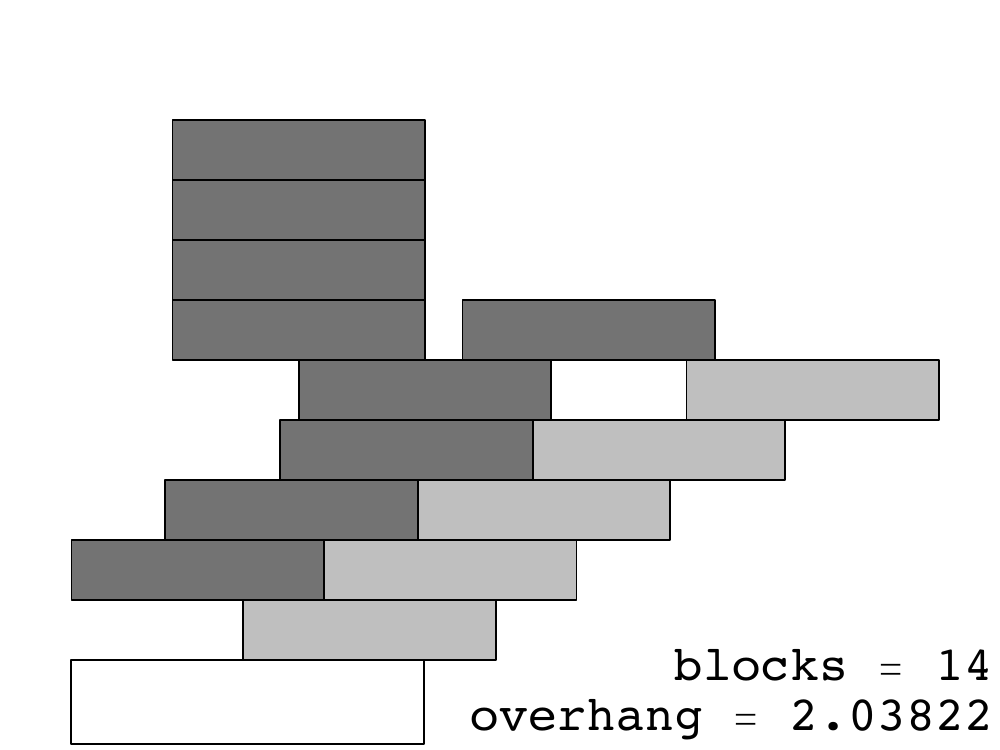}\hspace*{5mm}
\includegraphics[height=4cm]{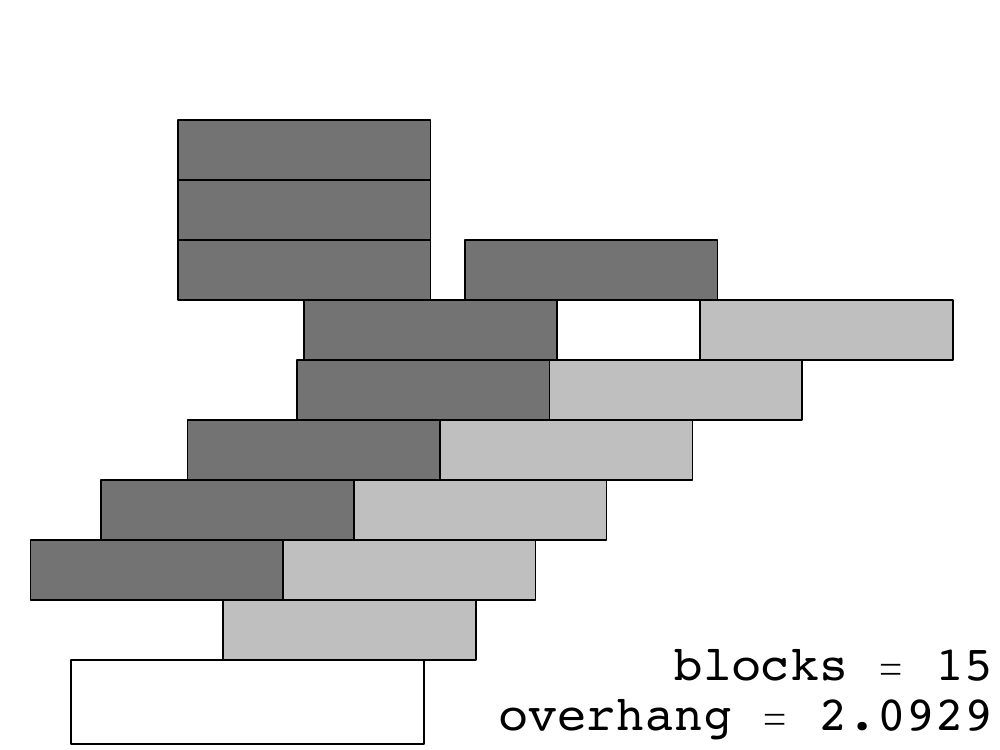}\hspace*{5mm}
\includegraphics[height=4cm]{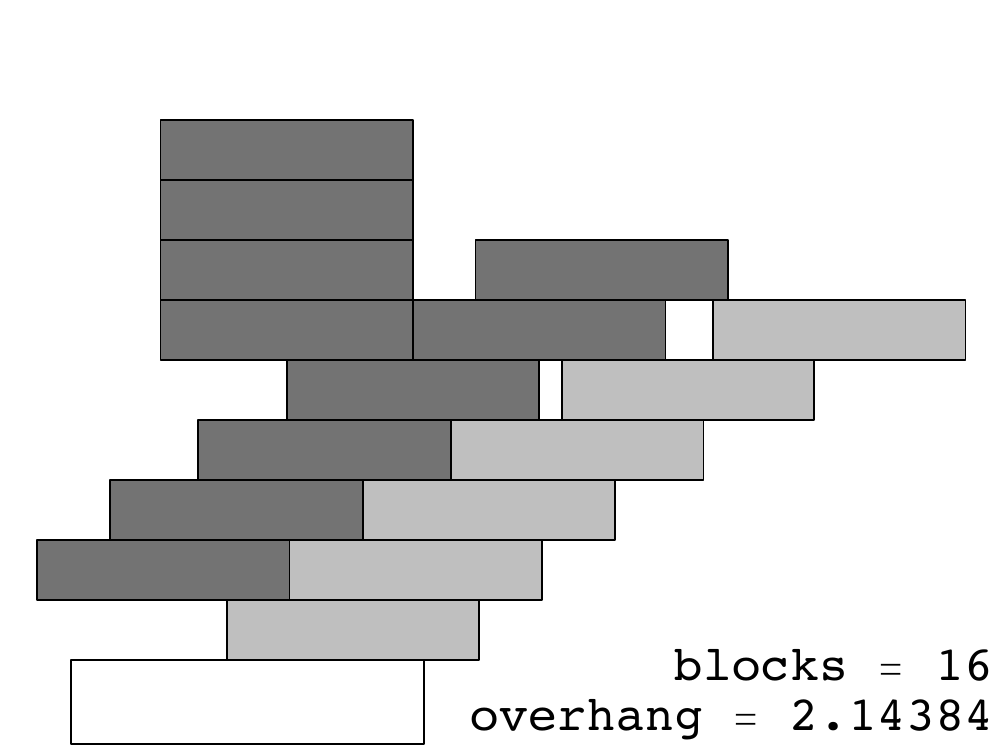}\hspace*{5mm}
}

\vspace*{0.7cm}

\centerline{
\includegraphics[height=4cm]{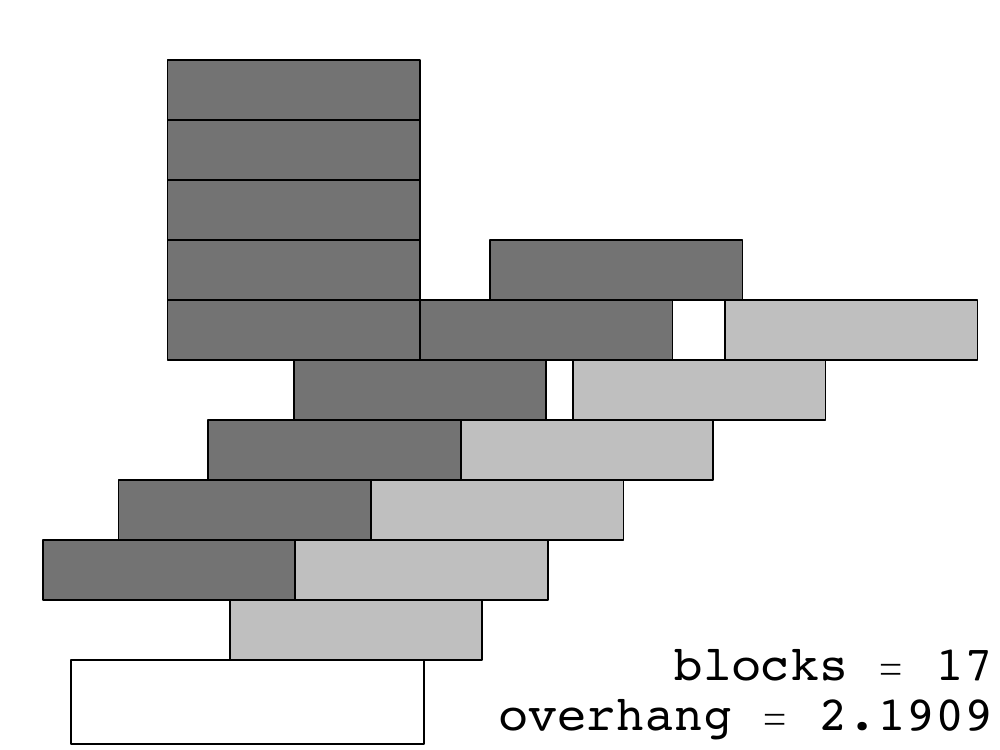}  \hspace*{3mm}
\includegraphics[height=4cm]{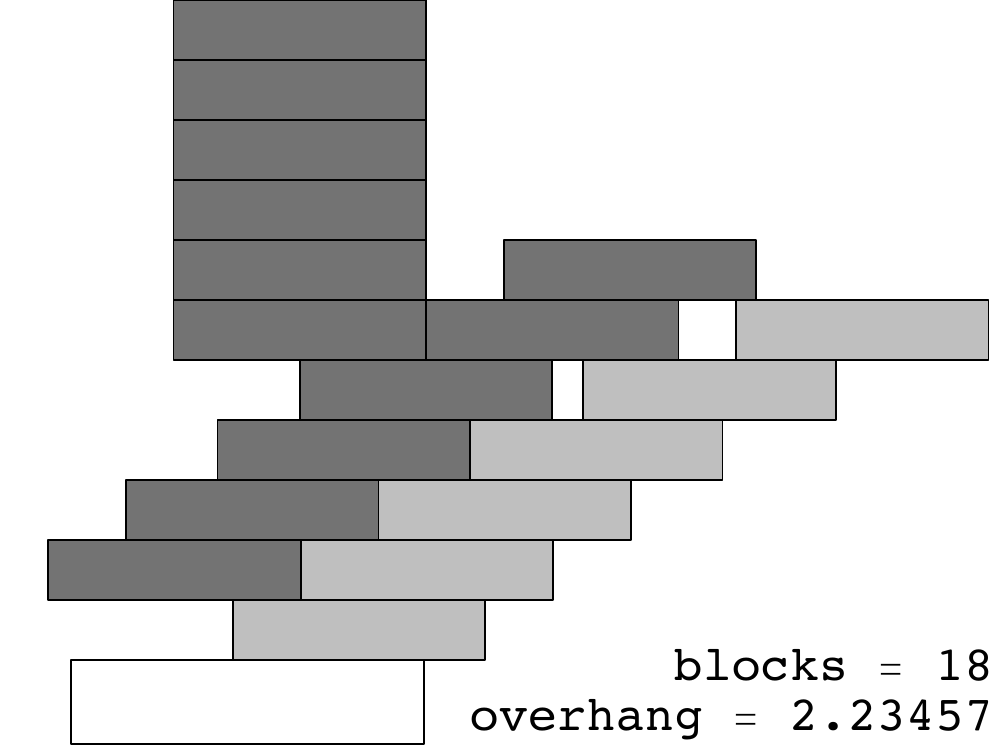}  \hspace*{3mm}
\includegraphics[height=4cm]{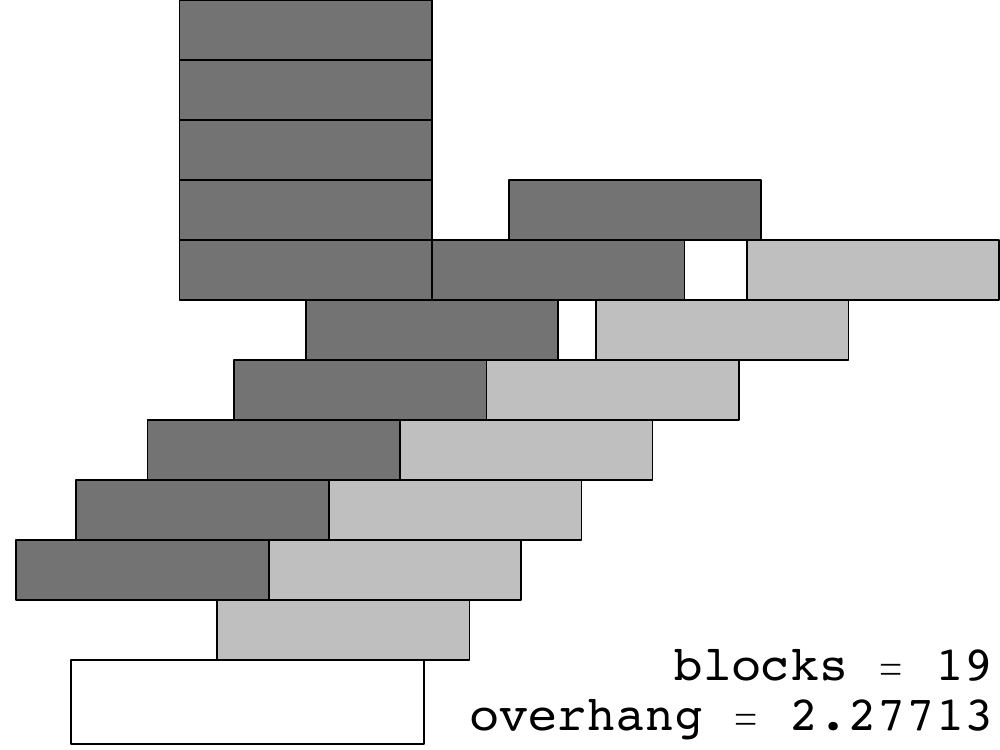} \hspace*{3mm}
}

\vspace*{0.5cm}

\caption{Optimal stacks with 11 up to 19 blocks.} \label{fig:11-19}
\end{figure*}

Optimal stacks with up to 19 blocks are shown in
Figures~\ref{fig:2-10} and~\ref{fig:11-19}. The lightly shaded blocks in
these stacks form the \emph{support set}, while the darker blocks
form the \emph{balancing set}. The \emph{principal block} of a stack
is defined to be the block which achieves the maximum overhang. (If
several blocks achieve the maximum overhang, the lowest one is
chosen.) The \emph{support set} of a stack is defined recursively as
follows: the principal block is in the support set, and if a block
is in the support set then any block on which this block rests is
also in the support set. The \emph{balancing set} consists of all
the blocks that do not belong to the support set. A stack is said to
be \emph{spinal} if its support set has a single block in each
level, up to the level of the principal block. All the stacks shown
in Figures~\ref{fig:2-10} and~\ref{fig:11-19} are thus spinal.

\begin{figure}[t]
\centerline{
\includegraphics[height=6.25cm]{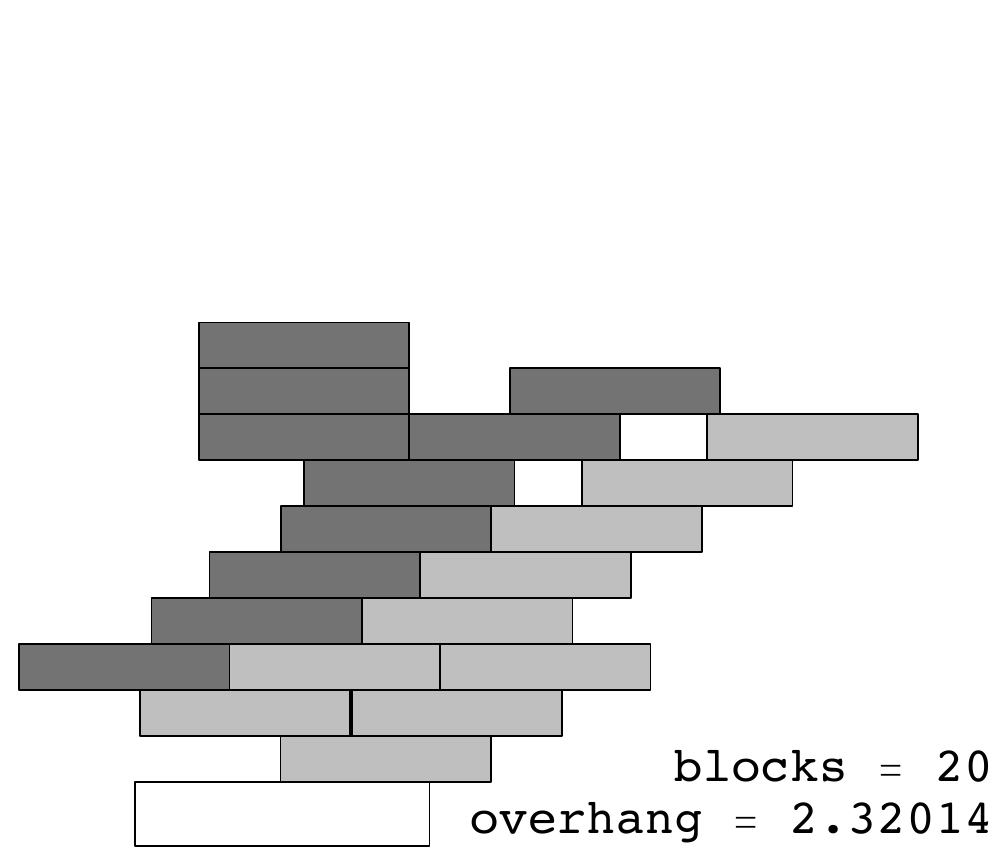}  \hspace*{3mm}
\includegraphics[height=6.25cm]{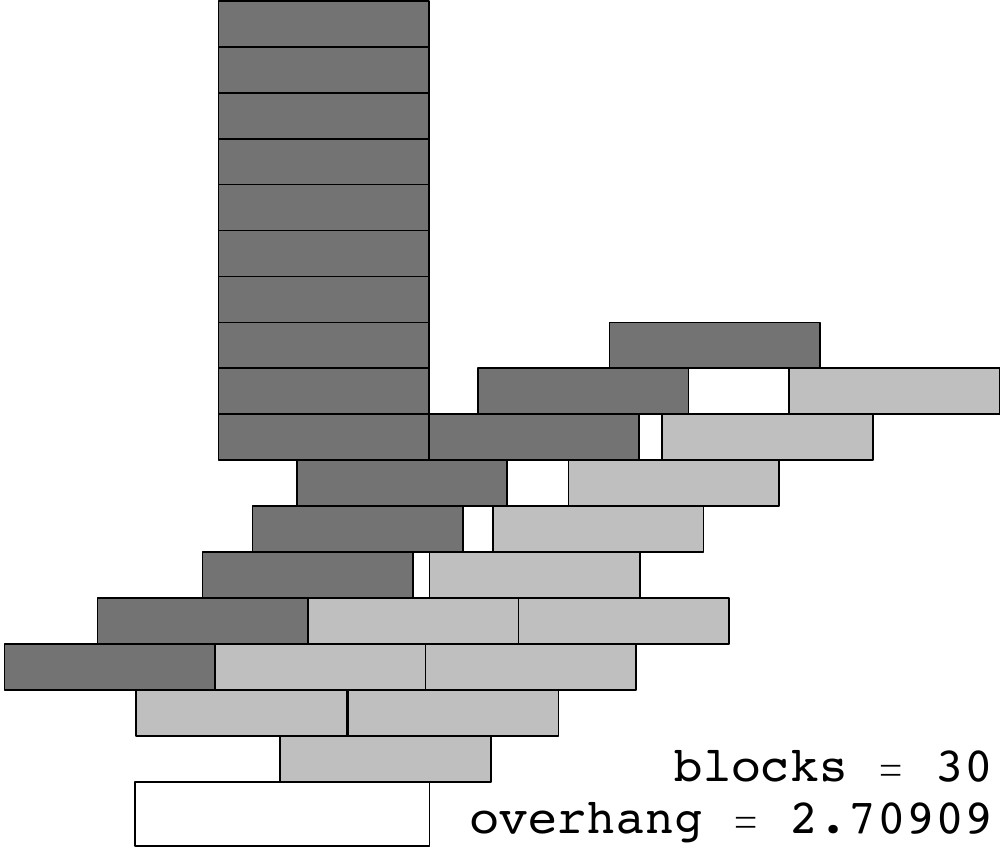}
} \vspace*{0.5cm}
\caption{Optimal stacks with 20 and 30 blocks.} \label{fig:20-30}
\end{figure}

It is very tempting to conclude, as done by Hall \cite{H05}, that
the optimal stacks are spinal. Surprisingly, the optimal stacks for
$n\geqslant 20$ are not spinal! Optimal stacks containing 20 and 30 blocks
are shown in Figure~\ref{fig:20-30}. Note that the right-hand contours of
these stacks are not monotone, which is somewhat counterintuitive.

For all $n\leqslant 30$, we have searched exhaustively through all
combinatorially distinct arrangements of $n$ blocks and found
optimal displacements numerically for each of these. The resulting
stacks, for $2\leqslant n\leqslant 19$ are shown in Figures~\ref{fig:2-10}
and~\ref{fig:11-19}. Optimal stacks with 20 and 30 blocks are shown
in Figure~\ref{fig:20-30}. We are confident of their optimality,
though we have no formal optimality proofs, as numerical techniques
were used.

While there seems to be a unique optimal placement of the blocks
that belong to the support set of an optimal stack, there is usually
a lot of freedom in the placement of the balancing blocks. Optimal
stacks seem not to be unique for $n\geqslant 4$.

\begin{figure*}[t]

\centerline{
\includegraphics[height=6.25cm]{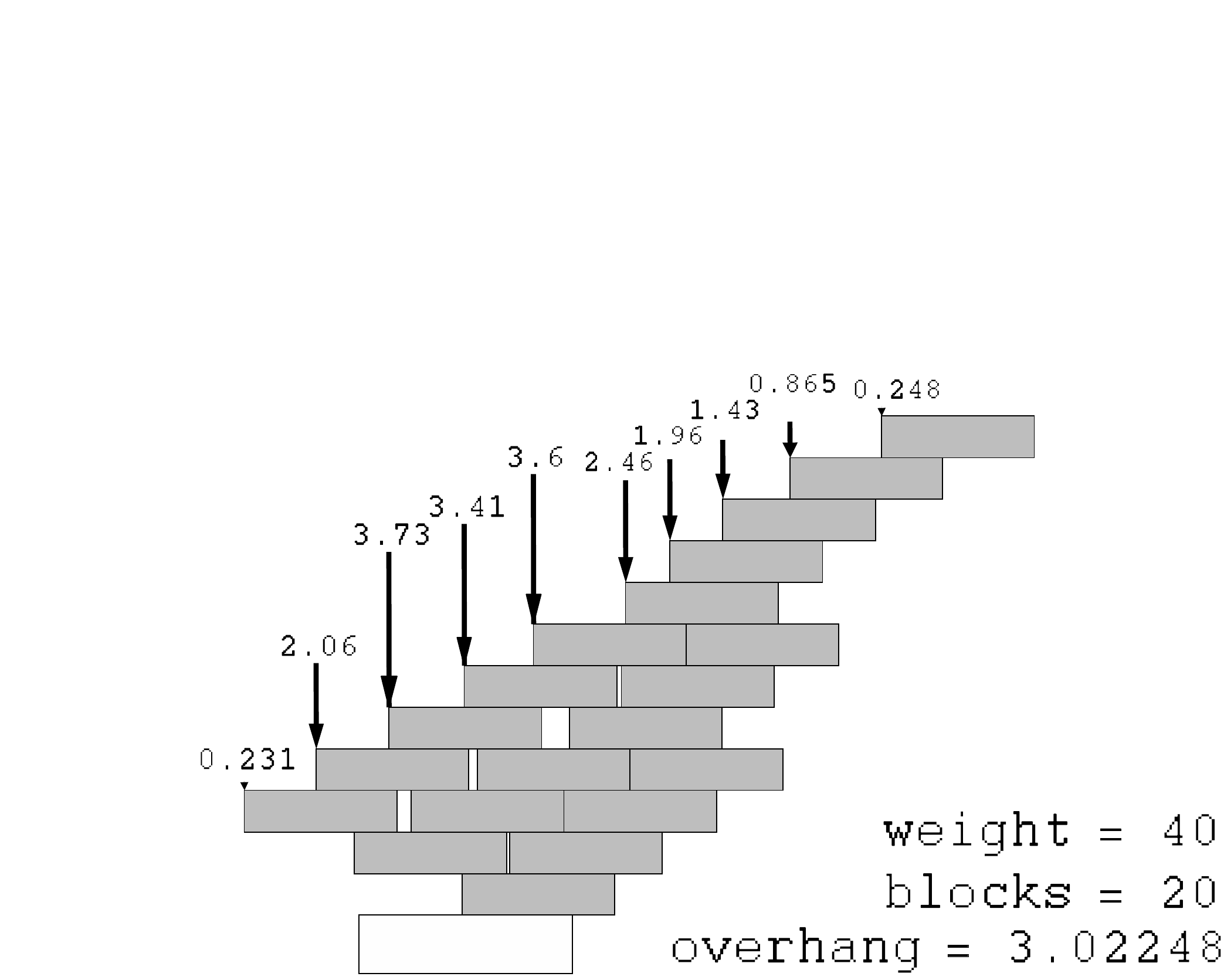}\hspace*{3mm}
\includegraphics[height=6.25cm]{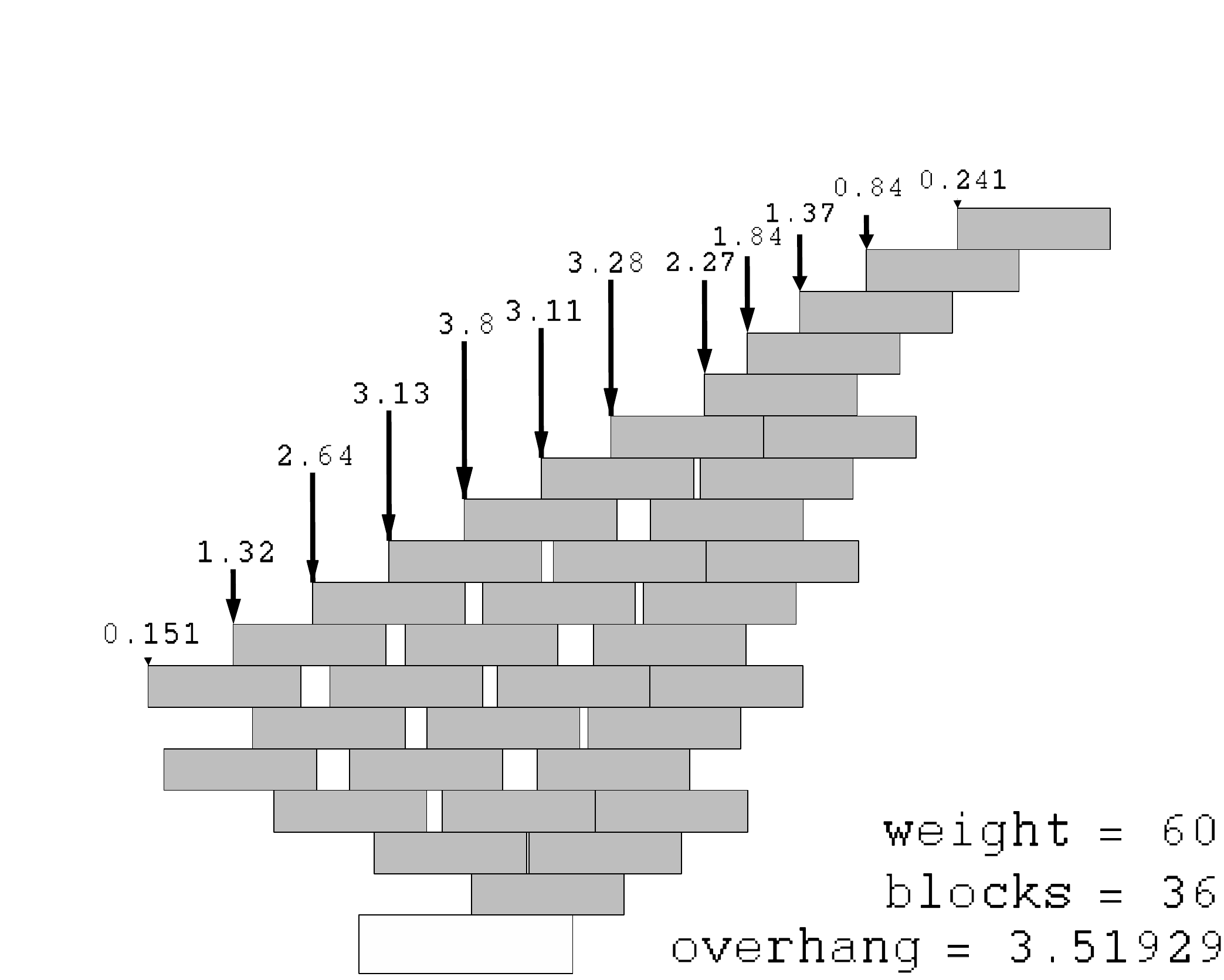}
} \vspace*{0.3cm}
 \centerline{
\includegraphics[height=6.25cm]{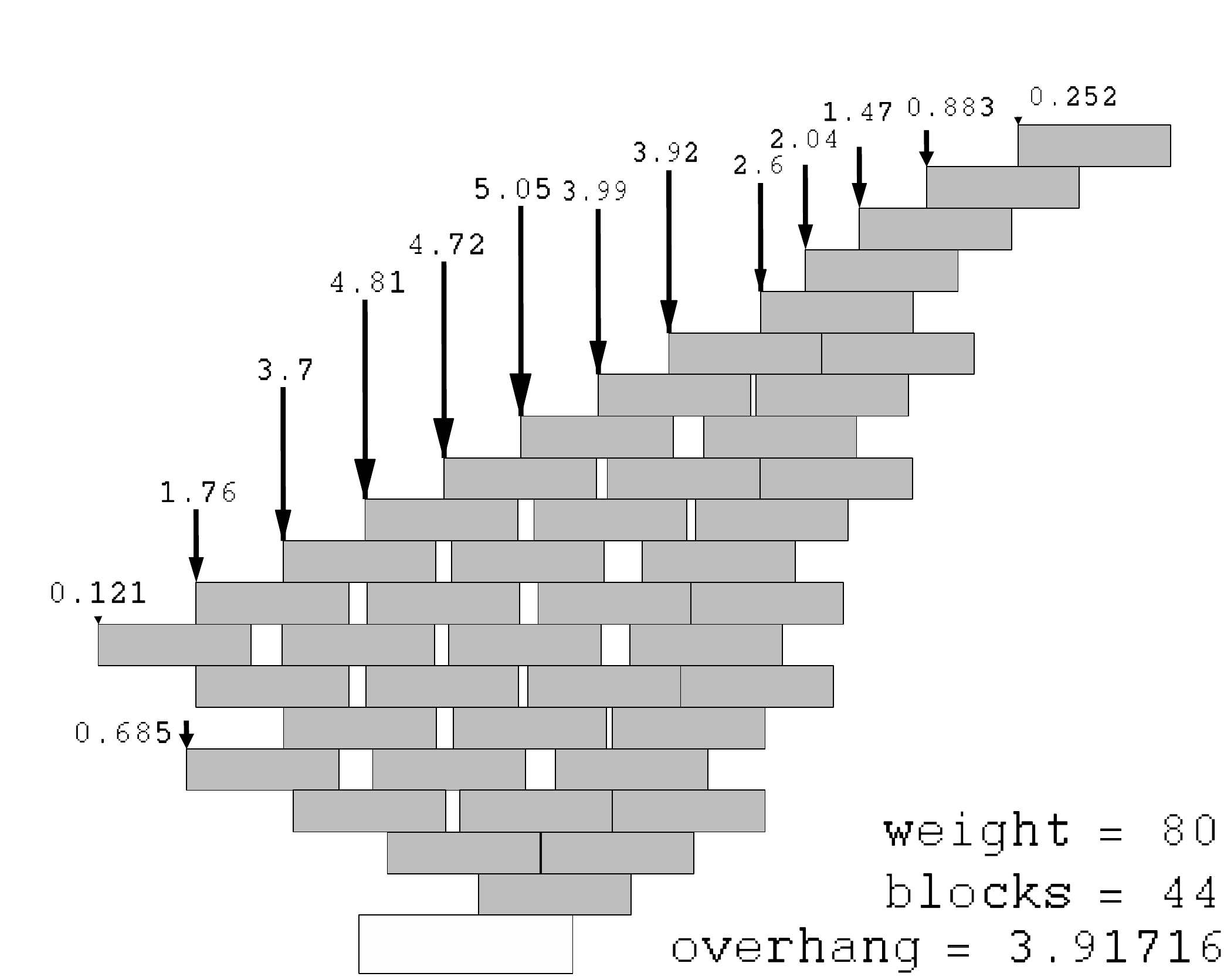}\hspace*{3mm}
\includegraphics[height=6.25cm]{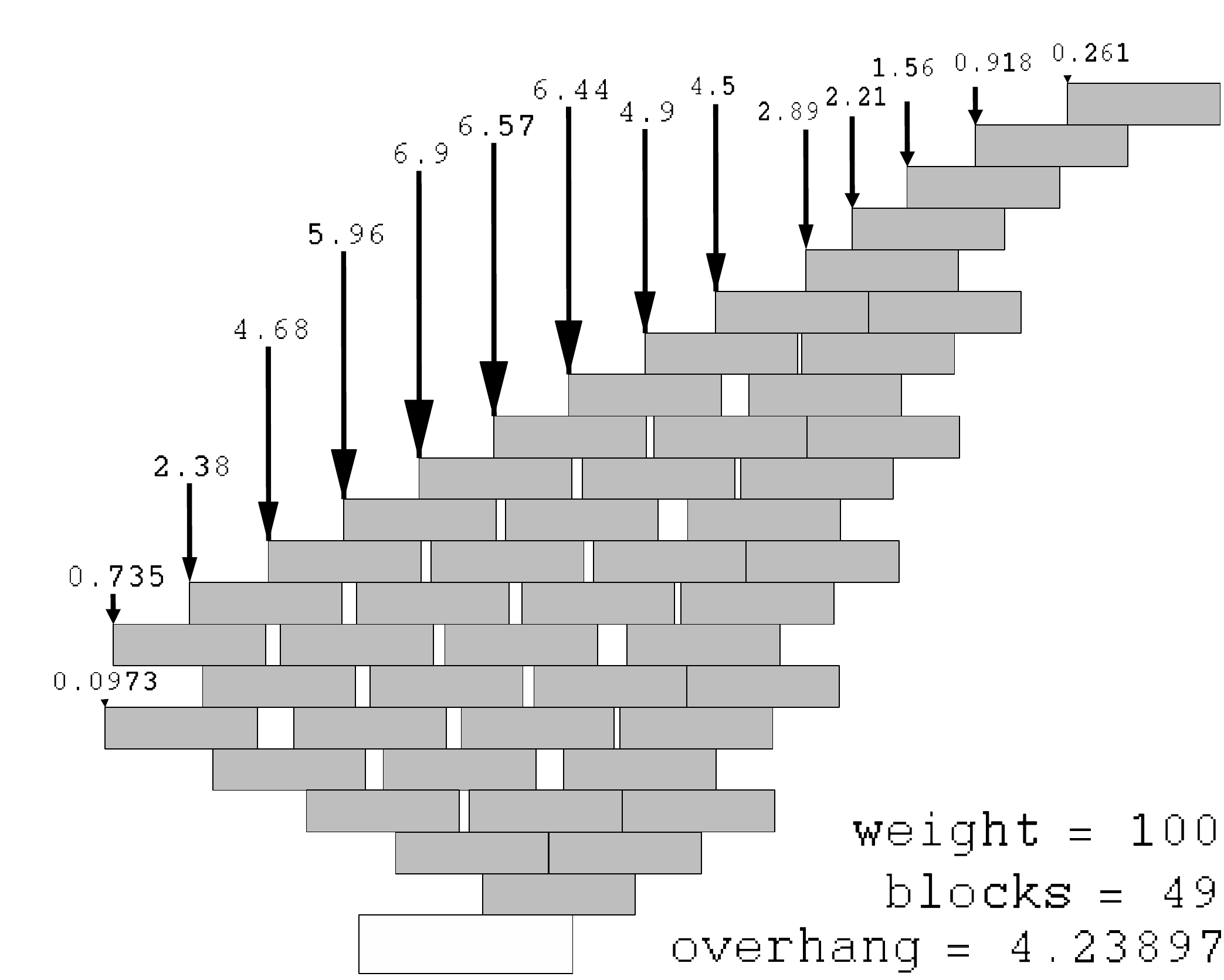}
}

\vspace*{0.4cm}

\caption{Loaded stacks, believed to be close to optimal, of weight
$40,60,80,$ and $100$.}

\label{fig:p40-100}
\end{figure*}

In view of the non-uniqueness and added complications caused
by balancing blocks, it is
natural to consider \emph{loaded stacks}, which consist only of a
support set with some \emph{external forces} (or \emph{point
weights}) attached to some of their blocks.
We will take the weight of each block to be~$1$; the size, or weight, of
a loaded stack is defined to be the number of blocks contained in it
plus the sum of all the point weights attached to it. The point
weights are not required to be integral. Loaded stacks of weight
$40,60,80,$ and $100$, which are believed to be close to optimal, are
shown in Figure~\ref{fig:p40-100}. The stack of weight~$100$, for
example, contains $49$ blocks in its support set. The sum of all the
external forces applied to these blocks is $51$. As can be seen, the
stacks become more and more non-spinal. It is also interesting to
note that the stacks of Figure~\ref{fig:p40-100} contain small gaps
that seem to occur at irregular positions. (There is also a scarcely
visible gap between the two blocks at the second level of the
20-block stack of Figure~\ref{fig:20-30}.)

That harmonic stacks are balanced can be verified using simple
center-of-mass considerations. These considerations, however, are
not enough to verify the balance of more complicated stacks,
such as those in
Figures~\ref{fig:2-10}, \ref{fig:11-19}, \ref{fig:20-30},
and~\ref{fig:p40-100}. A formal mathematical definition of ``balanced''
is given in the next section. Briefly, a stack is said to be balanced
if there is an appropriate set of forces acting between the blocks
of the stacks, and between the blocks at the lowest level and the
table, under which all blocks are \emph{in equilibrium}. A block is
in equilibrium if the sum of the forces and the sum of the moments
acting upon it are both 0. As shown in the next section, the
balance of a given stack can be determined by checking whether a
given set of \emph{linear inequalities} has a feasible solution.

\begin{figure}[t]
\begin{center}
\includegraphics[height=20mm]{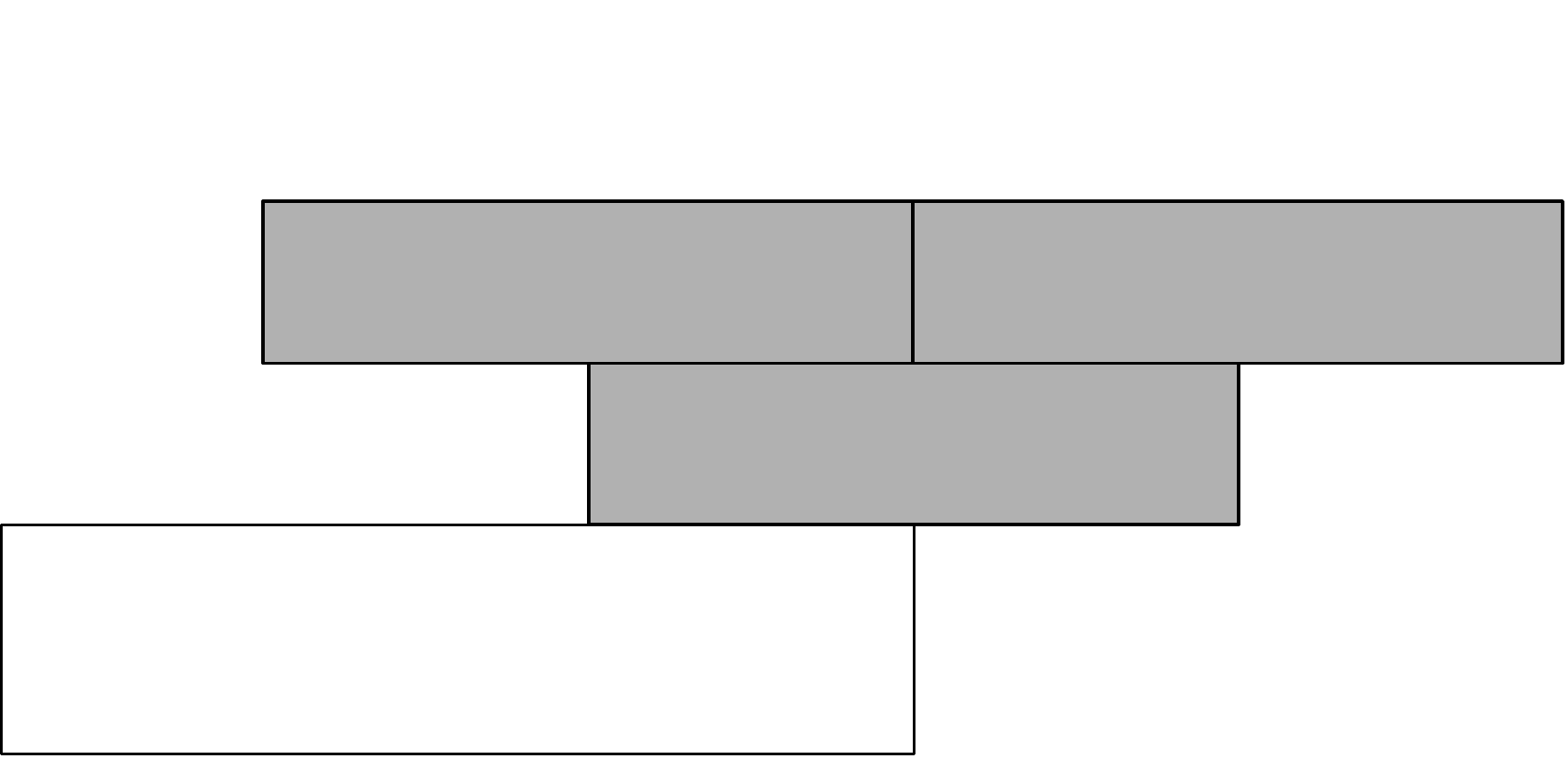}\hspace*{1cm}
\includegraphics[height=20mm]{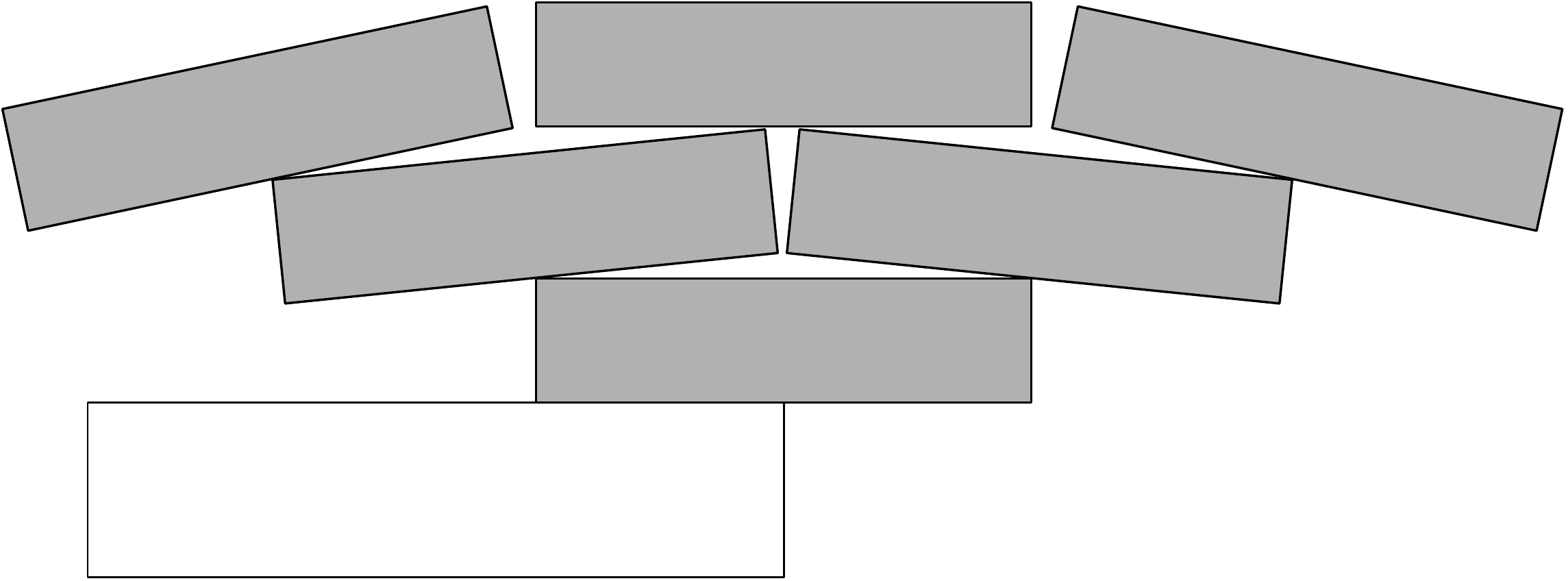}
\caption{The balanced inverted 2-triangle and the unbalanced
inverted 3-triangle.} \label{fig:p23}
\end{center}
\end{figure}

Given the fact that the 3-block stack that achieves the maximum
overhang is an \emph{inverted 2-triangle} (see
Figure~\ref{fig:opt34}), it is natural to enquire whether larger
inverted triangles are also balanced. Unfortunately, the next inverted
triangle is already unbalanced and would collapse in the way
indicated in Figure~\ref{fig:p23}. Inverted triangles show that
simple center-of-mass considerations are not enough to determine the
balance of stacks. As an indication that balance issues are not
always intuitive, we note that inverted triangles are falsely
claimed by Jargodzki and Potter~\cite{JP01} (Challenge~271: A
staircase to infinity, p.~246) to be balanced.

\begin{figure}[t]
\begin{center}
\includegraphics[height=30mm]{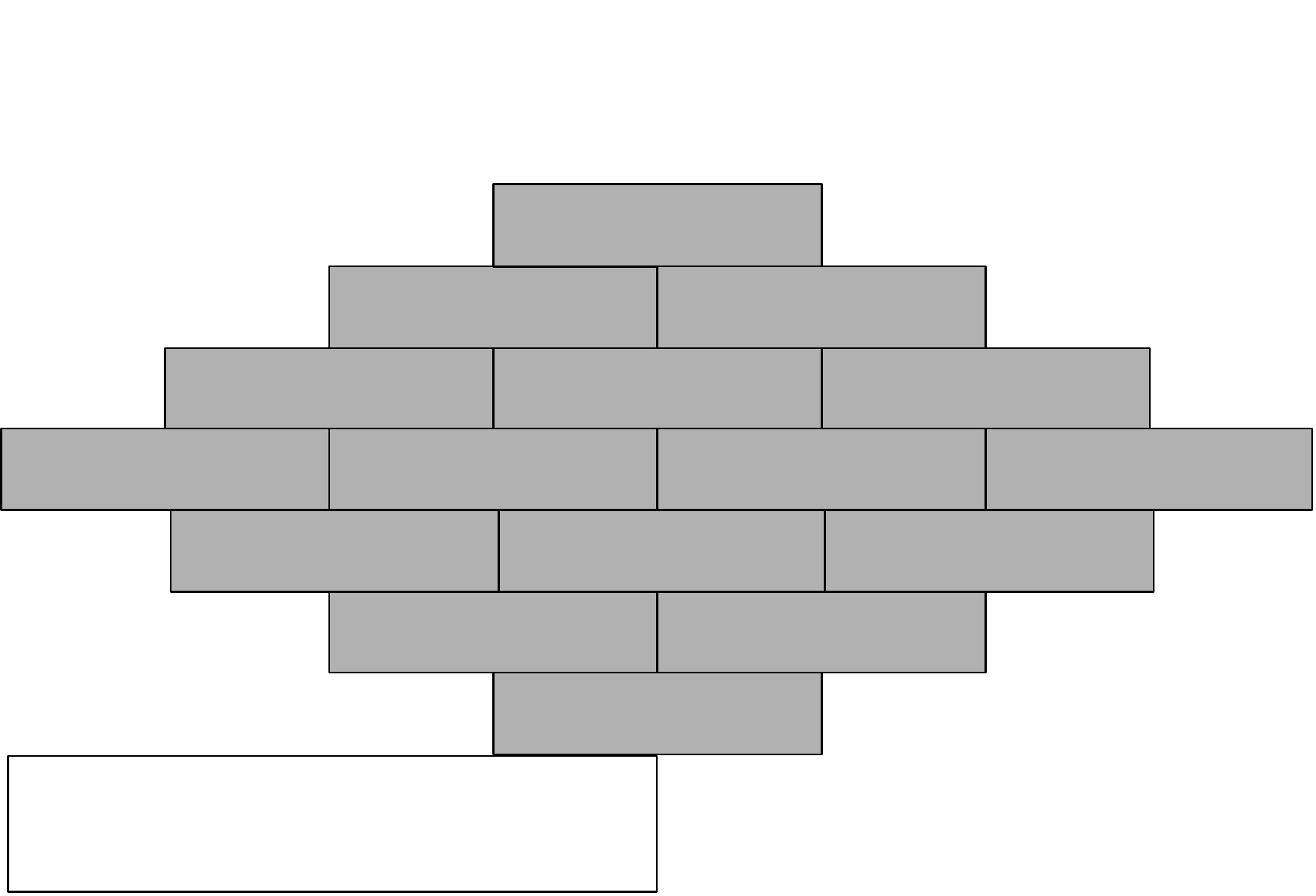}\hspace*{1cm}
\includegraphics[height=30mm]{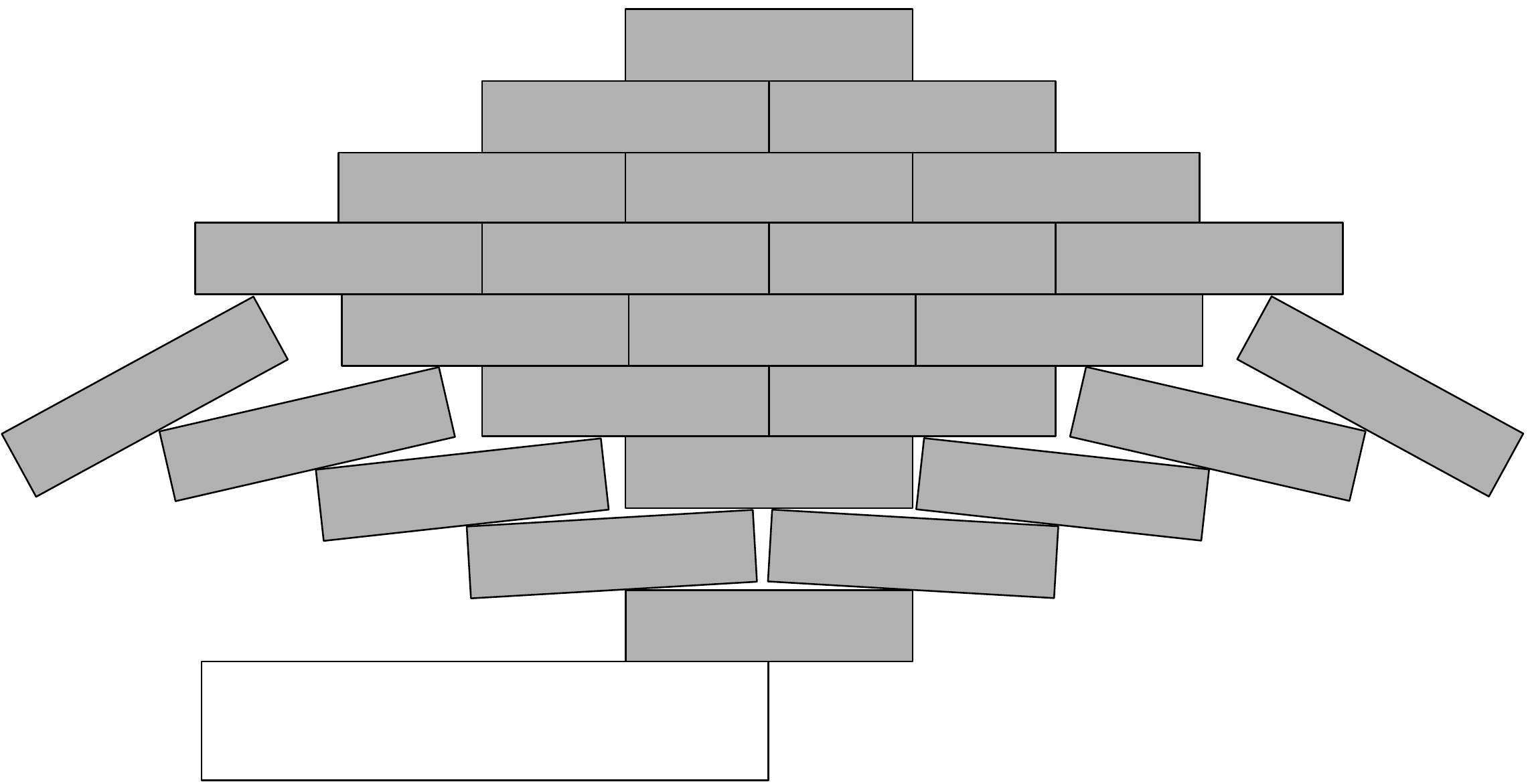}
\caption{The balanced 4-diamond and the unbalanced 5-diamond.}
\label{fig:d45}
\end{center}
\end{figure}

Another appealing structure, the \emph{$m$-diamond}, illustrated for
$m=4$ and~$5$ in Figure~\ref{fig:d45}, consists of a symmetric
diamond shape with rows of length $1,2,\ldots,m-1,m,m-1,\ldots,2,1$.
Small diamonds were considered by Drummond \cite{D81}. The
$m$-diamond uses~$m^2$ blocks and would give an overhang of $m/2$,
but unfortunately it is unbalanced for $m\geqslant 5$. A $5$-diamond
would collapse in the way indicated in the figure. An $m$-diamond
could be made balanced by adding a column of sufficiently many
blocks resting on the top block. The methodology introduced in
Section~\ref{sec:spinal} can be used to show that, for $m\geqslant
5$, a column of at least $2^m-m^2-1$ blocks would be needed. We can
show that this number of blocks is also sufficient, giving a stack
of $2^m-1$ blocks with an overhang of $m/2$. It is interesting to
note that these stacks are already better than the classical
harmonic stacks, as with $n=2^m-1$ blocks they give an overhang of
$\frac{1}{2}\log_2(n+1)\simeq 0.693\ln n$.

Determining the \emph{exact} overhang achievable using $n$ blocks,
for large values of~$n$, seems to be a formidable task. Our main
goal in this paper is to determine the \emph{asymptotic growth} of
this quantity. Our main result is that there exists a constant $c>0$
such that an overhang of $cn^{1/3}$ is achievable using $n$ blocks.
Note that this is an exponential improvement over the
$\frac{1}{2}\ln n+O(1)$ overhang of harmonic stacks and the $\ln
n+O(1)$ overhang of the best spinal stacks! In a subsequent paper
\cite{PPTWZ07}, with three additional coauthors, we show that our
improved stacks are asymptotically optimal, i.e., there exists a
constant $C>0$ such that the overhang achievable using $n$ blocks is
at most $Cn^{1/3}$.

\begin{figure}[t]
\vspace*{0.25cm}
\begin{center}
\includegraphics[width=90mm]{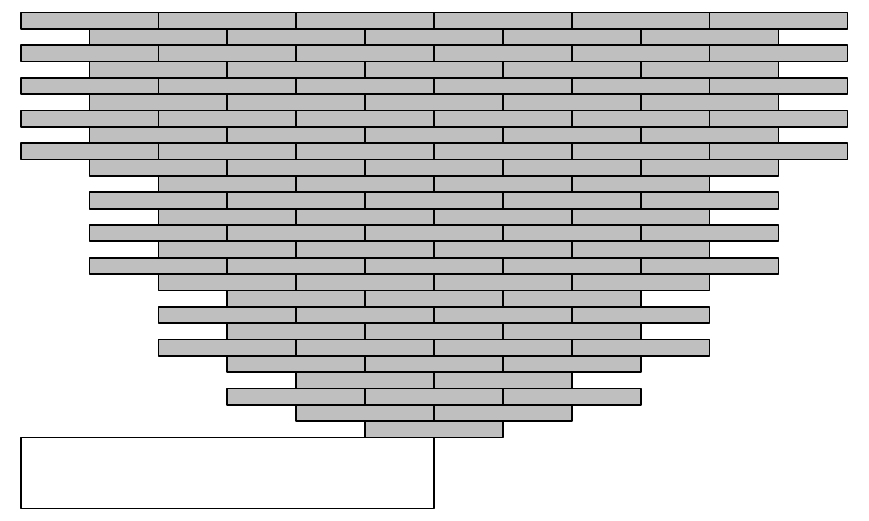}
\caption{A parabolic stack consisting of 111 blocks and
giving an overhang of $3$.} \label{fig:6stack}
\end{center}\vspace{-5mm}
\end{figure}

Our stacks that achieve an asymptotic overhang of $cn^{1/3}$, for
some $c>0$, are quite simple. We construct an explicit sequence of
stacks, called \emph{parabolic stacks}, with the $r$th stack
in the sequence containing about $2r^3/3$ blocks and achieving an
overhang of $r/2$. One stack in this sequence is shown in
Figure~\ref{fig:6stack}. The balance of the parabolic stacks is
established using an inductive argument.

The remainder of this paper is organized as follows. In the next
section we give formal definitions of all the notions used in this
paper. In Section~\ref{sec:spinal} we analyze \emph{spinal stacks}.
In Section~\ref{sec:parabolic}, which contains our main results, we
introduce and analyze our parabolic stacks. In
Section~\ref{sec:general} we describe some experimental results with
stacks that seem to improve, by a constant factor, the overhang
achieved by parabolic stacks. We end in Section~\ref{sec:conc} with
some open problems.

\section{Stacks and their balance} \label{sec:prelim}

As the maximum overhang problem is physical in nature, our first
task is to formulate it mathematically. We consider a 2-dimensional
version of the problem. This version captures essentially all the
interesting features of the overhang problem.

A \emph{block} is a rectangle of length~$1$ and height~$h$ with
uniform density and unit weight. (We shall see shortly that the
height~$h$ is unimportant.) We assume that the table occupies the
quadrant $x,y\leqslant 0$ of the 2-dimensional plane.
A \emph{stack} is a collection of blocks.
We consider only \emph{orthogonal} stacks in which the sides of the
blocks are parallel to the axes, with the length of each block
parallel to the $x$-axis. The position of a block is then determined
by the coordinate $(x,y)$ of its lower left corner. Such a block
occupies the box $[x,x+1]\times[y,y+h]$. A stack composed
of~$n$ blocks is specified by the sequence
$(x_1,y_1),\ldots,(x_n,y_n)$ of the coordinates of the lower left
corners of its blocks. We require each $y_i$ to be a nonnegative
integral multiple of~$h$, the height of the blocks. Blocks can touch
each other but are not allowed to overlap. The \emph{overhang} of
the stack is $1+\max_{i=1}^n x_i$.

A block at position $(x_1,y_1)$ \emph{rests on} a block in position
$(x_2,y_2)$ if $|x_1-x_2| < 1$ and $y_1-y_2=h$. The
\emph{interval of contact} between the two blocks is then
$[\max\{x_1,x_2\},1+\min\{x_1,x_2\}]\times \{y_1\}$. A block placed
at position $(x,0)$ \emph{rests on the table} if $x < 0$. The
interval of contact between the block and the table is
$[x,\min\{x+1,0\}]\times\{0\}$.

\begin{figure*}[t]
\begin{center}
\includegraphics[height=35mm]{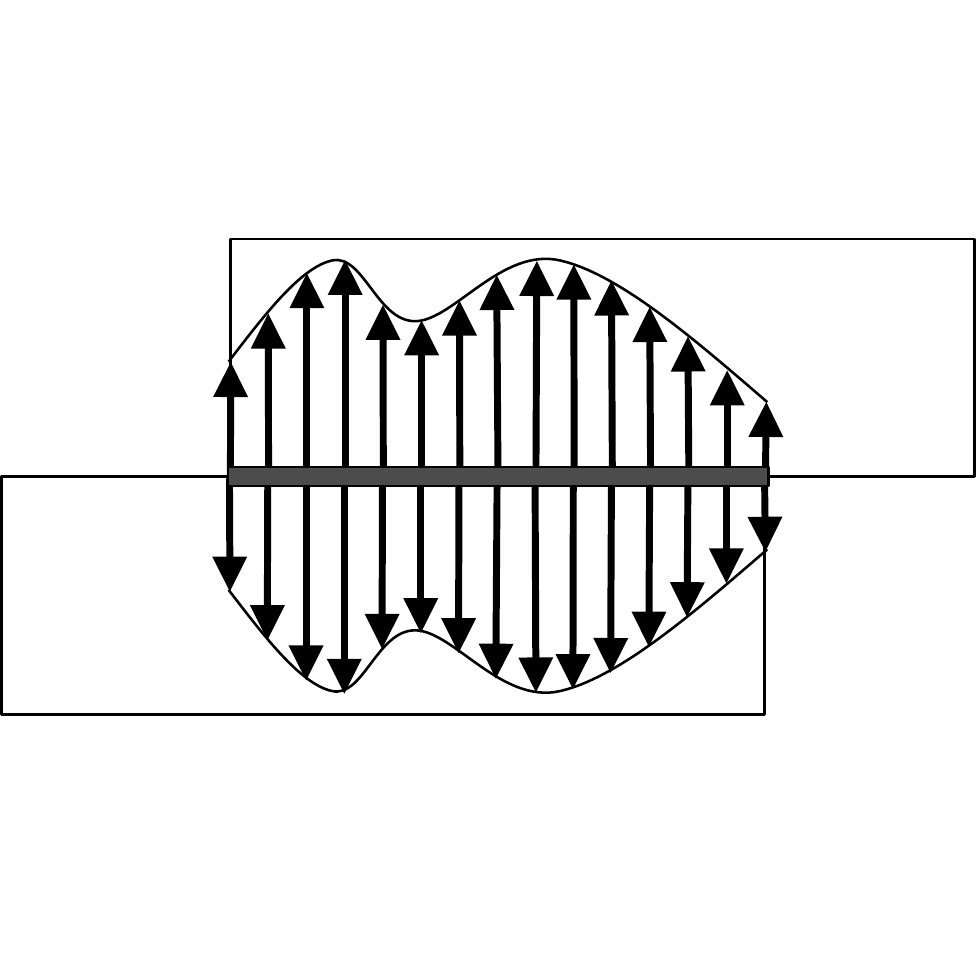}\hspace*{1cm}
\includegraphics[height=35mm]{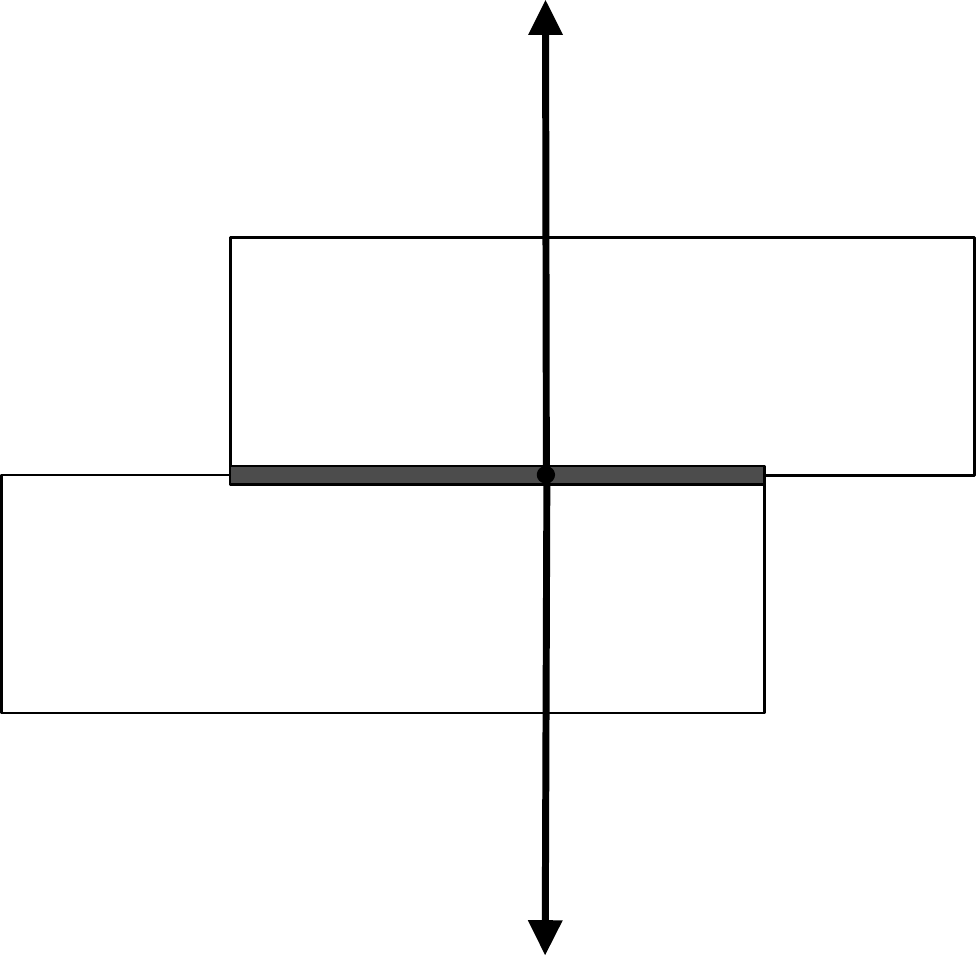}\hspace*{1cm}
\includegraphics[height=35mm]{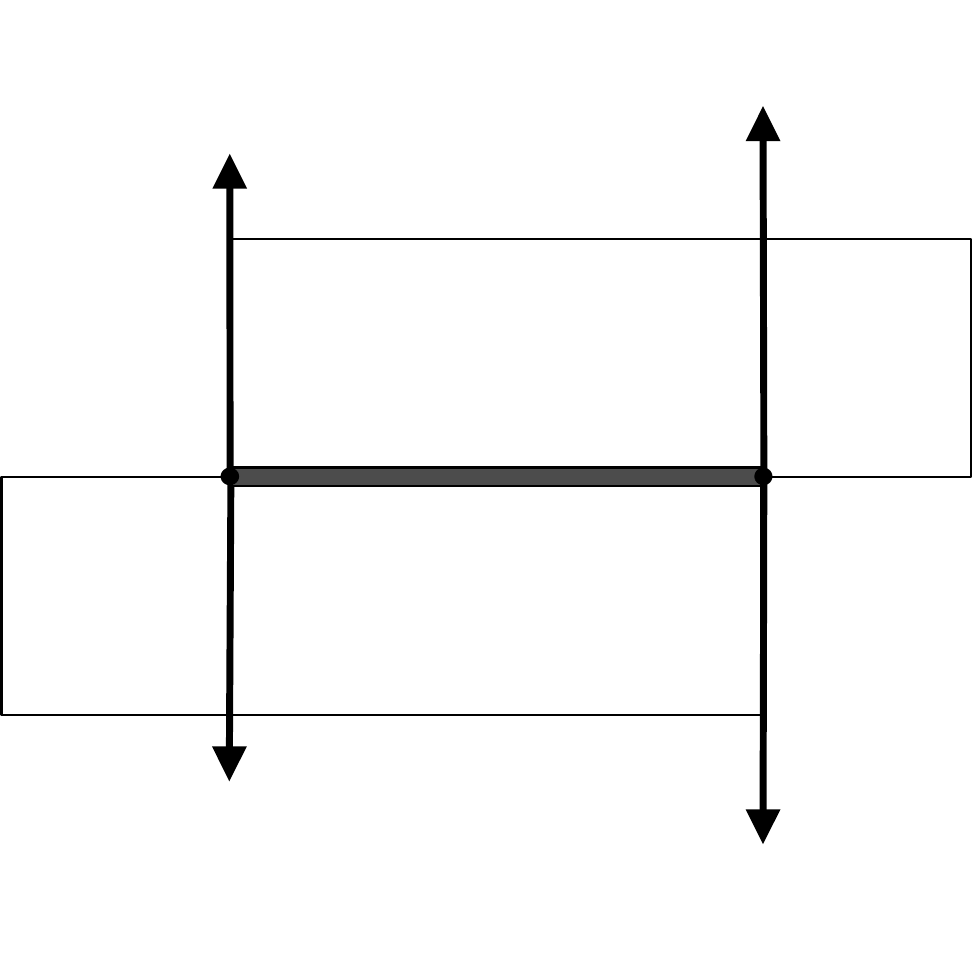}
\caption{Equivalent sets of forces acting between two blocks.}
\vspace*{-0.2cm} \label{fig:forces}
\end{center}
\end{figure*}

When block~$A$ rests on block~$B$, the two blocks may exert a
(possibly infinitesimal) force on each other at every point along
their interval of contact. A force is a vector acting at a specified
point. By Newton's third law, forces come in opposing pairs. If a
force~$f$ is exerted on block~$A$ by block~$B$, at $(x,y)$, then a
force~$-f$ is exerted on block~$B$ by block~$A$, again at $(x,y)$.
We assume that edges of all the blocks are completely smooth, so
that there is no \emph{friction} between them. All the forces
exerted on block~$A$ by block~$B$, and vice versa, are therefore
\emph{vertical} forces. Furthermore, as there is nothing that holds
the blocks together, blocks $A$ and $B$ can \emph{push}, but not
pull, one another. Thus, if block~$A$ rests on block~$B$, then all
the forces applied on block~$A$ by block~$B$ point upward, while all
the forces applied on block~$B$ by block~$A$ point downward, as
shown on the left in Figure~\ref{fig:forces}. Similar forces
are exerted between the table and the blocks that rest on it.

The distribution of forces acting between two blocks may be hard to
describe explicitly. Since all these forces point in the same
direction, they can always be replaced by a single \emph{resultant}
force acting at a some point within their interval of contact, as
shown in the middle drawing of Figure~\ref{fig:forces}. As an
alternative, they may be replaced by two resultant forces that act
at the endpoints of the contact interval, as shown on the right in
Figure~\ref{fig:forces}. Forces acting between blocks and between
the blocks and the table are said to be \emph{internal} forces.

Each block is also subjected to a downward \emph{gravitational
force} of unit size, acting at its center of mass. As the blocks are
assumed to be of uniform density, the center of mass of a block
whose lower left corner is at $(x,y)$ is at
$(x+\frac{1}{2},y+\frac{h}{2})$.

A rigid body is said to be in \emph{equilibrium} if the sum of the
forces acting on it, and the sum of the \emph{moments} they apply on
it, are both zero. A 2-dimensional rigid body acted upon by~$k$
vertical forces $f_1,f_2,\ldots,f_k$ at $(x_1,y_1),\ldots,(x_k,y_k)$
is in equilibrium if and only if $\sum_{i=1}^k {f}_i=0$ and
$\sum_{i=1}^k x_i f_i=0$. (Note that $f_1,f_2,\ldots,f_k$ are
\emph{scalars} that represent the magnitudes of vertical forces.)

A collection of internal forces acting between the blocks of a
stack, and between the blocks and the table, is said to be a
\emph{balancing} set of forces if the forces in this collection
satisfy the requirements mentioned above (i.e., all the forces are
vertical, they come in opposite pairs, and they act only between
blocks that rest on each other) and if, taking into account the
gravitational forces acting on the blocks, all the blocks are in
equilibrium under this collection of forces. We are now ready for a
formal definition of balance.

\begin{definition}[Balance]\label{def:balance}
A stack of blocks is \emph{balanced} if and only if it admits a
balancing set of forces.
\end{definition}

Static balance problems of the kind considered here are often
\emph{under-determined}, so that the resultants of balancing forces
acting between the blocks are usually not uniquely determined. It
was the consideration by one of us of balance issues that
arise in the game of \emph{Jenga}~\cite{Z02} which stimulated this
current work. The following theorem shows that the balance of a
given stack can be checked efficiently.

\begin{theorem} \label{thm:LP}
The balance of a stack containing $n$ blocks can be decided by
checking the feasibility of a collection of linear equations and
inequalities with $O(n)$ variables and constraints.
\end{theorem}

\begin{proof}
Let $(x_1,y_1),\ldots,(x_n,y_n)$ be the coordinates of the lower
left corners of the blocks in the stack. Let $B_i$, for $1\leqslant i\leqslant
n$, denote the $i$th block of the stack, and let $B_0$ denote
the table. Let $B_i/B_j$, where $0\leqslant i,j\leqslant n$, signify that $B_i$
rests on~$B_j$. If $B_i/B_j$, we let $a_{ij}=\max\{x_i,x_j\}$ and
$b_{ij}=\min\{x_i,x_j\}+1$ be the $x$-coordinates of the endpoints
of the interval of contact between blocks~$i$ and~$j$. (If $j=0$,
then $a_{i0}=x_i$ and $b_{i0}=\min\{x_i+1,0\}$.)

For all $i$ and $j$ such that $B_i/B_j$, we introduce two
variables $f^0_{ij}$ and $f^1_{ij}$ that represent the resultant
forces that act between $B_i$ and $B_j$ at $a_{ij}$ and $b_{ij}$. By
Definition~\ref{def:balance} and the discussion preceding it, the
stack is balanced if and only if there is a feasible solution to the
following set of linear equalities and inequalities:
$$\begin{array}{cl} \displaystyle\sum_{j\,:\,B_i/B_j} (f^0_{ij} + f^1_{ij}) \; -
\sum_{k\,:\, B_k/B_i} (f^0_{ki} + f^1_{ki}) \;=\; 1 \ , & {\rm for} \ 1\leqslant i\leqslant n ;\\
\displaystyle\sum_{j\,:\,B_i/B_j} (a_{ij}f^0_{ij} + b_{ij}f^1_{ij})
\; - \sum_{k\,:\, B_k/B_i} (a_{ki}f^0_{ki} + b_{ki}f^1_{ki}) \;=\;
x_i+\frac{1}{2} \ , & {\rm for} \ 1\leqslant i\leqslant n ;\\[17pt]
\displaystyle f^0_{ij},f^1_{ij}\;\geqslant\; 0 \ , & \mbox{for $i,j$
such that $B_i / B_j$.}
\end{array}$$
The first $2n$ equations require the forces applied on the blocks to
exactly cancel the forces and moments exerted on the blocks by
the gravitational forces. (Note that the table is not required to be
in equilibrium.) The inequalities $f^0_{ij},f^1_{ij}\geqslant 0$, for
every $i$ and $j$ such that~$B_i/B_j$, require the forces applied on
$B_i$ by $B_j$ to point upward. As a unit length block can rest on
at most two other unit length blocks,
the number of variables is at most~$4n$ and the number of
constraints is therefore at most~$6n$. The feasibility of such a
system of linear equations and inequalities can be checked using
\emph{linear programming} techniques. (See, e.g., Schrijver
\cite{S98}.)
\end{proof}

\begin{definition}[Maximum overhang, the function $D(n)$] The
maximum overhang that can be achieved using a balanced stack
comprising $n$ blocks of length~1 is denoted by $D(n)$.
\end{definition}

We now repeat the definitions of the principal block, the support
set, and the balancing set sketched in the introduction.

\begin{definitions}[Principal block, support set, balancing set]
The block of a stack that achieves the maximum overhang is the
\emph{principal block} of the stack. If several blocks achieve the
maximum overhang, the lowest one is chosen. The \emph{support set}
of a stack is defined recursively as follows: the principal block is
in the support set, and if a block is in the support set then any
block on which this block rests is also in the support set. The
\emph{balancing set} consists of any blocks that do not belong to
the support set.
\end{definitions}

The blocks of the support sets of the stacks in
Figures~\ref{fig:2-10}, \ref{fig:11-19}, and \ref{fig:20-30} are
shown in light gray while the blocks in the balancing sets are shown
in dark gray. The purpose of blocks in the support set, as its name
indicates, is to support the principal block. The blocks in the
balancing set, on the other hand, are used to counterbalance the
blocks in the support set.

As already mentioned, there is usually a lot of freedom in the
placement of the blocks of the balancing set. To concentrate on the
more important issue of where to place the blocks of the support set,
it is useful to introduce the notion of loaded stacks.

\begin{definitions}[Loaded stacks, the function $D^*(w)$]
A \emph{loaded stack} consists of a set of blocks with some
\emph{point weights} attached to them. The \emph{weight} of a loaded
stack is the sum of the weights of all the blocks and point weights
that participate in it, where the weight of each block is taken to
be~$1$. A loaded stack is said to be balanced if it admits a
balancing set of forces, as for unloaded stacks, but now also taking
into account the point weights. The maximum overhang that can be
achieved using a balanced loaded stack of weight~$w$ is denoted
by~$D^*(w)$.
\end{definitions}

Clearly $D^*(n)\geqslant D(n)$, as a standard stack is a trivially loaded
stack with no point weights. When drawing loaded stacks, as in
Figure~\ref{fig:p40-100}, we depict point weights as external forces
acting on the blocks of the stack, with the length of the arrow
representing the force proportional to the weight of the point
weight.
(Since forces can be transmitted vertically downwards through any block,
we may assume that point weights are applied only to upper edges
of blocks outside any interval of contact.)

As the next lemma shows, balancing blocks can always be replaced by
point weights, yielding loaded stacks in which all blocks belong to
the support set.

\begin{lemma} For every balanced stack that contains $k$ blocks in its support
set and $n-k$ blocks in its balancing set, there is a balanced loaded
stack composed of $k$ blocks, all in the support set, and additional
point weights of total weight $n-k$ that achieves the same overhang.
\end{lemma}

\begin{proof}
Consider the set of forces exerted on the support set of the stack
by the set of balancing blocks. From the definition of the support
set, no block of the support set can rest on any balancing block,
therefore the effect of the balancing set can be represented by a
set of \emph{downward} vertical forces on the support set, or
equivalently by a finite set of point weights attached to the
support set with the same total weight as the set of balancing
blocks.
\end{proof}

\begin{figure*}[t]
\begin{center}
\includegraphics[height=36mm]{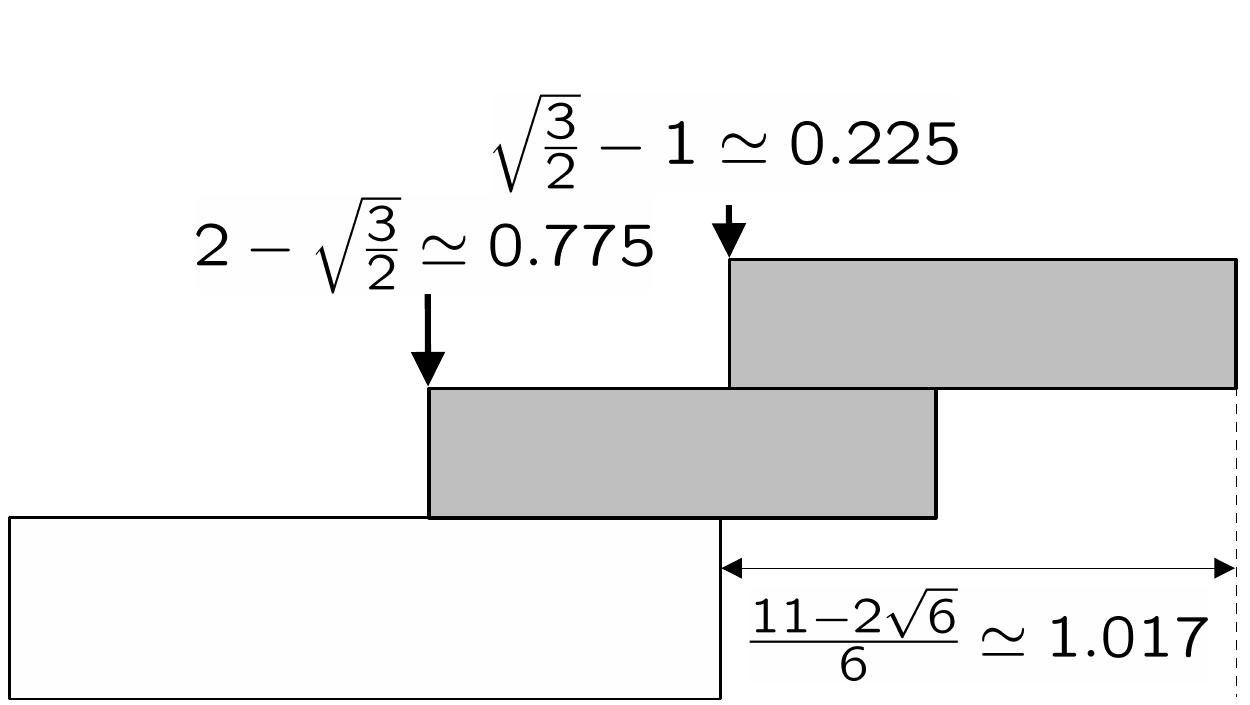} \hspace*{7mm}
{\includegraphics[height=36mm]{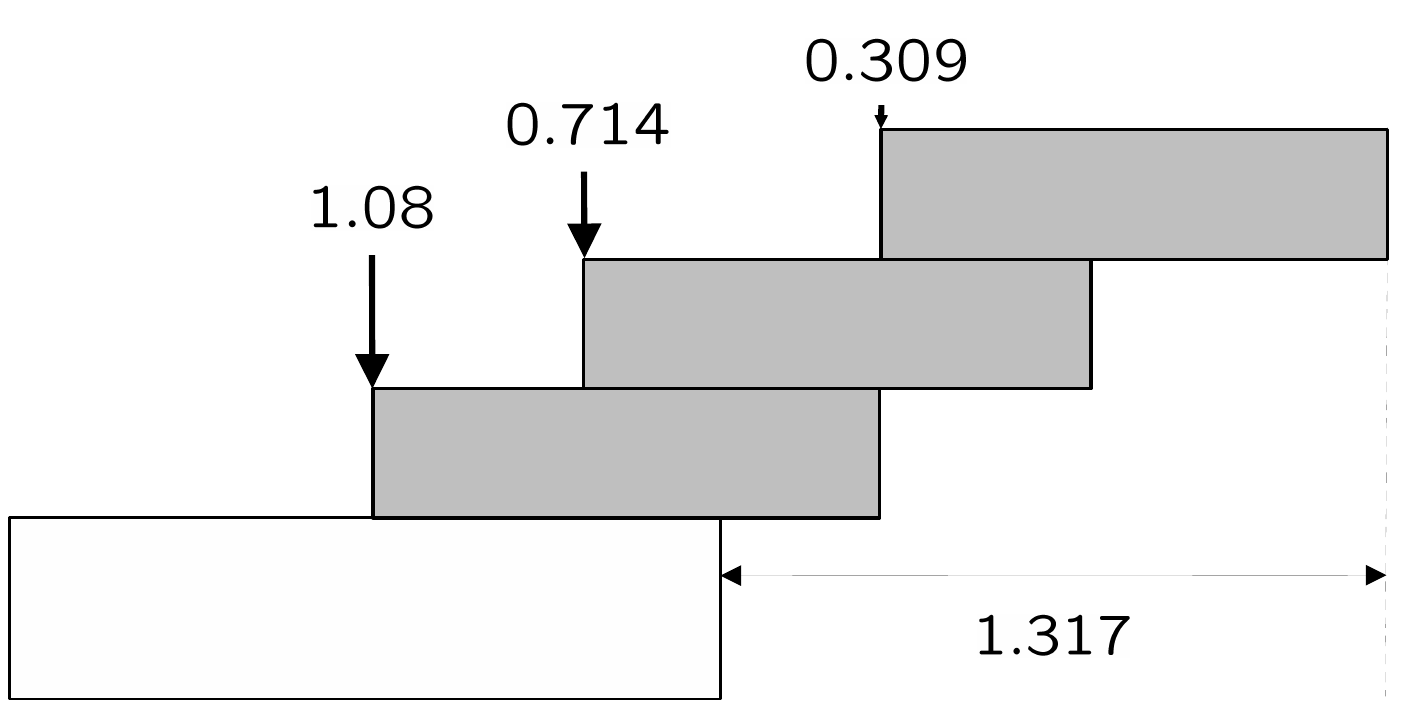}} \caption{Optimal
loaded stacks of weight 3 and 5.} \label{fig:l3}
\end{center}
\end{figure*}

Given a loaded stack of integral weight, it is in many cases
possible to replace the set of point weights by a set of
appropriately placed balancing blocks. In some cases, however, such
a conversion of a loaded stack into a standard stack is not
possible. The optimal loaded stacks of weight~3,~5, and~7
cannot be converted into standard stacks without decreasing
the overhang, as the number of point weights needed is larger than
the number of blocks remaining. (The cases of weights~3 and~5 are
shown in Figure~\ref{fig:l3}.) In particular, we get that
$D^*(3)=\frac{11-2\sqrt{6}}{6}>D(3)=1$. Experiments with optimal
loaded stacks lead us, however, to conjecture that the difference
$D^*(n)-D(n)$ tends to~$0$ as $n$ tends to infinity.

\begin{conjecture}\label{con:D}
$D(n)= D^*(n) - o(1)$.
\end{conjecture}

\section{Spinal stacks} \label{sec:spinal}

In this section we focus on a restricted, but quite natural, class
of stacks which admits a fairly simple analysis.

\begin{definitions}[Spinal stacks, spine]
A stack is \emph{spinal} if its support set has just a single block
at each level. The support set of a spinal stack is referred to as
its \emph{spine}.
\end{definitions}

The optimal stacks with up to 19 blocks, depicted in
Figures~\ref{fig:2-10} and~\ref{fig:11-19}, are spinal. The stacks
of Figure~\ref{fig:20-30} are \emph{not} spinal.
A stack is said to be \emph{monotone} if the $x$-coordinates of the
rightmost blocks in the various levels, starting from the bottom,
form an increasing sequence. It is easy to see that every monotone
stack is spinal.

\begin{definitions}[The functions $S(n)$, $S^*(w)$ and $S^*_k(w)$]
Let $S(n)$ be the maximum overhang achievable using a spinal stack
of size~$n$.
Similarly, let $S^*(w)$ be the maximum overhang achievable using a
loaded spinal stack of weight~$w$, and let $S^*_k(w)$ be the maximum
overhang achievable using a spinal stack of weight~$w$ with
exactly~$k$ blocks in its spine.
\end{definitions}

\begin{figure*}[t]
\begin{center}
\includegraphics[height=80mm]{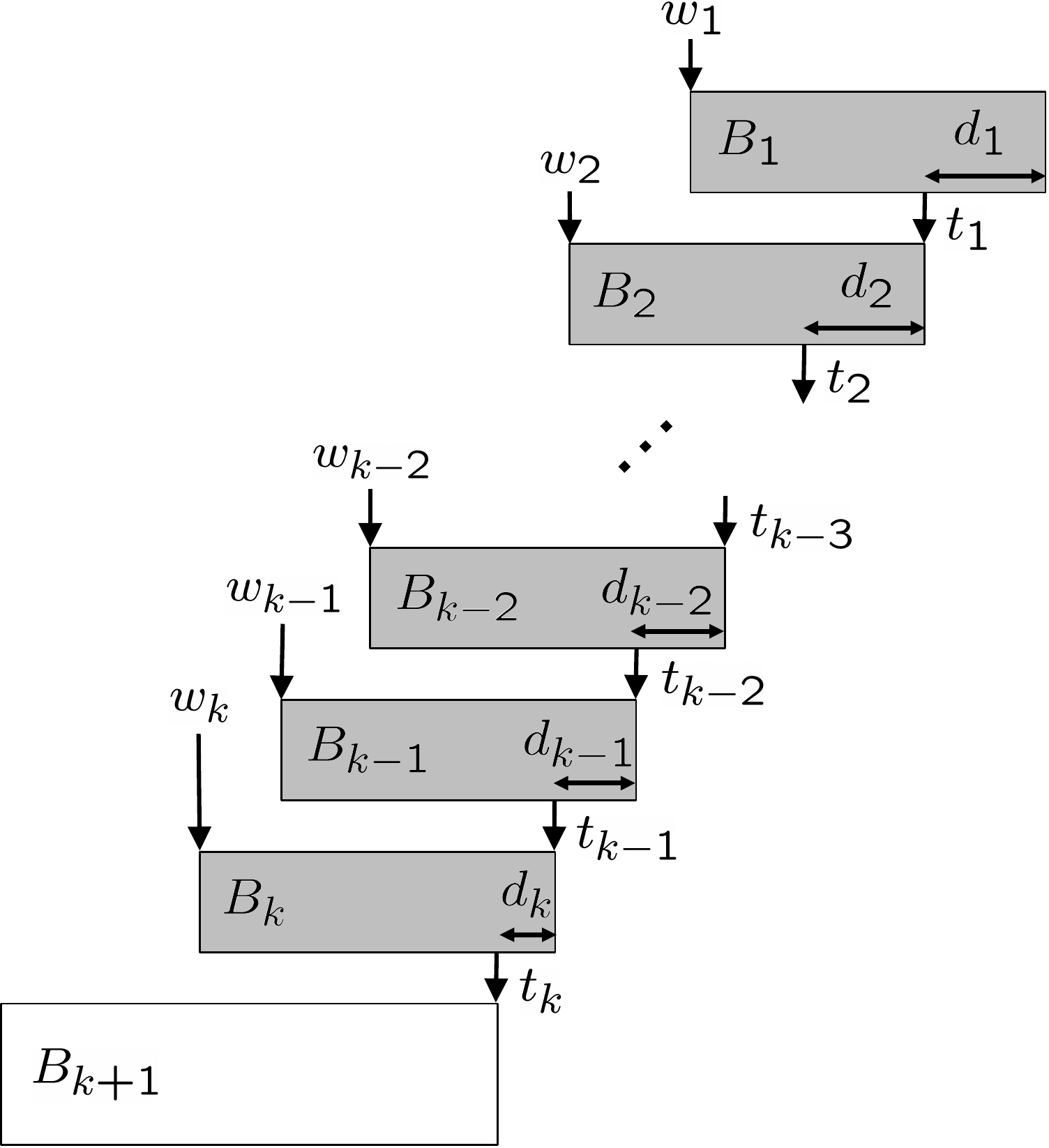}
\caption{A generic loaded spinal stack.} \label{fig:spinal}
\end{center}
\end{figure*}

It is tempting to make the (false) assumption that optimal stacks
are spinal. (As mentioned in the introduction, this assumption is
implicit in \cite{H05}.) The assumption holds, however, only for
$n\leqslant 19$. (See the discussion following Theorem~\ref{thm:lower}.)
As spinal stacks form a very natural class of stacks, it
is still interesting to investigate the maximum overhang achievable
using stacks of this class.

A generic loaded spinal stack with $k$ blocks in its spine is shown
in Figure~\ref{fig:spinal}. We denote the blocks from top to bottom
as $B_1,B_2, \ldots ,B_k$, with $B_1$ being the principal block. We
regard the tabletop as $B_{k+1}$. For $1\leqslant i\leqslant k$, the
weight attached to the left edge of~$B_i$ is denoted by~$w_i$, and
the relative overhang of~$B_i$ beyond~$B_{i+1}$ is denoted by~$d_i$.
We define $t_i = \sum_{j=1}^i (1+w_j)$, the total downward force
exerted upon $B_{i+1}$ by block $B_i$. We also define $t_0=0$. Note
that $t_i=t_{i-1}+w_i+1$, for $1\leqslant i\leqslant k$, and that
$t_k=w=k+\sum_{i=1}^k w_i$, the total weight of the loaded stack.

The assumptions made in Figure~\ref{fig:spinal}, that each
block is supported by a force that acts along the right-hand edge of the
block underneath it and that all point weights are attached to the
left-hand ends of blocks, are justified by the following lemma.

\begin{lemma}\label{lem:L1} In an optimal loaded spinal stack:
$(i)$ Each block is supported by a force acting along the right-hand edge
of the block underneath it. In particular, the stack is monotone.
$(ii)$ All point weights are attached to the left-hand ends of
blocks.
\end{lemma}

\begin{proof} For~$(i)$,
suppose there were some block $B_{i+1}$ ($1\leqslant i\leqslant k$)
where the resultant force exerted on it from $B_i$ does not go
through its right-hand end. If $i<k$  then $B_{i+1}$ could be shifted
some distance to the left and $B_i$ together with all the blocks
above it shifted to the right in such a way that the resultant force
from $B_{i+1}$ on $B_{i+2}$ remains unchanged in position and the
stack is still balanced. In the case of $i=k$ (where $B_{k+1}$ is the
tabletop), the whole stack could be moved to the right. The result
of any such change is a balanced spinal stack with an increased
overhang, a contradiction. As an immediate consequence, we get that
optimal spinal stacks are monotone.

For~$(ii)$, suppose that some block has weights attached other than
at its left-hand end. We may replace all such weights by the same
total weight concentrated at the left end. The result will be to
move the resultant force transmitted to any lower block somewhat to
the left. Since the stack is monotone, this change cannot unbalance
the stack, and indeed would then allow the overhang to be increased
by slightly shifting all blocks to the right; again a contradiction.
\end{proof}

We next note that for any nonnegative point weights
$w_1,w_2,\ldots,w_k\geqslant 0$, there are appropriate positive
displacements $d_1,d_2,\ldots,d_k> 0$ for which the generic spinal
stack of Figure~\ref{fig:spinal} is balanced.

\begin{lemma}\label{lem:L2} A loaded spinal stack with $k$ blocks in its spine that
satisfies the two conditions of Lemma~\ref{lem:L1} is balanced if and only
if $d_i = \frac{w_i+\frac{1}{2}}{t_i} = 1-
\frac{t_{i-1}+\frac{1}{2}}{t_i} $, for $1\leqslant i\leqslant k$.
\end{lemma}

\begin{proof} The lemma is verified using a simple calculation. The
net downward force acting on~$B_i$ is $(w_i+t_{i-1}+1)-t_i=0$, by
the definition of~$t_i$. (Recall that $t_i = \sum_{j=1}^i (1+w_j)$.)
The net moment acting on~$B_i$, computed relative to the right-hand
edge of~$B_i$, is $d_it_i-(\frac{1}{2}+w_i)$, which vanishes if and
only if $d_i=\frac{\frac{1}{2}+w_i}{t_i} = 1-
\frac{t_{i-1}+\frac{1}{2}}{t_i} $, as required.
\end{proof}

Note, in particular, that if $w_i=0$ for $1\leqslant i\leqslant k$, then
$t_i=i$ and $d_i=\frac{1}{2i}$, and we are back to the classic
harmonic stacks.

We can now also justify the claim made in the introduction
concerning the instability of diamond stacks. Consider the spine of
an $m$-diamond. In this case, $d_i=\frac12$ for all~$i$ and so the
balance conditions give the equations $t_i= 2t_{i-1}+1$ for
$1\leqslant i\leqslant m$. As $t_0 = 0$, we have $t_i\geqslant 2^i
-1$ for all~$i$ and hence $t_m\geqslant 2^m -1$. Since $t_m$ is the
total weight of the stack, the number of extra blocks required to be
added for stability is at least $2^m -1-m^2$, which is positive for
$m\geqslant 5$.

Next, we characterize the choice of the weights
$w_1,w_2,\ldots,w_k$, or alternatively of the total loads
$t_1,t_2,\ldots,t_k$, that maximizes the overhang achieved by a
spinal stack of total weight~$w$. (Note that $w_i=t_i-t_{i-1}-1$,
for $1\leqslant i\leqslant k$.)

\begin{lemma} \label{lem:L3}
If a loaded spinal stack with total weight $w$ and with $k$ blocks
in its spine achieves the maximal overhang of $S^*_k(w)$, then for
some $j$ $(1\leqslant j\leqslant k)$ we have $t_i^2=(t_{i-1}+\frac{1}{2})t_{i+1}$, for
$1\leqslant i<j$, and $w_i=0$, for $j< i\leqslant k$.
\end{lemma}

\begin{proof} Let $w_1,w_2,\ldots,w_k$ be the point weights attached
to the blocks of an optimal spinal stack with overhang $S^*_k(w)$.
For some $i$ satisfying $1\leqslant i<k$ and a small $x$, consider the stack obtained by
increasing the point weight at the left-hand end of block~$B_i$ from $w_i$
to $w_i+x$, and decreasing the point weight on~$B_{i+1}$ from $w_{i+1}$
to $w_{i+1}-x$, assuming
that $w_{i+1}\geqslant x$. Note that this small perturbation does not
change the total weight of the stack.
The overhang of the perturbed stack is
$$V(x) \;=\; \left(1-\frac{t_{i-1}+\frac12}{t_i+x}\right) +
\frac{w_{i+1}-x+\frac12}{t_{i+1}} + \sum_{j\ne i,i+1}
\frac{w_j+\frac12}{t_j}\;.$$
The first two terms in the expression
above are the new displacements $d_i(x)$ and $d_{i+1}(x)$. Note that
all other displacements are unchanged. Differentiating $V(x)$ we get
$$V'(x) \;=\; \frac{t_{i-1}+\frac12}{(t_i+x)^2} -
\frac{1}{t_{i+1}}\quad{\rm\ and}\quad V'(0) \;=\;
\frac{t_{i-1}+\frac12}{t_i^2} - \frac{1}{t_{i+1}}\;.$$ If $w_i=0$
while $w_{i+1}>0$, then $t_{i-1}=t_i-1$ and $t_{i+1}>t_i+1$, which
in conjunction with $t_i\geqslant 1$ implies that $V'(0)>0$, contradicting
the optimality of the stack. Thus, if in an optimal stack we have
$w_i=0$, then also $w_{i+1}=w_{i+2}=\ldots=w_k=0$. If
$w_i,w_{i+1}>0$, then we must have $V'(0)=0$, or equivalently
$t_i^2=(t_{i-1}+\frac{1}{2})t_{i+1}$, as claimed.
\end{proof}

\begin{figure}[t]
\begin{center}
\includegraphics[width=8cm]{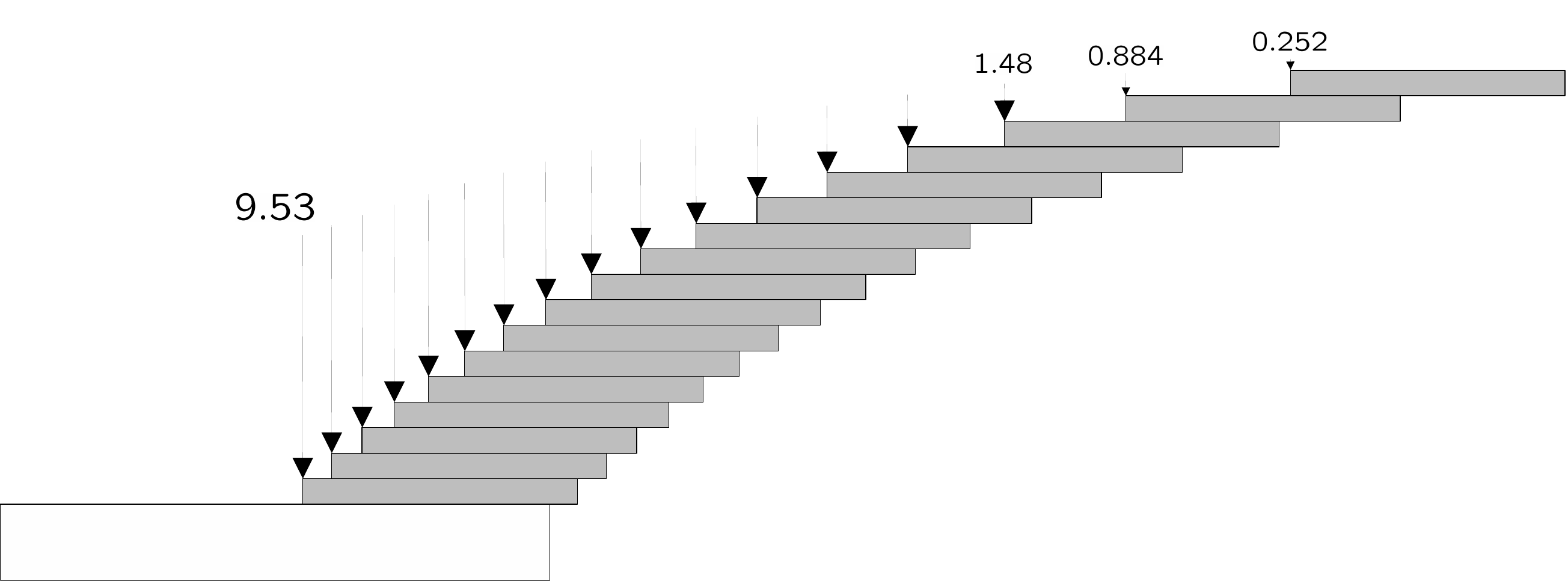}
\includegraphics[width=8cm]{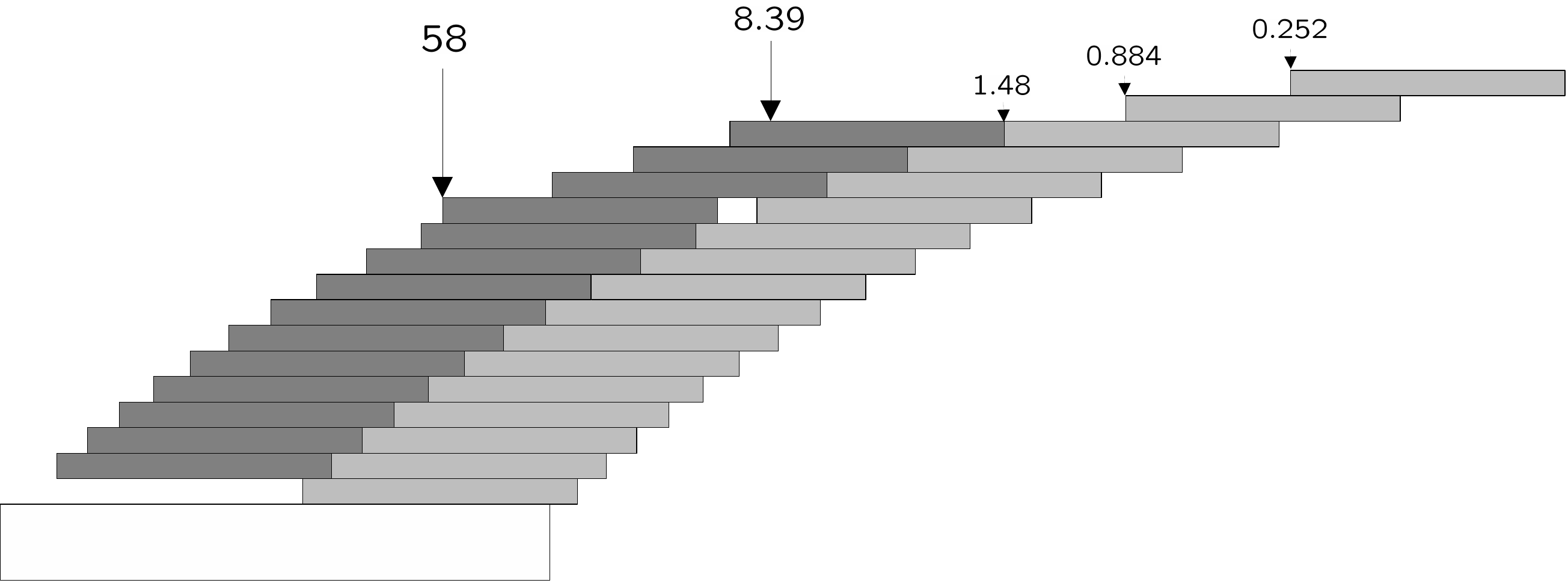}\\[0.8cm]
\includegraphics[width=8cm]{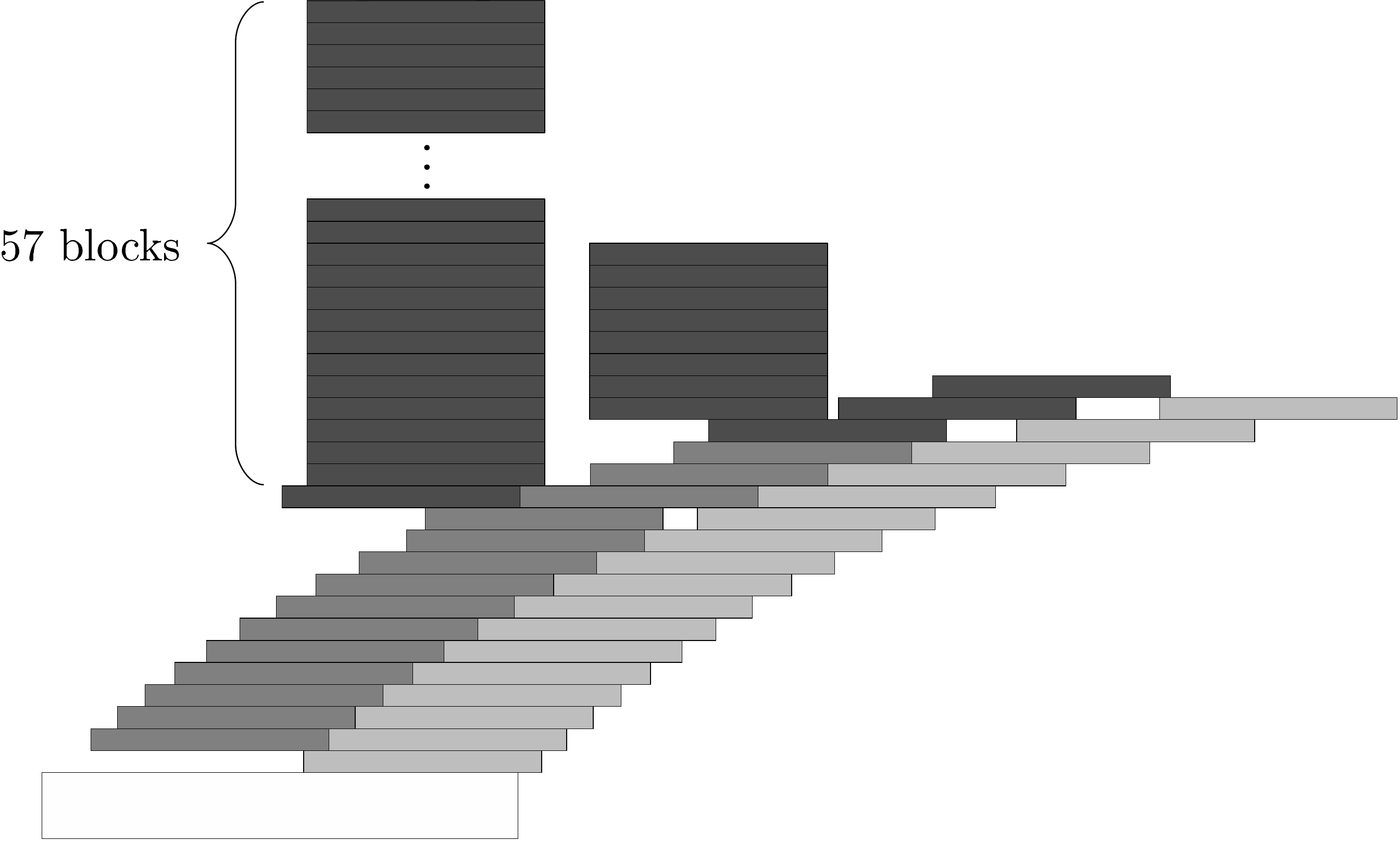}
\caption{Optimal loaded spinal stack of weight 100 (top left), with
shield added (top right) and with a complete balancing set added
(bottom).} \label{fig:s100}
\end{center}\vspace{-5mm}
\end{figure}

The optimality equations given in Lemma~\ref{lem:L3} can be solved
numerically to obtain the values of~$S^*_k(w)$ for specific
values of $w$ and $k$. The value of $S^*(w)$ is then found by
optimizing over~$k$. The optimal loaded spinal stacks of weight~3
and 5, which also turn out to be the optimal loaded stacks of these
weights, are shown in Figure~\ref{fig:l3}. The optimality equations
of Lemma~\ref{lem:L3} were also used to compute the spines of the
optimal stacks with up to 19 blocks shown in Figures~\ref{fig:2-10}
and~\ref{fig:11-19}. The spines of the stacks with~3 and~5 blocks
were obtained by adding the requirement that no point weight be
attached to the topmost block of the spine.
A somewhat larger example is given on the top left of
Figure~\ref{fig:s100} where the optimal loaded spinal stack of
weight~100 is shown. It is interesting to note that the
point weights in optimal spinal stacks form an almost arithmetical
progression. This observation is used in the proof of
Theorem~\ref{thm:spinal-lower}.

Numerical experiments suggest that for every $w\geqslant 1$, all the point
weights in the spinal stacks with overhang~$S^*(w)$ are nonzero.
There are, however, non-optimal values of~$k$ for which some of the
bottom blocks in the stack that achieves an overhang of $S^*_k(w)$
have no point weights attached to them. We next show, without
explicitly using the optimality conditions of Lemma~\ref{lem:L3},
that $S^*(w)=\ln w + \Theta(1)$.

\begin{theorem} \label{thm:spinal-upper}
$S^*(w) < \ln w +1$.
\end{theorem}

\begin{proof}
For fixed total weight $w=t_k$ and fixed $k$, the largest possible
overhang $S^*_k(w)= \sum_{i=1}^k d_i$ is attained when the
conditions of Lemmas~\ref{lem:L1} and~\ref{lem:L2}
(and~\ref{lem:L3}) hold. Thus, as $t_{0}=0$,
$$\sum_{i=1}^k d_i \;=\; \sum_{i=1}^{k}\left(1- \frac{t_{i-1}+\frac{1}{2}}{t_i}
\right) \;<\; k - \sum_{i=2}^{k} \frac{t_{i-1}}{t_i}\;.$$
Putting $x_i=\frac{t_{i-1}}{t_i}$, we see that
$$S^*_k(w) \;<\; k - \sum_{i=2}^{k} x_i
   \;\; {\rm\ and\;\; } \prod_{i=2}^{k}x_i \;=\; \frac{t_1}{t_k} \;\geqslant\; \frac{1}{w}\;.$$
The minimum sum for a finite set of positive real numbers with fixed
product is attained when the numbers are equal, hence
$$S^*_k(w) \;<\; k - (k-1)w^{-\frac{1}{k-1}}\;.$$
Let $z=\frac{k-1}{\ln w}$, so that $k-1=z\ln w$ and
$w^{-\frac{1}{k-1}} = e^{-1/z}$. Then
$$S^*_k(w) \;<\; 1 + z\ln w(1-e^{-1/z}) < 1 + \ln w\;,$$
as $z(1-e^{-1/z})\leqslant 1$, for every $z>0$.
\end{proof}

\begin{corollary} \label{cor:spinal-upper}
$S(n) <  \ln n + 1$.
\end{corollary}

We can now describe a construction of loaded spinal stacks
which achieves an overhang agreeing asymptotically with the
upper bound proved in Theorem~\ref{thm:spinal-upper}.

\begin{theorem}\label{thm:spinal-lower}
$S^*(w) > \ln w - 1.313$.
\end{theorem}

\begin{proof}
We construct a spine with $k=\lfloor \sqrt{w}\rfloor$ blocks in it
with $w_i=2(i-1)$, for $1\leqslant i\leqslant k$. It follows easily by
induction that $t_i=i^2$, for $1\leqslant i\leqslant k$. In particular, the
total weight of the stack is $t_k=k^2\leqslant w$, as required. By
Lemma~\ref{lem:L2}, we get that
$$d_i \;=\; \frac{w_i+\frac{1}{2}}{t_i} \;=\;
\frac{2(i-1)+\frac{1}{2}}{i^2} \;=\; \frac{2}{i} -
\frac{3}{2i^2}\;.$$ Thus,
$$S^*(w) \;\geqslant\; \sum_{i=1}^k d_i \;=\; 2\sum_{i=1}^k\frac{1}{i}
-\frac{3}{2}\sum_{i=1}^k\frac{1}{i^2}\;=\;
2H_{\lfloor\sqrt{w}\rfloor}-\frac{3}{2}\sum_{i=1}^{\lfloor\sqrt{w}\rfloor}\frac{1}{i^2}\;>\;
\ln w +2\gamma - \frac{\pi^2}{4} \;>\; \ln w - 1.313\;.$$ In the
above inequality, $\gamma\simeq 0.5772156$ is Euler's gamma.
\end{proof}

\begin{figure*}[t]
\begin{center}
\includegraphics[height=80mm]{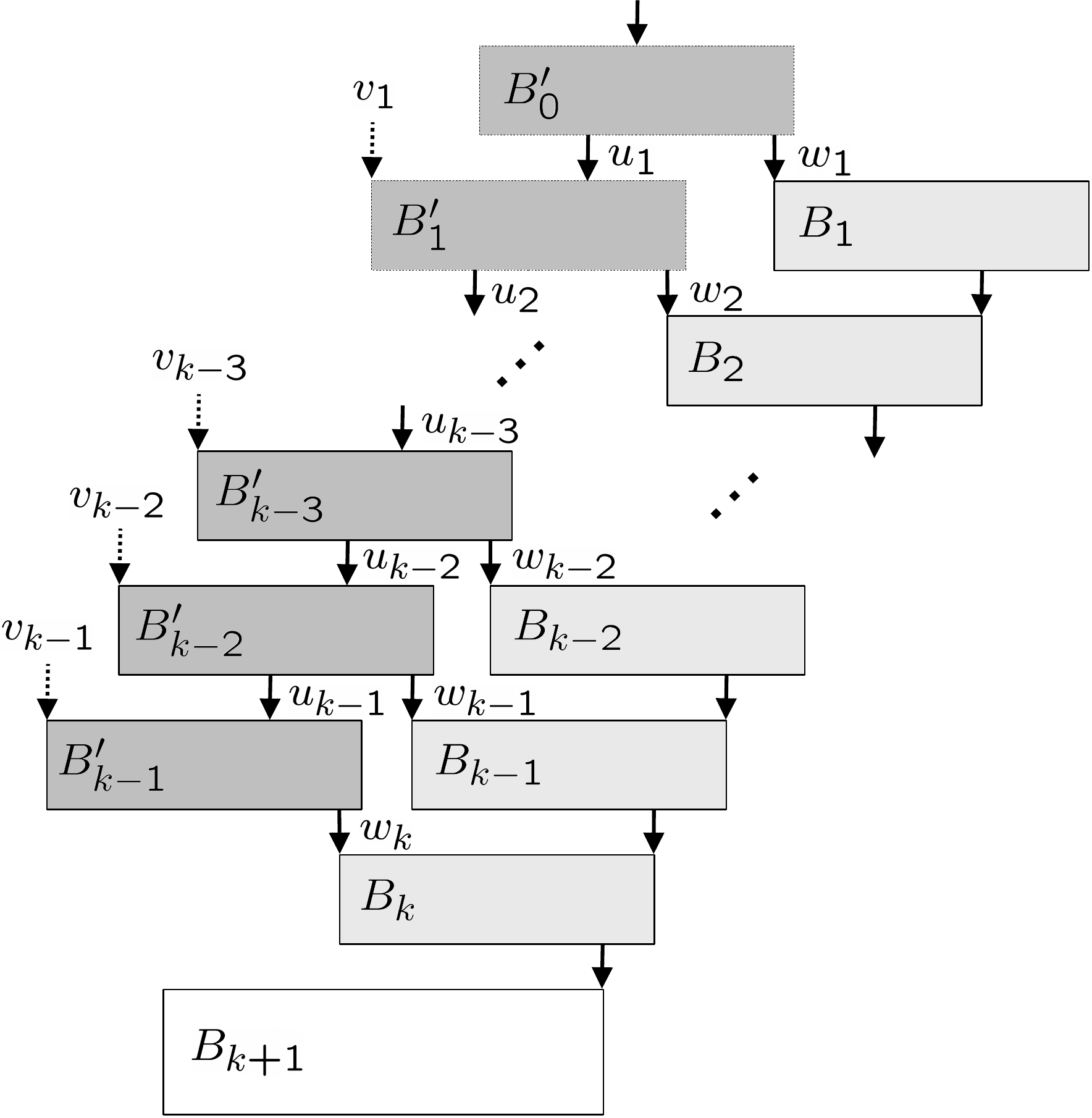}
\caption{A spinal stack with a shield.} \label{fig:shadow}
\end{center}
\end{figure*}

We next discuss a technique that can be used to convert loaded
spinal stacks into standard stacks. This is of course done by
constructing balancing sets that apply the required forces on the
left-hand edges of the spine blocks. The first step is the placement of
\emph{shield} blocks on top of the spine blocks, as shown in
Figure~\ref{fig:shadow}. We let $B'_i$, for $0\leqslant i\leqslant k-1$, be the
shield block placed on top of spine block $B_{i+1}$ and alongside
spine block $B_i$ for $i>0$. We let $y_i$ be the $x$-coordinate of the left
edge of~$B'_i$, for $1\leqslant i\leqslant k-1$. Note that $x_{i+1}-1 < y_i\leqslant
x_i-1$, where $x_i$ is the $x$-coordinate of the left edge of~$B_i$.

Shield block~$B'_i$ applies a downward force of $w_{i+1}$ on
$B_{i+1}$. The force is applied at $x_{i+1}$, i.e., at the left
edge of $B_{i+1}$. Block~$B'_i$ also applies a downward force of
$u_{i+1}$ on $B'_{i+1}$ at $z_{i+1}$, where $y_i\leqslant z_{i+1}\leqslant
y_{i+1}+1$. Similarly, block $B'_{i-1}$ applies a downward force of
$u_i$ on $B'_i$ at $z_i$. Finally a downward external force of $v_i$
is applied on the left edge of $B'_i$. The goal of the shield blocks
is to aggregate the forces that should be applied on the spine
blocks and to replace them by a set of fewer \emph{integral} forces
that are to be applied on the shield blocks. We will therefore place
the shield blocks and choose the forces $u_i$ and their positions in
such a way that most of the $v_i$ will be~$0$. (This is why we use
dashed arrows to represent the $v_i$ forces in Figure~\ref{fig:shadow}.)

The shield blocks are in equilibrium provided that the following
balance conditions are satisfied: $$\begin{array}{c} u_i+v_i+1 \;=\;
u_{i+1}+w_{i+1}\;,\\
z_iu_i+y_iv_i+(y_i+\frac12) \;=\;
z_{i+1}u_{i+1}+x_{i+1}w_{i+1}\;,\\
\end{array}$$
for $1\leqslant i\leqslant k-1$. (We define $u_k=0$.) It is easy to see
that if $u_{i+1},w_{i+1},x_{i+1},y_{i+1},$ and $z_{i+1}$ are set,
then any choice of $v_i$ uniquely determines $u_i$ and $z_i$. The
choice is feasible if $u_i,v_i\geqslant 0$ and $y_{i-1}\leqslant z_i$.

In our constructions, we used the following heuristic to place the
shield blocks and specify the forces between them. We start placing
the shield blocks from the bottom up. In most cases, we choose
$y_i=x_i-1$ and $v_i=0$, i.e., $B'_i$ is adjacent to $B_i$ and no
external force is applied to it. Eventually, however, it may happen
that $z_{i+1}<x_{i}-1$, which makes it impossible to place~$B'_{i}$
adjacent to~$B_{i}$ and still apply the force~$u_{i+1}$ down
on~$B_{i+1}$ at~$z_{i+1}$. In that case we choose $y_i=z_{i+1}$. A
more significant event, that usually occurs soon after the previous
event, is when $z_{i+1}\leqslant x_{i+1}-1$, in which case no placement of
$B'_i$ allows it to apply the forces $u_{i+1}$ and $v_{i+1}$ on
$B'_{i+1}$ and $B_{i+1}$ at the required positions, as they are at
least a unit distance apart. In this case, we introduce a nonzero,
integral, external force~$v_{i+1}$ as follows. We let
$v_{i+1}=\lfloor (1-z_{i+1}+y_{i+1})u_{i+1} \rfloor$ and then
recompute $u_{i+1}$ and $z_{i+1}$. It is easy to check that
$u_{i+1},v_{i+1}\geqslant 0$ and that $y_i\leqslant z_{i+1}\leqslant y_{i+1} +1$. If we
now have $z_{i+1}>x_{i+1}-1$, then the process can continue. Otherwise
we stop. In our experience, we were always able to use this process
to place all the shield blocks, except for a very few top ones.
The~$v_i$ forces left behind tend to be few and far apart. When this
process is applied, for example, on the optimal loaded spinal stack
of weight~$100$, only one such external force is needed, as shown in
the second diagram of Figure~\ref{fig:s100}.

The nonzero $v_i$'s can be easily realized by erecting appropriate
towers, as shown at the bottom of Figure~\ref{fig:s100}. The top
part of the balancing set is then designed by solving a small linear
program. We omit the fairly straightforward details.

The overhang achieved by the spinal stack shown at the bottom of
Figure~\ref{fig:s100} is about $3.6979$, which is a considerable
improvement on the $2.5937$ overhang of a 100-block harmonic stack,
but is also substantially less than the $4.23897$ overhang of the
non-spinal loaded stack of weight 100 given in
Figure~\ref{fig:p40-100}.

Using the heuristic described above we were able to fit appropriate
balancing sets for all optimal loaded spinal stacks of integer
weight~$n$, for every $n\leqslant 1000$, with the exception of $n=3,5,7$.
We conjecture that the process succeeds for every $n\ne 3,5,7$.

\begin{conjecture}
$S(n)= S^*(n)$ for $n\neq 3,5,{\rm or\ }7$.
\end{conjecture}

\begin{figure}[t]
\begin{center}
\includegraphics[width=85mm]{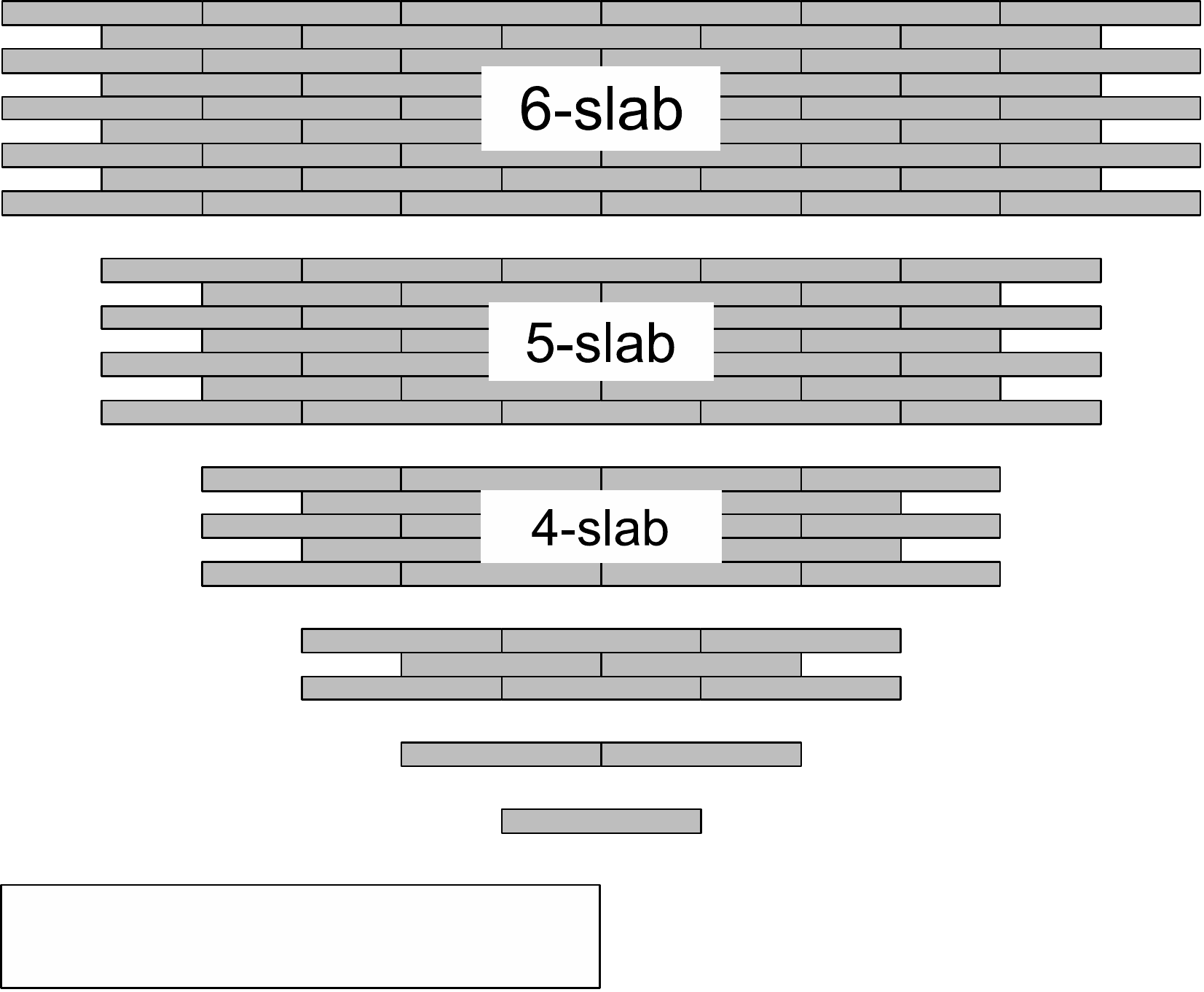}
\caption{A $6$-stack composed of $r$-slabs, for $r=2,3,\ldots,6$,
and an additional block.} \label{fig:slabs2-6}
\end{center}\vspace{-5mm}
\end{figure}

\section{Parabolic stacks} \label{sec:parabolic}

We now give a simple explicit construction of $n$-block stacks with
an overhang of about $(3n/16)^{1/3}$, an \emph{exponential}
improvement over the $O(\log n)$ overhang achievable using spinal
stacks in general and the harmonic stacks in particular. Though the
stacks of this sequence are not optimal (see the empirical results
of the next section), they are within a constant factor of
optimality, as will be shown in a subsequent paper \cite{PPTWZ07}.

The stacks constructed in this section are what we term
\emph{brick-wall} stacks. The blocks in each row are contiguous, and
each is centered over the ends of blocks in the row beneath. This
resembles the simple ``stretcher-bond'' pattern in real-life
bricklaying.
Overall the stacks have a symmetric roughly parabolic shape, hence
the name, with vertical axis at the table edge and a brick-wall
structure. An illustration of a $111$-block \emph{parabolic}
$6$-stack with overhang~$3$ was given in Figure~\ref{fig:6stack}.

\begin{figure}[t]
\begin{center}
\vspace*{0.3cm}
\includegraphics[scale=0.6]{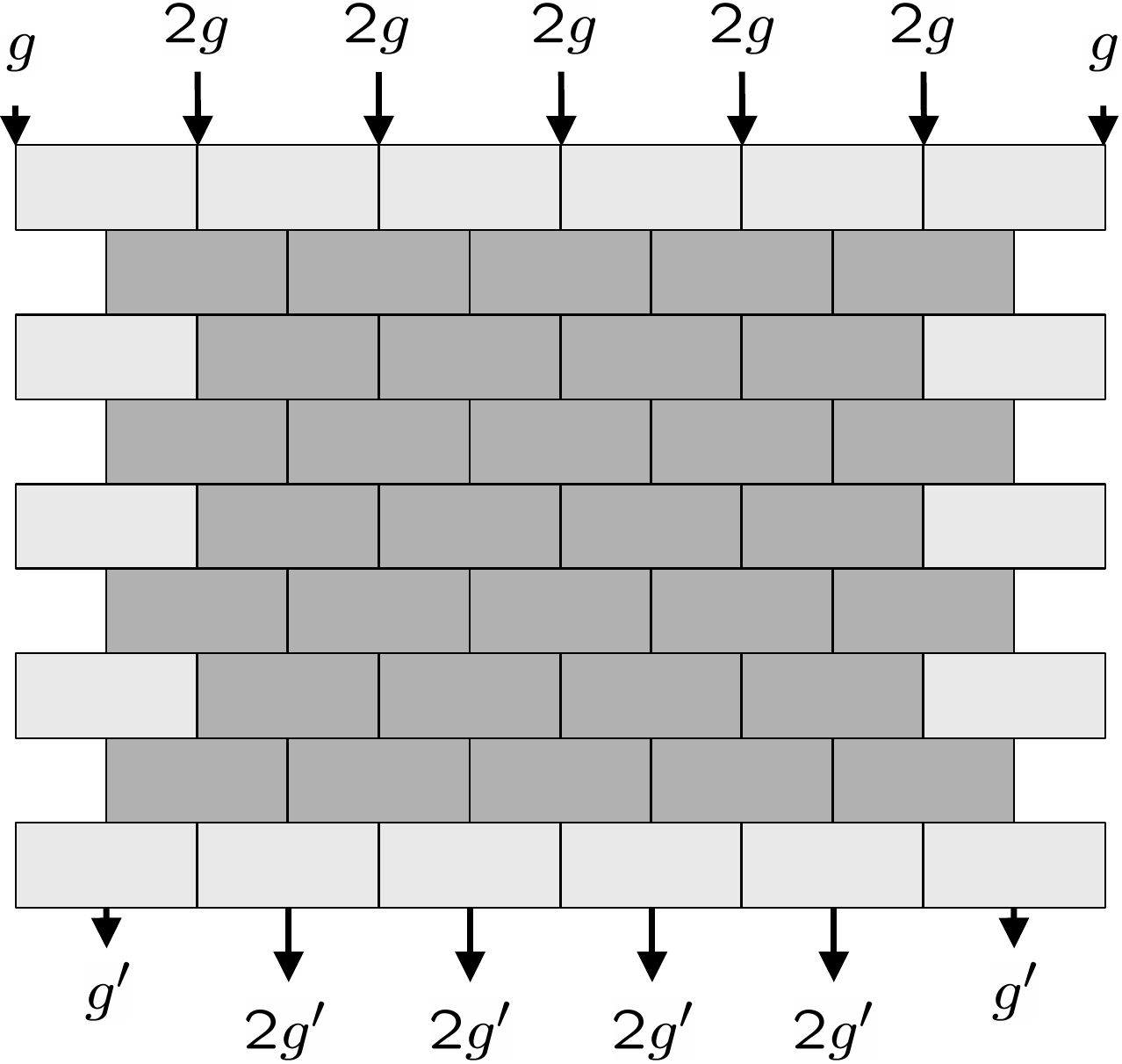}
\caption{A $6$-slab with a grey $5$-slab contained in it.}
\label{fig:slab6}
\end{center}
\end{figure}

An \emph{$r$-row} is a row of $r$ adjacent blocks, symmetrically
placed with respect to $x=0$. An \emph{$r$-slab}, for $r\geqslant 2$, has height $2r-3$
and consists of alternating $r$-rows and $(r-1)$-rows, starting and
finishing with $r$-rows. An $r$-slab therefore contains
$r(r-1)+(r-1)(r-2)=2(r-1)^2$ blocks.  Figure~\ref{fig:slabs2-6} shows
$r$-slabs, for $r=2,3,\ldots,6$. A \emph{parabolic $d$-stack}, or
just \emph{$d$-stack}, for short, is a $d$-slab on a $(d-1)$-slab on
\ldots\ on a 2-slab on a single block. The slabs shown
in Figure~\ref{fig:slabs2-6} thus compose a $6$-stack.

\begin{lemma} A parabolic $d$-stack contains $\frac{d(d-1)(2d-1)}{3} + 1$
blocks and, if balanced, has an overhang of~$\frac{d}{2}$.
\end{lemma}

\begin{proof} The number of blocks contained in a $d$-stack is
$1+\sum_{r=2}^d 2(r-1)^2 = 1 + \frac{d(d-1)(2d-1)}{3}$. The overhang
achieved, if the stack is balanced, is half the width of the top row,
i.e., $\frac{d}{2}$.
\end{proof}

In preparation for proving the balance of parabolic stacks, we
show in the next lemma that a slab can concentrate a set of forces
acting on its top together with the weights of its own blocks down
into a narrower set of forces acting on the row below it. The lemma
is illustrated in Figure~\ref{fig:slab6}.

\begin{lemma} \label{lem:slab}
For any $g\geqslant 0$, an $r$-slab with forces of $g,2g,2g,\ldots
,2g,g$ acting downwards onto its top row at positions
$-\frac{r}{2},-\frac{r-2}{2},-\frac{r-4}{2}, \dots
,\frac{r-2}{2},\frac{r}{2}$, respectively, can be stabilized by
applying a set of upward forces $g',2g',2g',\ldots ,2g',g'$, where
$g'=\frac{r}{r-1}g+r-1$, on its bottom row at positions
$-\frac{r-1}{2},-\frac{r-3}{2}, \dots ,\frac{r-3}{2},\frac{r-1}{2}$,
respectively.
\end{lemma}

\begin{proof}
The proof is by induction on $r$. For $r=2$, a $2$-slab is just a
$2$-row, which is clearly balanced with downward forces of $g,2g,g$ at
$-1,0,1$ and upward forces of $2g+1,2g+1$ at
$-\frac{1}{2},\frac{1}{2}$, when half of the downward force $2g$
acting at $x=0$ is applied on the right-hand edge of the left block
and the other half applied on the left-hand edge of the right block.

For the induction step, we first observe that for any $r\geqslant 2$ an
$(r+1)$-slab can be regarded as an $r$-slab with an $(r+1)$-row
added above and below and with an extra block added at each end of
the $r-2$ rows of length $r-1$ of the $r$-slab.  The 5-slab
(shaded) contained in a 6-slab together with the added blocks is
shown in Figure~\ref{fig:slab6}.

\begin{figure}[t]
\begin{center}
\includegraphics[scale=0.7]{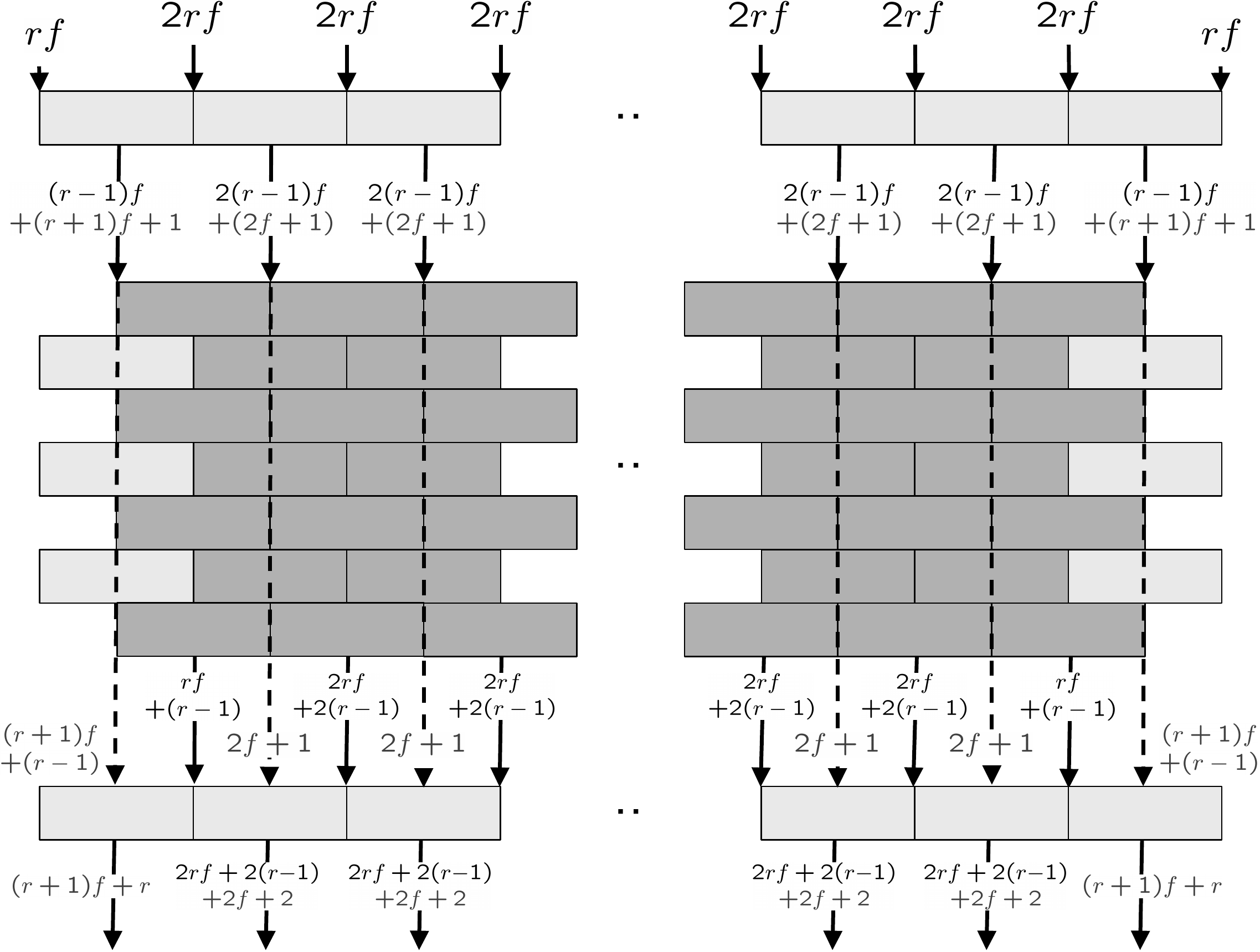}
\caption{The proof of Lemma~\ref{lem:slab}.}
\label{fig:slab-induction}
\end{center}
\end{figure}

Suppose the statement of the lemma holds for $r$-slabs and
consider an $(r+1)$-slab with the supposed forces acting on its top row.
Let $f=g/r$, so that $g=rf$. As in the basis of the induction, the
top row can be balanced by $r+1$ equal forces of $2rf+1$ from below (the 1 is
for the weight of the blocks in the top row) acting at positions
$-\frac{r}{2},-\frac{r-2}{2},\ldots ,\frac{r-2}{2},\frac{r}{2}$. As
$$2rf+1 \;=\; (r-1)f+((r+1)f+1) \;=\; 2(r-1)f+2f+1\;,$$ we can
express this constant sequence of $r+1$ forces as the sum of the
following two force sequences: $$\begin{array}{ccccccccccccc} (r-1)f
&,& 2(r-1)f &,& 2(r-1)f &,& \ldots &,&
2(r-1)f &,& 2(r-1)f &,& (r-1)f\\
(r+1)f+1 &,& 2f+1 &,& 2f+1 &,& \ldots &,& 2f+1 &,& 2f+1 &,& (r+1)f+1\\
\end{array}$$
The forces in the first sequence can be regarded as acting on the
$r$-slab contained in the $(r+1)$-slab, which then, by the induction
hypothesis, yield downward forces on the bottom row of
$$rf+r-1\;\;,\;\;2rf+2(r-1)\;\;,\;\;\ldots \;\;,\;\;2rf+2(r-1)\;\;,\;\;rf+r-1$$ at positions
$-\frac{r-1}{2},-\frac{r-3}{2},\dots\ ,\frac{r-3}{2},\frac{r-1}{2}$.

The forces of the second sequence, together with the weights of the
outermost blocks of the $(r+1)$-rows, are passed straight down
through the rigid structure of the $r$-slab to the bottom row.
The combined forces acting down on the bottom row are now
$$(r+1)f+r{-}1\;,\;rf+r{-}1\;,\;2f+1\;,\;2rf+2(r{-}1)\;,\;2f+1\;,\;\ldots
\;,\;2f+1\;,\;rf+r{-}1\;,\;(r+1)f+r{-}1$$
at positions $-\frac{r}{2},-\frac{r-1}{2},\dots\
,\frac{r-1}{2},\frac{r}{2}$. The bottom row is in equilibrium when
the sequence of upward forces
$$(r+1)f+r\;\;,\;\;2(r+1)f+2r\;\;,\;\;2(r+1)f+2r\;\;,\;\;\ldots
\;\;,\;\;2(r+1)f+2r\;\;,\;\;2(r+1)f+2r\;\;,\;\;(r+1)f+r$$ is applied
on the bottom row at positions $-\frac{r}{2},-\frac{r-2}{2},\dots\
,\frac{r-2}{2},\frac{r}{2}$,
as required.
\end{proof}

\begin{theorem} \label{thm:construction}
For any $d\geqslant 2$, a parabolic $d$-stack is balanced, contains
$\frac{d(d-1)(2d-1)}{3} + 1$ blocks, and has an overhang of
$\frac{d}{2}$.
\end{theorem}

\begin{proof} The balance of a parabolic $d$-stack follows by a
repeated application of Lemma~\ref{lem:slab}.
For $2\leq r\leq d$, let $g(r)$ denote the value of $g$
in Lemma~\ref{lem:slab} for the $r$-slab in the $d$-stack.
Although the
argument does not rely on the specific values that $g(r)$ assumes,
it can be verified that $g(r)=\frac{1}{r}\sum_{i=r}^{d-1} i^2$. Note
that $g(d)=0$, as no downward forces are exerted on the top row of
the $d$-slab, which is also the top row of the $d$-stack, and that
$g(r-1)=\frac{r}{r-1}g(r)+r-1$, as required by Lemma~\ref{lem:slab}.
\end{proof}

\begin{theorem} \label{thm:lower}
$\D(n)\geqslant (\frac{3n}{16})^{1/3} - \frac14$ for all $n$.
\end{theorem}

\begin{proof}
Choose $d$ so that
$\frac{(d-1)d(2d-1)}{3}+1 \leqslant n \leqslant
\frac{d(d+1)(2d+1)}{3}$. Then
Theorem~\ref{thm:construction} shows that a $d$-stack yields an
overhang of $d/2$ and can be constructed using $n$ or fewer blocks.
Any extra blocks can be just placed in a vertical pile in the center
on top of the stack without disturbing balance (or arbitrarily
scattered on the table). Hence
$$n < \frac{2(d+\frac{1}{2})^3}{3}\ {\rm\quad and\ so\quad }
\ \D(n) \geqslant d/2> \left(\frac{3n}{16}\right)^{1/3}-\frac{1}{4}\;.
\vspace*{-5pt}$$
\end{proof}

In Section~\ref{sec:spinal} we claimed that optimal stacks are spinal
\emph{only} for $n\leqslant 19$. We can justify this claim for $n\leqslant 30$
by exhaustive search, while comparison of the lower bound from
Theorem~\ref{thm:lower} with the upper bound
of $S(n) < 1+\ln n$ from Corollary~\ref{cor:spinal-upper} deals with the
range $n\geqslant 5000$. The intermediate values of $n$ can be covered
by combining a few explicit constructions, such as the stack shown in
Figure~\ref{fig:vase95}, with numerical bounds using Lemma~\ref{lem:L3}.

\begin{figure}[t]
\begin{center}
\includegraphics[scale=0.3]{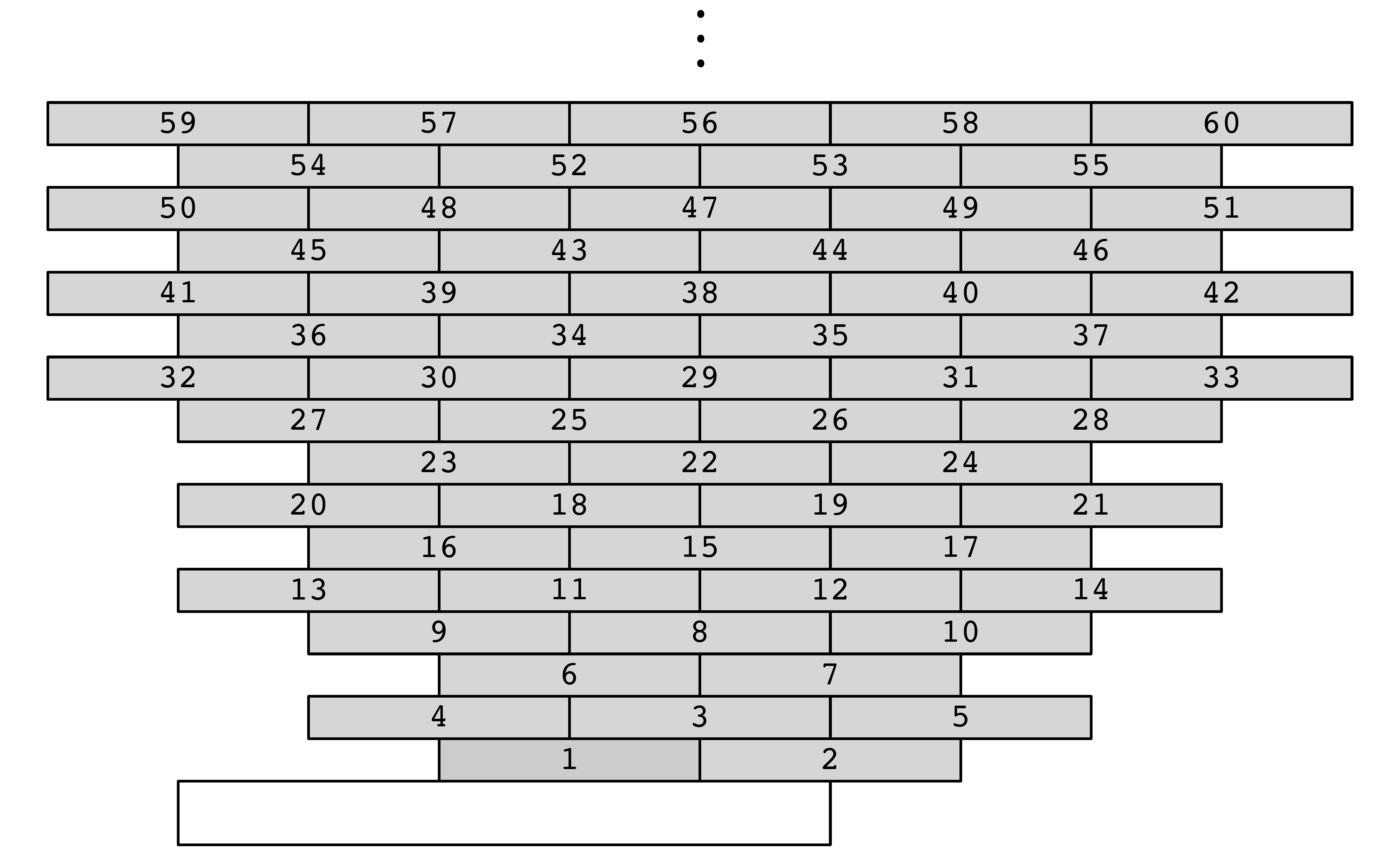}
\caption{Incremental block-by-block construction of modified
parabolic stacks.} \label{fig:modified}
\end{center}
\end{figure}

Can parabolic $d$-stacks be built incrementally by laying one brick
at a time? The answer is no, as the bottom three rows of a parabolic
stack form an unbalanced inverted 3-triangle. The inverted
3-triangle remains unbalanced when the first block of the fourth row
is laid down. Furthermore, the bottom six rows, on their own, are
also not balanced. These, however, are the only obstacles to an
incremental row-by-row and block-by-block construction of parabolic
stacks and they can be overcome by the modified parabolic stacks
shown in Figure~\ref{fig:modified}. We simply omit the lowest block
and move the whole stack half a block length to the left. The bricks
can now be laid row by row, going in each row from the center
outward, alternating between the left and right sides, with the left
side, which is over the table, taking precedence. The numbers in
Figure~\ref{fig:modified} indicate the order in which the blocks are
laid. Thus, unlike with harmonic stacks, it is possible to construct
an arbitrarily large overhang using sufficiently many blocks,
\emph{without} knowing the desired overhang in advance.

\section{General stacks} \label{sec:general}

We saw in Section~\ref{sec:prelim} that the problem of checking
whether a given stack is balanced reduces to checking the feasibility
of a system of linear equations and inequalities. Similarly, the
minimum total weight of the point weights that are needed to
stabilize a given loaded stack can be found by solving a linear
program.

Finding a stack with a given number of blocks, or a loaded stack
with a given total weight, that achieves maximum overhang seems,
however, to be a much harder computational task. To do so, one
should, at least in principle, consider all possible combinatorial
stack structures and for each of them find an optimal placement
of the blocks. The \emph{combinatorial structure} of a stack
specifies the contacts between the blocks of the stack, i.e., which
blocks rest on which, and in what order (from left to right), and
which rest on the table.

The problem of finding a (loaded) stack with a given combinatorial
structure with maximum overhang is again not an easy problem. As
both the forces and their locations are now unknowns, the problem is
not linear, but rather a constrained quadratic programming problem.
Though there is no general algorithm for efficiently finding the
global optimum of such constrained quadratic programs, such
problems can still be solved in practice using nonlinear
optimization techniques.

For stacks with a small number of blocks, we enumerated all
possible combinatorial stack structures and numerically optimized
each of them. For larger numbers of blocks this approach is
clearly not feasible and we had to use various heuristics to cut
down the number of combinatorial structures considered.
The stacks of Figures~\ref{fig:2-10}, \ref{fig:11-19},
\ref{fig:20-30}, and~\ref{fig:p40-100} were found using extensive
numerical experimentation. The stacks of
Figures~\ref{fig:2-10}, \ref{fig:11-19}, and \ref{fig:20-30} are
optimal, while the stacks of Figure~\ref{fig:p40-100} are either
optimal or very close to being so.

\begin{figure}[t]
\centerline{
\includegraphics[height=80mm]{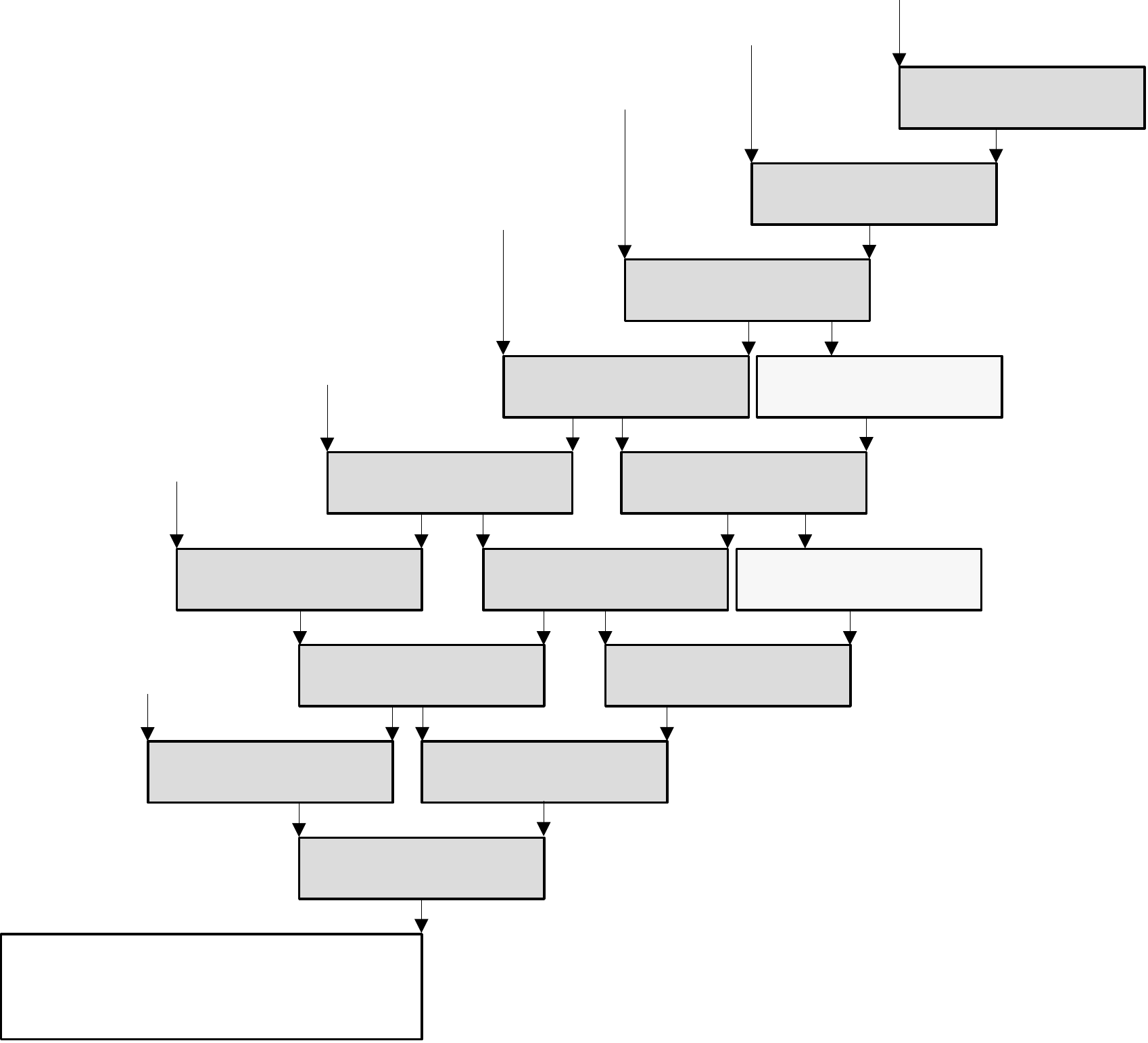}}
\caption{A schematic description of a well-behaved set of
stabilizing forces.} \label{fig:schematic}
\end{figure}

The collections of forces that stabilize the loaded stacks of
Figure~\ref{fig:p40-100} (and the loaded stacks contained in the
stacks of Figures~\ref{fig:2-10}, \ref{fig:11-19}, and
\ref{fig:20-30}) have the following interesting properties. First,
the stabilizing collections of forces of these stacks are unique.
Second, almost all downward forces in these collections are applied
at the edges of blocks. The only exceptions occur when a downward force
is applied on a \emph{right-protruding} block, i.e., a rightmost block in
a level that protrudes beyond the rightmost block of the level
above it. In addition, all point weights are placed on the
left-hand edges of \emph{left-protruding} blocks, where
left-protruding blocks are defined in an analogous way. The table, of
course, supports the (only) block that rests on it at its right-hand
edge. A collection of stabilizing forces that satisfies these
conditions is said to be \emph{well-behaved}. A schematic
description of well-behaved collections of stabilizing forces is
given in Figure~\ref{fig:schematic}. The two right-protruding blocks
are shown with a slightly lighter shading. A right-protruding block
is always adjacent to the block on its left.
We conjecture that forces that balance optimal loaded stacks are
always well-behaved.

A useful property of well-behaved collections of stabilizing forces
is that the total weight of the stack and the positions of its
blocks uniquely determine all the forces in the collection. This
follows from the fact that each block has either two downward forces
acting upon it at specified positions, namely at its two edges, or
just a single force in an unspecified position. Given the upward
forces acting on a block, the downward force or forces acting upon
it can be obtained by solving the force and moment equations of the
block. All the forces in the collection can therefore be determined
in a bottom-up fashion.
We conducted most of our experiments, on blocks with more than 30
blocks, on loaded stacks balanced by well-behaved sets of balancing
forces.

We saw in Section~\ref{sec:prelim} that loaded stacks of total
weight~3, 5, and~7 achieve a larger overhang than the corresponding
unloaded stacks, simply because the number of blocks available for
use in their balancing sets is smaller than the number of point
weights to be applied. The loaded stacks of
Figure~\ref{fig:p40-100} exhibit another trivial impediment
to the conversion of loaded stacks into standard ones: the point weight
to be applied in the lowest position has magnitude less than~$1$. Thus,
these stacks can be converted into standard ones only after making
some small adjustments. These adjustments have only a very small
effect on the overhang achieved. Thus, although we believe that the
difference between the maximum overhangs achieved by loaded and
unloaded stacks is bounded by a small universal constant, we also
believe that for most sizes, loaded stacks yield slightly larger
overhangs.

Although the placements of the blocks in the optimal, or close to
optimal, stacks of Figure~\ref{fig:p40-100} are somewhat irregular,
with some small (essential) gaps between blocks of the same layer,
at a high level, these stacks seem to resemble brick-wall stacks, as
defined in Section~\ref{sec:parabolic}. This, and the fact that
brick-wall stacks were used to obtain the $\Omega(n^{1/3})$ lower
bound on the maximum overhang, indicate that it might be interesting
to investigate the maximum overhang that can be achieved using
brick-wall stacks.

\begin{figure}[t]
\centerline{
\includegraphics[height=80mm]{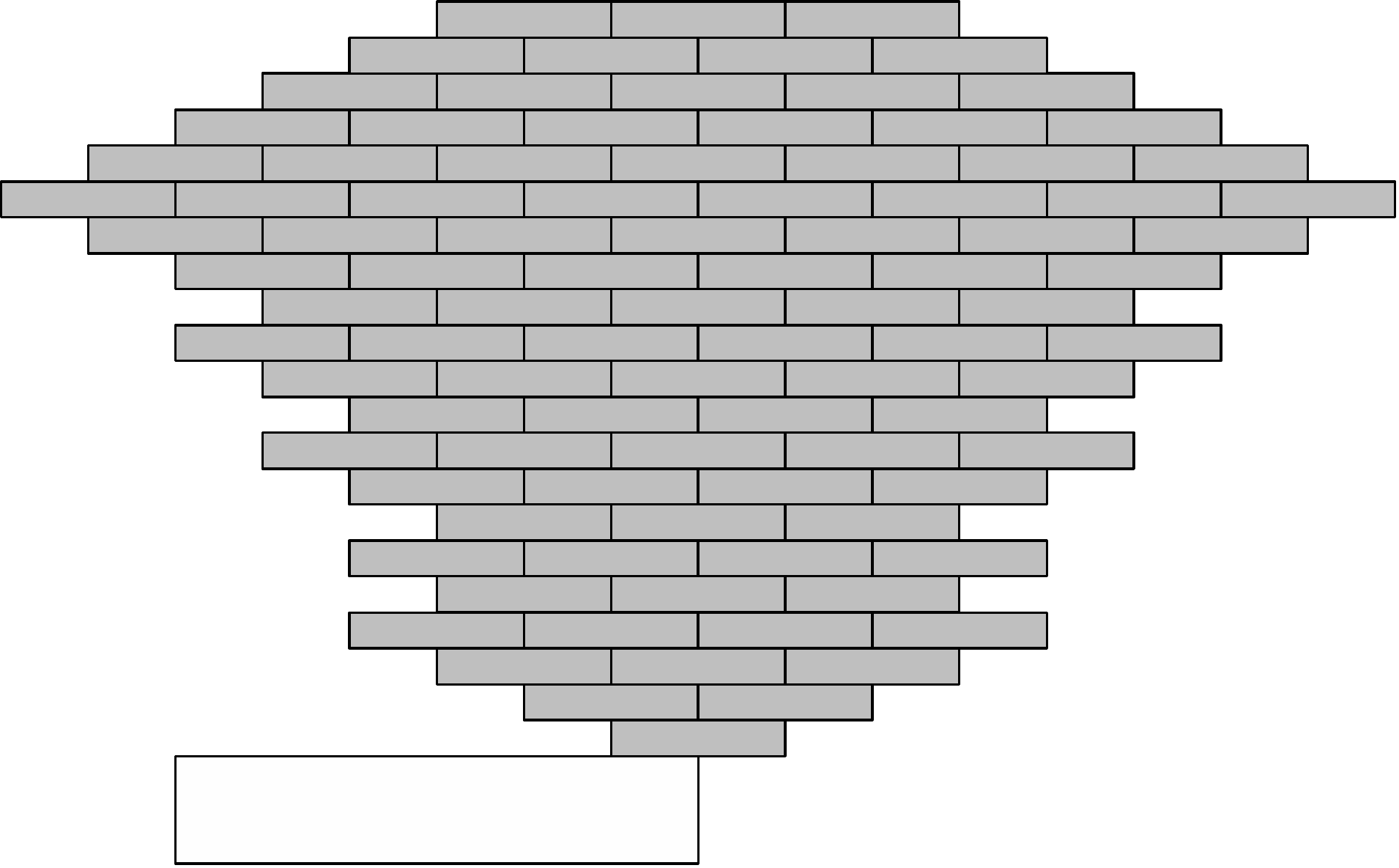}}
\caption{A 95-block symmetric brick-wall stack with overhang 4.}
\label{fig:vase95}
\end{figure}

The parabolic brick-wall stacks of Section~\ref{sec:parabolic} were
designed to enable a simple inductive proof of their balance.
Parabolic stacks, however, are far from being optimal brick-wall
stacks. The balanced 95-block symmetric brick-wall stack with an
overhang of 4 depicted in Figure~\ref{fig:vase95}, for example,
contains fewer blocks and achieves a larger overhang than that
achieved by the 111-block overhang-3 parabolic stack of
Figure~\ref{fig:6stack}.

\begin{figure}[t]
\centerline{
\includegraphics[height=80mm]{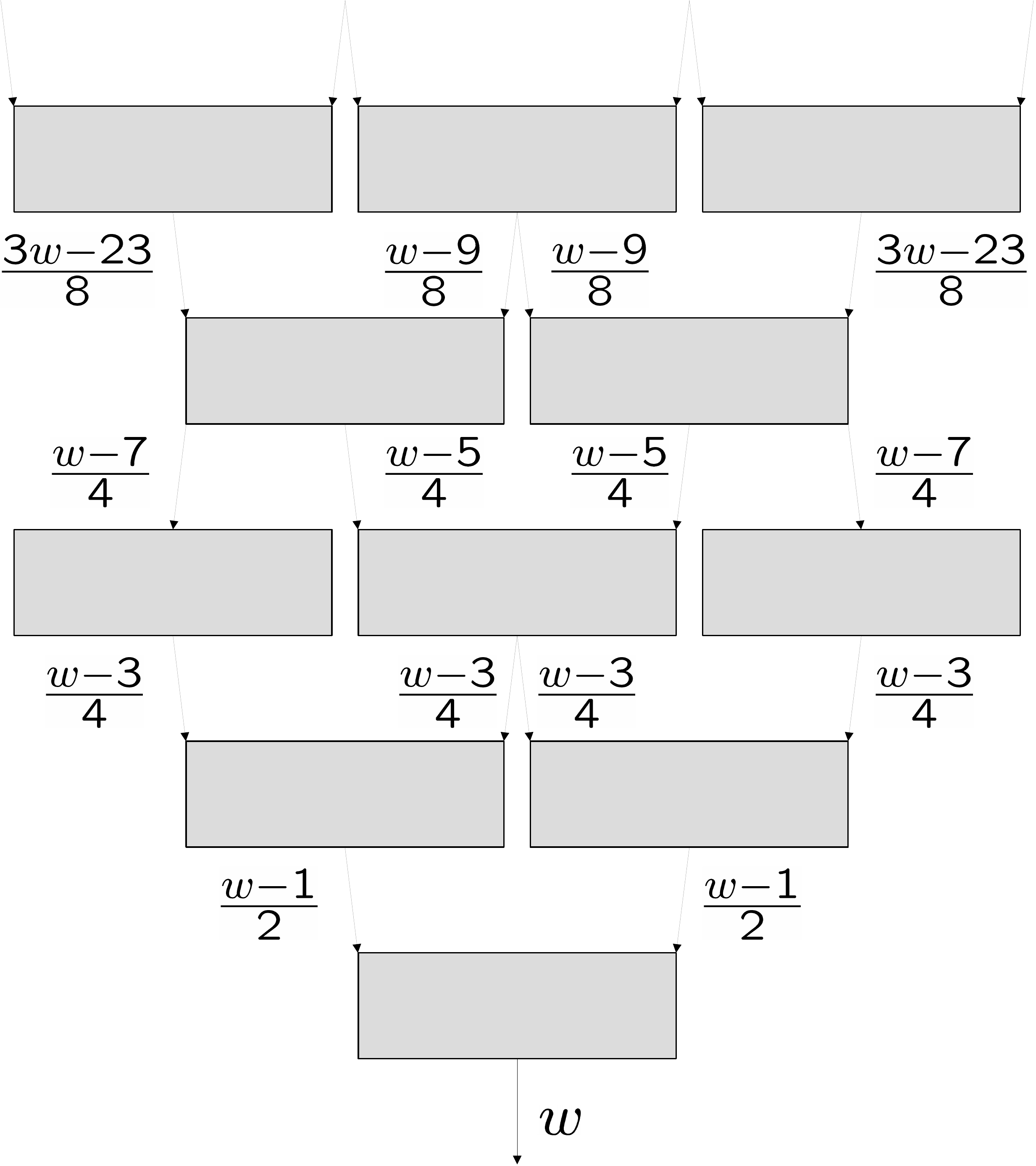}\hspace*{1cm}
\includegraphics[height=80mm]{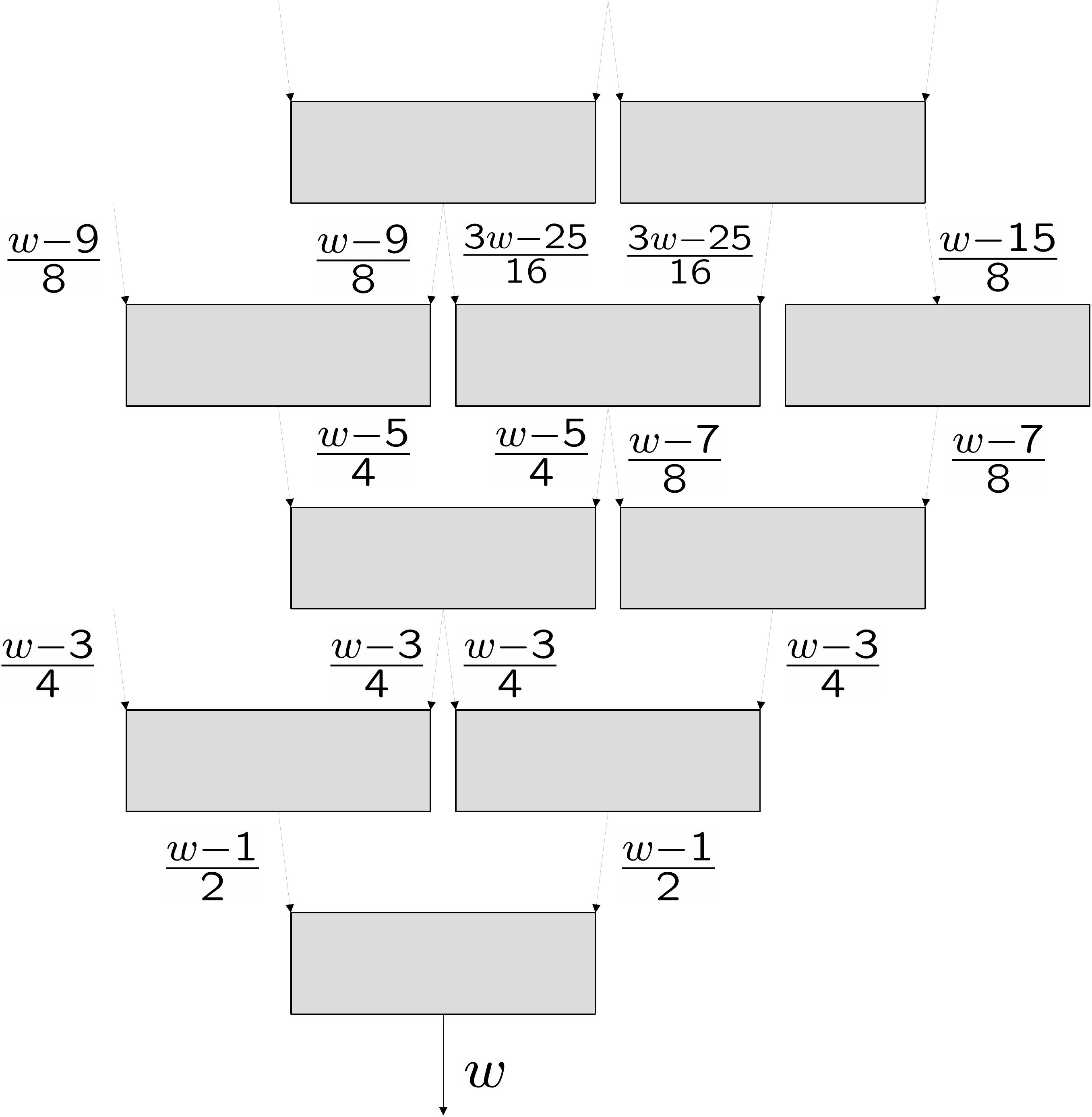}}
\caption{A schematic description of well-behaved collections of forces
that stabilize symmetric and asymmetric brick-wall stacks.}
\label{fig:symfor}
\end{figure}

Loaded brick-wall stacks are especially easy to experiment with.
Empirically, we have again discovered that the minimum weight
collections of forces that balance them turn out to be well-behaved,
in the formal sense defined above. When the brick-wall stacks are
\emph{symmetric} with respect to the $x=0$ axis, and have a flat
top, point weights are attached only to blocks at the top layer of
the stack. Protruding blocks, both on the left and on the right,
then simply serve as \emph{props}, while all other blocks are
perfect \emph{splitters}, i.e., they are supported at the center of
their lower edge and they support other blocks at the two ends of
their upper edge. In non-symmetric brick-wall stacks it is usually
profitable to use the left-protruding blocks as splitters and not as
props, attaching point weights to their left ends. A schematic
description of well-behaved forces that stabilize symmetric and
asymmetric brick-wall stacks is shown in Figure~\ref{fig:symfor}.
As can be seen, all forces in such well-behaved collections are
linear functions of~$w$, the total weight of the stack. This allows
us, in particular, to find the minimum total weight needed to
stabilize a brick-wall loaded stack \emph{without} solving a linear
program. We simply choose the smallest total weight~$w$ for which
all forces are nonnegative. This observation enabled us to
experiment with huge symmetric and asymmetric brick-wall stacks.

\begin{figure}[t]
\centerline{
\includegraphics[height=80mm]{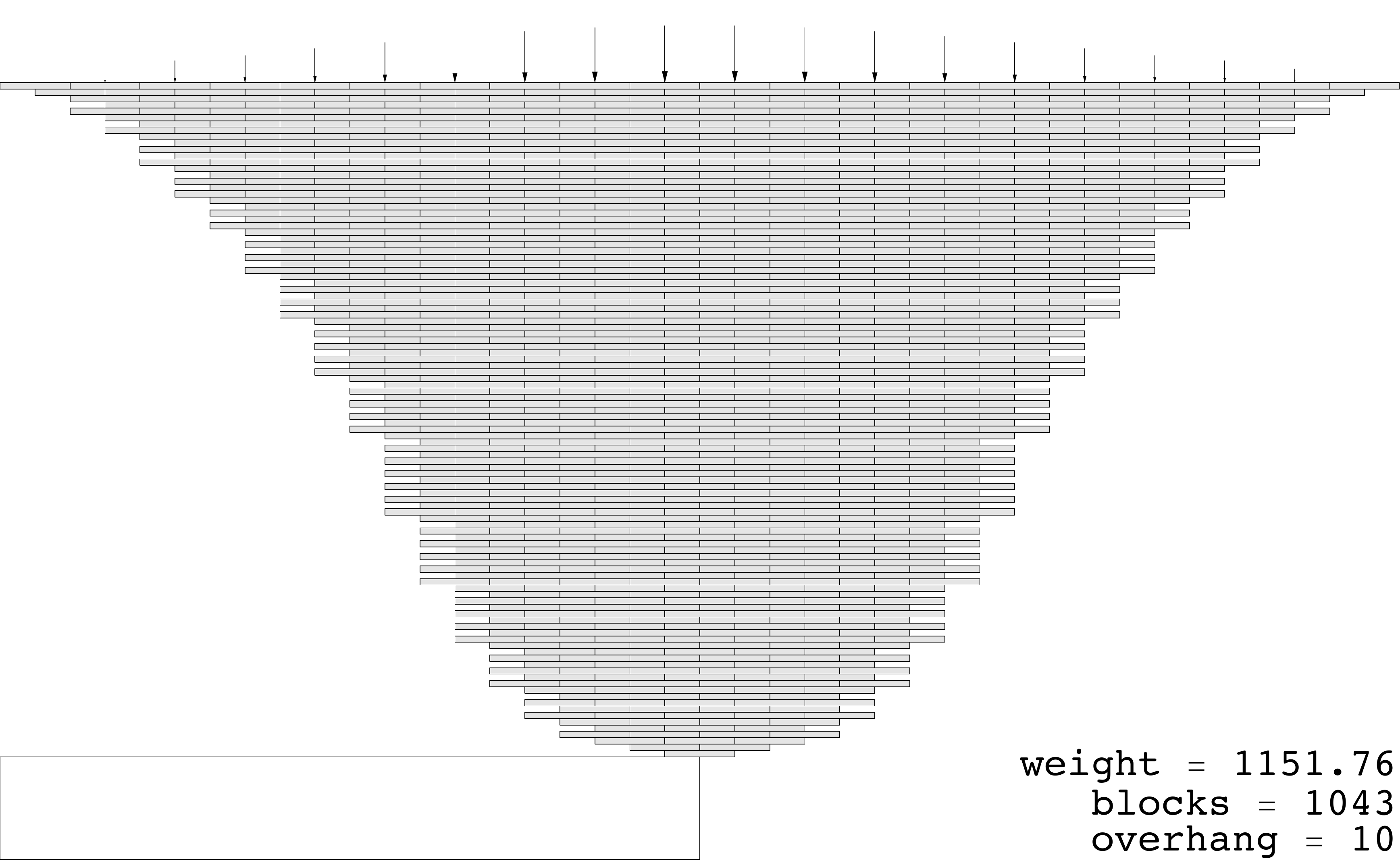}}
\caption{A symmetric loaded brick-wall stack with an overhang of
10.}
\label{fig:s10}
\end{figure}

\begin{figure}[t]
\centerline{
\includegraphics[height=8cm]{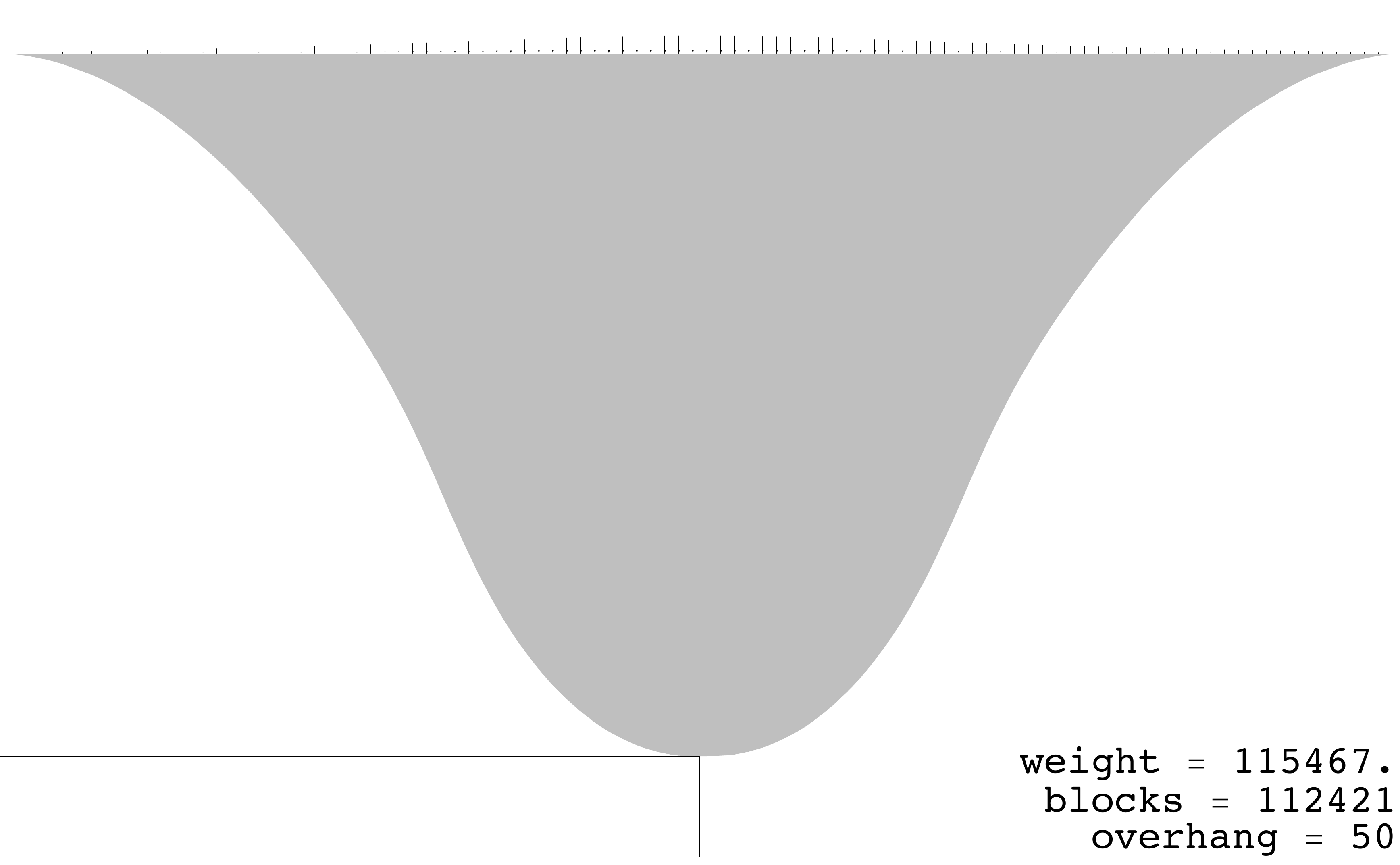}}
\caption{A scaled outline of a loaded brick-wall stack with an
overhang of 50.}
\label{fig:width100}
\end{figure}

The best symmetric loaded brick-wall stacks with overhangs 10 and
50 that we have found are shown in
Figures~\ref{fig:s10} and~\ref{fig:width100}. Their total weights
are about 1151.76 and 115,467, respectively. The blocks in the
larger stack are so small that they are not shown individually. We
again believe that these stacks are close to being the optimal
stacks of their kind. They were found using a local search approach.
In particular, these stacks cannot be improved by widening or
narrowing layers, or by adding or removing single layers.
Essentially the same symmetric stacks were obtained by starting from
almost any initial stack and repeatedly improving it by widening,
narrowing, adding, and removing layers.

As can be seen from Figures~\ref{fig:s10}
and~\ref{fig:width100}, the shapes of optimal symmetric loaded
stacks, after suitable scaling, seem to tend to a limiting curve.
This curve, which we have termed the \emph{vase}, is similar to but
different from that of an inverted normal distribution. We have as yet
no conjecture for its equation.

\begin{figure}[t]
\centerline{
\includegraphics[height=80mm]{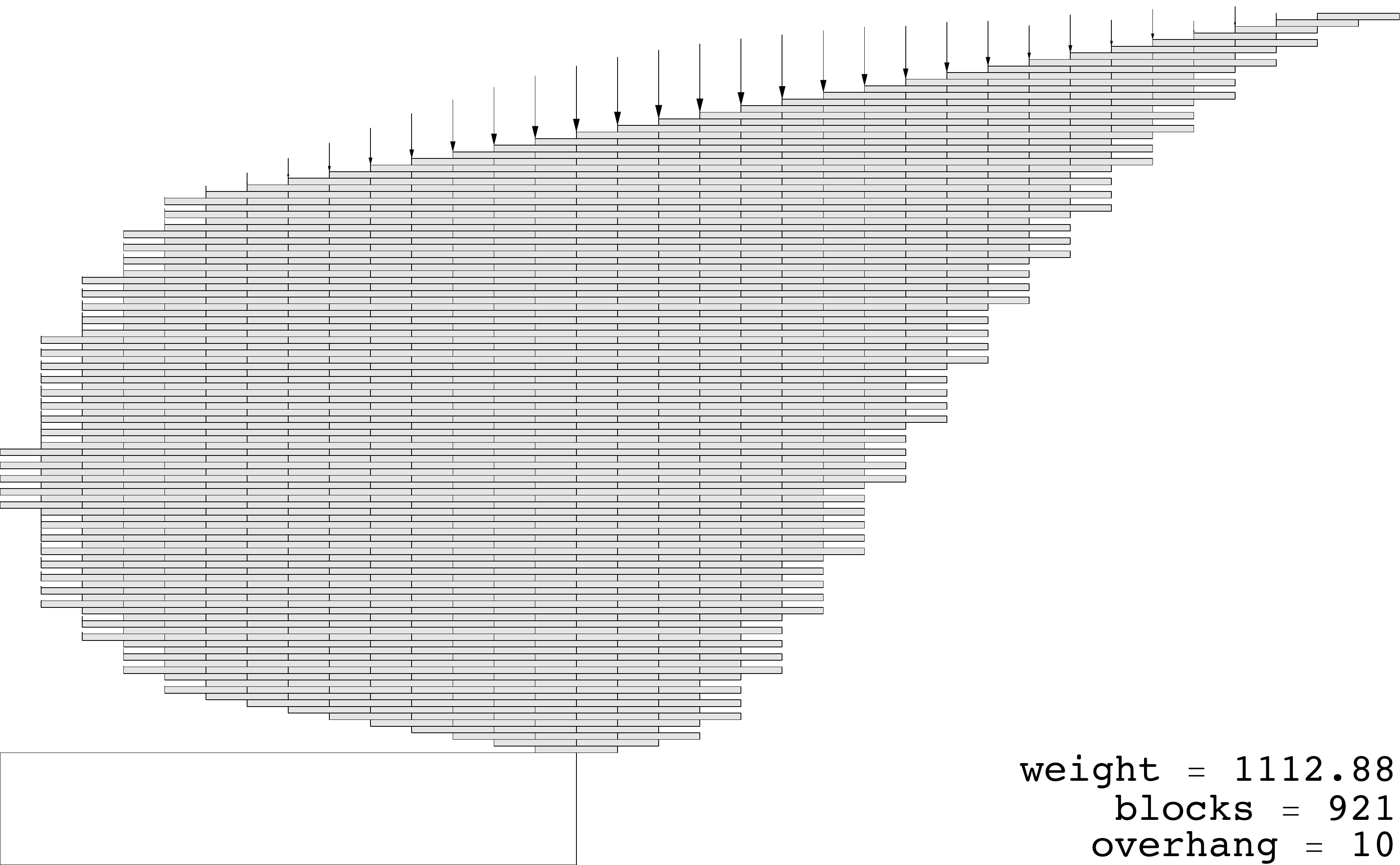}}
\caption{An asymmetric loaded brick-wall stack with an overhang of
10.}
\label{fig:a10}
\end{figure}

We have conducted similar experiments with asymmetric loaded
brick-wall stacks. The best such stack with overhang~10 that we have
found is shown in Figure~\ref{fig:a10}. Its total weight of about
1128.84 is about $3.38\%$ less than the weight of the symmetric
stack of Figure~\ref{fig:s10}. The scaled shapes of optimal
asymmetric loaded brick-wall stacks seem again to tend to a limiting
curve which we have termed the \emph{oil lamp}. We again have no
conjecture for its equation.

\section{Open problems} \label{sec:conc}

\ignore{
We have revisited the classic overhang problem and answered some of
the questions that were latent there. We have shown that the
overhang achievable with~$n$ blocks is exponentially larger than was
previously supposed.
}

Some intriguing problems still remain open. In a subsequent paper
\cite{PPTWZ07}, we show that the $\Omega(n^{1/3})$ overhang lower
bound presented here is optimal, up to a constant factor, but it
would be interesting to determine the largest constant $c_{over}$
for which overhangs of $(c_{over}-o(1))n^{1/3}$ are possible. Can
this constant $c_{over}$ be achieved using stacks that are simple to
describe, e.g., brick-wall stacks, or simple modifications of them,
such as brick-wall stacks with adjacent levels having a displacement
other than~$\frac12$, or small gaps left between the blocks of the
same level?

What are the limiting vase and oil lamp curves? Do they yield,
asymptotically, the maximum overhangs achievable using symmetric and
general stacks?

Another open problem is the relation between the maximum overhangs
achievable using loaded and unloaded stacks. We believe, as
expressed in Conjecture~\ref{con:D}, that the difference between
these two quantities tends to~$0$ as the size of the stacks tends to
infinity. We also conjecture that $D^*(n)-D(n)\leqslant
D^*(3)-D(3)=\frac{5-2\sqrt{6}}{6}\simeq 0.017$ , for every $n\geqslant 1$.

Our notion of balance, as defined formally in
Section~\ref{sec:prelim}, allows stacks to be precarious:
stacks that achieve maximum overhang are always on the verge of
collapse. It is not difficult, however, to define more robust
notions of balance, where there is some \emph{stability}.
In one such natural definition, a stack is \emph{stable} if there is a
balancing set of forces in which none of the forces
acts at the edge of any block.

We note in passing that Farkas' lemma, or the theory of linear
programming duality (see \cite{S98}), can be used to derive an
equivalent definition of stability: a stack is stable if and only if
every feasible infinitesimal motion of the blocks of the stack
increases the total potential energy of the system.

This requirement of stability raises some technical difficulties
but does not substantially change the
nature of the overhang problem. Our parabolic $d$-stacks, for example,
can be made stable by adding a $(d-1)$-row symmetrically placed on top.
The proof of this is straightforward but not trivial. We believe that for
any $n\neq 3$, the loss in the overhang due to this stricter definition is
infinitesimal.

Our analysis of the overhang problem was made under the \emph{no
friction} assumption. All the forces considered were therefore
vertical. The presence of friction introduces horizontal forces and
thus changes the picture completely, as also observed by Hall
\cite{H05}.
We can show that there is a fixed coefficient of friction such that
the inverted triangles are all balanced, and so achieve overhang of
order $n^{1/2}$.

\bibliographystyle{plain}

\end{document}